\documentclass[11pt,reqno]{amsart} 
\usepackage{amsmath,tikz}
    \usetikzlibrary{calc}
\usetikzlibrary{decorations.pathreplacing}
\usetikzlibrary{shapes.geometric}
\usepackage[utf8]{inputenc}
\usepackage{hyperref}
\usepackage[normalem]{ulem}
\usepackage{amsmath}
\usepackage{amssymb}
\usepackage{amsthm}
\usepackage{mathtools}
\usepackage{verbatim}
\usepackage{xcolor}
\usepackage{fullpage}
\usepackage[skip=3pt, indent=20pt]{parskip}
\usepackage{caption}
\usepackage{comment}
\usepackage{marginnote}
\usepackage{enumitem}
\usepackage{caption}
\usepackage[
backend=biber,
style=alphabetic,
sorting=nty
]{biblatex}
\theoremstyle{plain}
\newtheorem{thm}{Theorem}[section]
\newtheorem{lem}[thm]{Lemma}

\newtheorem{cor}[thm]{Corollary}

\theoremstyle{definition}
\newtheorem{defn}[thm]{Definition}

\newtheorem{ex}[thm]{Example}

\theoremstyle{remark}
\newtheorem{rem}[thm]{Remark}

\newcommand{\ZZ}{\mathbb{Z}}
\newcommand{\NN}{\mathbb{N}}

\newcommand{\Sym}{\mathsf{Sym}}

\title{Inflection sets of permutations}
\author{Gordon Rojas Kirby}
\address[G.~Rojas Kirby]{Department of Mathematics and Statistics, San Diego State University, CA 92182}
\email{\textcolor{blue}{\href{mailto:gkirby@sdsu.edu}{gkirby@sdsu.edu}}}
\thanks{}
\author{Mariana Smit Vega Garcia}
\address[M.~Smit Vega Garcia]{Department of Mathematics, Western Washington University, Bellingham, WA 98225}
\email{\textcolor{blue}{\href{mailto:smitvem@wwu.edu}{smitvem@wwu.edu}}}
\thanks{M.S.V.G \ was partially supported by the NSF grant DMS-2348739 and a Karen EDGE fellowship}

\begin{document}

\begin{abstract} 
In 2020, Brualdi and Dahl introduced the discrete derivative $\sigma'(i)=\sigma(i+1)-\sigma(i)$ of a permutation $\sigma\in\Sym_n$. In this paper, we focus on the second derivative of a permutation, $\sigma''(i)=\sigma(i+1)+\sigma(i-1)-2\sigma(i)$, and the inflections it determines. We show that for even $n$, the only distinct permutations with the same second derivative are the identity and its reverse, while for odd $n$, the second derivative is not injective even after excluding them. We then study which sets of pairs of integers $(i,j)$ are the inflection set of some permutation and whether a given set is witnessed in $\Sym_n$ for infinitely many $n$. We show that exactly $14$ permutations of $[n]$ have no inflections when $n\ge5$. For a single inflection $\{(k,m)\}$, we determine all $n$ for which it is witnessed in $\Sym_n$: we prove it is always witnessed, and it is witnessed for only finitely many $n$ if and only if $m\ge 2k+1$. In that case, all such $n$ satisfy $2m-2k+1\le n\le 2m-k+2$ (apart from $k=2$, $n=m+1$), and, perhaps surprisingly, this range has gaps once $k\ge 10$, the smallest occurring at $(k,m,n)=(10,21,30)$. Finally, we characterize maximally inflected permutations, whose second derivative alternates in sign, and show that their number $|\mathcal M_n|$ satisfies $2E_n\le|\mathcal M_n|\le(2/3)^{\lfloor n/3\rfloor-1}\,n!$, where $E_n$ is the Euler zigzag number.
\end{abstract}
\maketitle


\section{Introduction}

In this paper, we consider derivatives and inflections of permutations in $\Sym_n$.  
Throughout, we let $n\in\NN=\{1,2,3,\ldots\}$, $[n]=\{1,2,\ldots,n\}$, $\pm[n]=[n]\cup\{-1,-2,\dots,-n\}$, and $[i,j]=\{i,i+1,\dots,j\}$ for integers $i\le j$, with $[i,j]=\emptyset$ if $i>j$.
We let $\Sym_n$ denote the set of permutations of $[n]$, and given a permutation $\pi\in\Sym_n$ we write $\pi=\pi(1)\pi(2)\cdots \pi(n)$ in one-line notation. 
Our goal is to understand inflections of permutations. To do so, we define the first and second derivative of a permutation.

\begin{defn}
Let $\pi\in \mathsf{Sym}_n$. We define the \emph{derivative} of $\pi$ to be the function $\pi':[n-1]\rightarrow\mathbb{Z}$ given by $\pi'(j)=\pi(j+1)-\pi(j)$. We sometimes use $D$ to denote the map $D:\Sym_n\to \ZZ^{n-1}$ defined by $D(\pi)=\pi'$. Given an integer $2\leq\ell\leq n-1$, we define the \textit{second derivative of $\pi$ at $\ell$} to be 
\[
\pi''(\ell)\coloneqq\pi(\ell+1)+\pi(\ell-1)-2\pi(\ell).
\]
\end{defn}
With a notion of second derivatives at hand, we can define convexity/concavity.

\begin{defn}
We say $\pi\in \Sym_n$ is \textit{convex (concave) at}  $m\in[2,n-1]$ if  $\pi(m+1)+\pi(m-1)-2\pi(m)>0$ (respectively, $\pi(m+1)+\pi(m-1)-2\pi(m)<0$).
\end{defn}

\begin{defn} Let $n\ge 4$, $\pi\in \mathsf{Sym}_n$, and $m,p\in[2,n-1]$ with $m<p$. We say that $\pi$ has an \textit{inflection} between $m$ and $p$ if 
\begin{itemize}
    \item $\pi''(m)\pi''(p)<0$ and 
    \item $\pi''(\ell)=0$ for all $m<\ell<p$.   
\end{itemize} 
We denote the inflection set of a permutation $\sigma$ by $\mathcal{I}(\sigma)$, consisting of all pairs of integers $(i,j)$, $2\leq i<j\leq n-1$ such that $\sigma$ has an inflection between $i$ and $j$.
\end{defn}

\begin{ex}
Consider the permutation $\pi=23468715\in\Sym_8$. We have
\begin{center}    
\begin{tabular}{|c|c|c|c|c|c|c|c|c|}\hline
     $\ell$&2&3&4&5&6&7 \\\hline
$\pi''(\ell)$ &0&1&0&-3&-5&10   \\\hline
\end{tabular}
\end{center}
This shows that $\pi$ has inflections between 3 and 5 and between 6 and 7, which we denote as $\mathcal{I}(\pi)=\{(3,5), (6,7)\}$.

    \center
    \begin{tikzpicture}[yscale=.4, xscale=.8]

    \def \n{9}
    \draw[-latex] (0,0) to (0,\n+1);
    \draw[-latex] (0,0) to (\n+1,0);

    \foreach \i in {1,2,...,\n}{
        \draw[] (-0.1, \i) to (0.1,\i) ;
        \draw[] (\i,-0.1) to (\i,0.1);         
    }

    \def\yvalues{2,3,4,6,8,7,1,5}

    \foreach \y [
        count=\x,
        remember=\x as \previousx,
        remember=\y as \previousy
    ] in \yvalues{
        \node at (\x,\y) {\small$\bullet$};
        \ifnum \x>1
            \draw[thick] (\previousx,\previousy) to (\x,\y);
        \fi
    }
    \end{tikzpicture} 
    
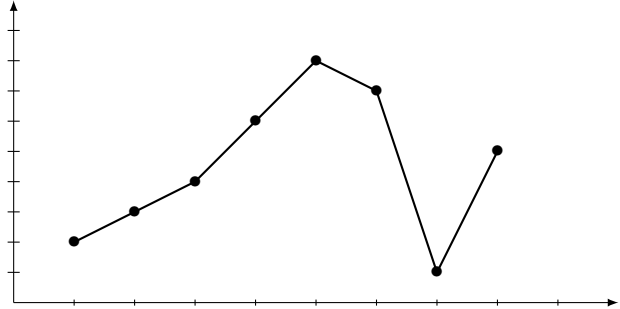
\captionof{figure}{The graph of $\pi=23468715$ with inflection set $\mathcal{I}(\pi)=\{(3,5),(6,7)\}$.}
    \label{fig:InflectionSet(3,5),(6,7)}
    \end{ex}

\begin{defn}\label{defn:fundamental symmetr}
    Fix $w \in \Sym_n$. We call the \emph{reverse} of $w$ the permutation $w^r := w(n) \cdots w(1) \in \Sym_n$, and the \emph{complement} of $w$ the permutation $\overline{w} \in \Sym_n$ defined by $\overline{w}(i) = n+1-w(i)$. The \emph{fundamental symmetries} of $w$ are the (possibly non-distinct) permutations
    $$\{w, w^r, \overline{w}, \overline{w}^r = \overline{w^r}\}.$$
\end{defn}

The first derivative $\pi'$ of a permutation was introduced by Brualdi and Dahl \cite{brualdi-dahl} under the name \emph{discrete derivative}. They show that $D:\Sym_n\rightarrow \mathbb{Z}^{n-1}$ is injective and characterize the integer vectors that arise as first derivatives. We recall these facts in Lemmas~\ref{L:f'injective} and~\ref{lem:firstderv}. They also classify \emph{convex} permutations, which are those for which $\pi'$ is non-decreasing. Since a permutation has no inflections exactly when it or its complement is convex, their classification is equivalent to Theorem~\ref{thm:no inflection points}, for which we give an independent proof. Using this result, we conclude that if $n\ge 5$, the number of permutations with no inflections is $14$.

The second derivative $\pi''$ and the notion of inflection are not considered in \cite{brualdi-dahl}. In particular, the question of when two permutations share a second derivative, and the study of permutations with a prescribed inflection set, are new.

Our approach to inflection sets is modeled on the study of other position and value-based statistics. Billey, Burdzy and Sagan \cite{billey-burdzy-sagan} enumerate permutations with a given peak set, and Davis, Nelson, Petersen and Tenner \cite{davis-nelson-petersen-tenner} introduce the pinnacle set, characterize the admissible pinnacle sets, and study how many permutations realize each one. Analogous questions have been studied by the authors, together with Gonz{\'a}lez, Harris, and Tenner, in the context of
 signed permutations of types $B$ and $D$ \cite{PinnaclesTypeBD} and for Stirling permutations, where pinnacles are replaced by mesas \cite{gonzalez2023mesas}. See also \cite{diazlopez-pinnacle}. We ask the same two questions for inflections: which sets occur as the inflection set of some $\sigma\in\Sym_n$, and how many witnesses each has. A feature with no counterpart for peaks or pinnacles is that admissibility of an inflection set depends on $n$; a set may be realized in $\Sym_n$ for some values of $n$ and for no others.

Maximally inflected permutations generalize alternating permutations (see \cite{stanley2010survey} for a survey), which are exactly the permutations whose first derivative alternates in sign. Here, we ask instead that the second derivative alternates in sign. In Theorem~\ref{T:odd}, we obtain a necessary condition on run lengths,  which places maximally inflected permutations within the family of permutations with restricted run structure, whose enumeration is studied in \cite{fewster-siemssen,zhuang-runs}. Our bounds on the number of maximally inflected permutations of $\Sym_n$ are in the spirit of the asymptotics of consecutive pattern avoidance \cite{elizalde-noy,ehrenborg-kitaev-perry}, where the number of permutations in $\Sym_n$ avoiding a fixed consecutive pattern grows like $c\,\lambda^n n!$ for constants $c>0$ and $0<\lambda<1$. Asymptotics of the same shape, obtained through the spectrum of an integral operator, also govern permutations with a periodic pattern of rises and falls \cite{bender-helton-richmond}. Omar and Troyka \cite{OmarTroyka2026} study growth rates of permutations with a fixed descent or peak set. However, maximal inflection is not a condition on relative order alone, since the sign of $\sigma''(i)$ depends on the values $\sigma(i-1),\sigma(i),\sigma(i+1)$, so none of these results applies directly.

Finally, $\sigma''(i)=0$ exactly when $\sigma(i-1),\sigma(i),\sigma(i+1)$ are in an arithmetic progression. Permutations with $\sigma''$ nowhere zero are thus those avoiding three consecutive terms in arithmetic progression, counted by \cite{oeis-A295370}. This is different from the classical problem of avoiding three-term arithmetic progressions at arbitrary positions \cite{degs}.

The outline of the paper is as follows.
In Section \ref{S:injectivityfirst}, we recall results from \cite{brualdi-dahl} pertaining to the injectivity of the first derivative. Our first main result analyzes the injectivity of the second derivative in $\Sym_n$ for $n$ even and odd. When $n$ is even, we show that  the only two distinct permutations with the same
second derivative are the identity (denoted by $e$) and its reverse $e^r$. When $n$ is odd, we provide examples of permutations, distinct from $e$ and $e^r$, with the same second derivative.
In Section \ref{S:Admissible-Inflections}, we define admissible inflection sets. Those are  divided into two types, based on whether there are infinitely many $n$ such that there is a witness in $\Sym_n$.

In Section \ref{S:NoInflection}, we provide a different proof of \cite[Theorem~6]{brualdi-dahl}, which characterizes permutations with no inflections. We also prove that the number of permutations of $[n]$ with no inflections is eventually constant for all sufficiently large $n$.

In Section \ref{S:Perms-with-one-inflection}, 
we show that any inflection set of size one is admissible, that is, given $2\le k<m$, $I=\{(k,m)\}$ is an admissible inflection set. Moreover, our second and third main results, Theorems \ref{T:windows} and \ref{T:short}, determine, in terms of $k$, $m$, and $n$, exactly when $I$ has a witness in $\Sym_n$. As a consequence, we show that given $2\le k<m$ and $I=\{(k,m)\}$, $I$ is vanishing if and only if $m\ge 2k+1$.

In Section \ref{S:max}, we study the size of $\mathcal M_n$, the set of maximally inflected permutations in $\Sym_n$. We first  characterize maximally inflected permutations (Theorem~\ref{T:char}) and show that alternating permutations are maximally inflected (Corollary~\ref{C:alt}). Using  classical asymptotics of the Euler zigzag numbers, this gives $\lim_{n\to\infty}|\mathcal M_n|^{1/n}=\infty$. The rest of the section addresses how $|\mathcal M_n|$ compares to $n!$. We give two upper bounds on the fraction $|\mathcal M_n|/n!$, both probabilistic. We first prove a necessary condition for a permutation to be maximally inflected: every window of $k$ consecutive positions must standardize to a permutation whose interior runs all have odd length (Theorems~\ref{T:odd} and \ref{T:firstboundsMgood}). Via Fekete's Lemma, we are able to bound $\limsup_n(|\mathcal M_n|/n!)^{1/n}$ by $\lim_n p_n^{1/n}$, where $p_n$ is the proportion of permutations in $\Sym_n$ all of whose interior runs have odd length. Our last main result, Theorem~\ref{T:estimatemax}, 
samples a uniformly random $\sigma$ and shows that each block of three consecutive values, conditioned on the past, is consistent with maximal inflection with probability at most $2/3$. Iterating over the $\lfloor n/3\rfloor$ blocks gives the explicit bound $|\mathcal M_n|/n!\le(2/3)^{\lfloor n/3\rfloor-1}$, hence $|\mathcal M_n|/n!\to 0$ with $\limsup_n(|\mathcal M_n|/n!)^{1/n}\le(2/3)^{1/3}$.

\subsection*{Acknowledgements}

The authors are very thankful to Nicolle Gonz\'{a}lez, Pamela E. Harris and Bridget Eileen Tenner for helpful discussions.

\section{Injectivity of Derivatives}\label{S:injectivityfirst}
In this section, we discuss the injectivity of the first and second derivatives in $\Sym_n$.

\subsection{Injectivity of the first derivative}
The injectivity of the first derivative has been settled. 
\begin{lem}[{\cite[Proposition~2]{brualdi-dahl}}]\label{L:f'injective}
    The map $D:\Sym_n\to \ZZ^{n-1}$ defined $D(\pi)=\pi'$ is injective.
\end{lem}

The following result characterizes which $(n-1)$-tuples in $(\pm [n-1])^{n-1}$ are derivatives of a permutation.

\begin{lem}[{\cite[Proposition~4]{brualdi-dahl}}]\label{lem:firstderv}
    Given $\beta\in(\pm[n-1])^{n-1}$, there exists $\alpha\in \Sym_n$ such that $\alpha'=\beta$ if and only if $0<|\sum_{i=j}^k\beta_i|\le n-1$ for all $j,k\in[n-1]$ with $j\le k$.
\end{lem}

\begin{ex} Lemma \ref{lem:firstderv} yields the following $12$ first derivatives of permutations $\pi\in\Sym_4$ with $\pi'(1)>0$ (the other $12$ are their negatives). We use $\overline{i}$ to denote $-i$.
    \[
    111 \quad 3\bar1\bar1  \quad 3\bar21 \quad 21\bar2 \quad 2\bar12\quad 2\bar1\bar2 \quad 2\bar32 \quad 12\bar1 \quad 11\bar3 \quad 1\bar23 \quad 1\bar2\bar1 \quad 1\bar3 1
    \]
\end{ex}
\subsection{(Lack of) Injectivity of the second derivative}
\begin{lem}
    Let $\sigma,\pi\in \Sym_{2k+1}$ for some nonnegative integer $k$. If $\sigma''=\pi''$, then $\sigma(k+1)=\pi(k+1)$.
\end{lem}
\begin{proof}
    Notice that 
     \begin{align*}
    \binom{2}{2}\sigma''(2)+\cdots+\binom{k}{2}\sigma''(k)+\binom{k+1}{2}\sigma''(k+1)+\binom{k}{2}\sigma''(k+2)+\cdots+\binom{2}{2}\sigma''(2k)\\
    =\sigma(1)+\cdots+\sigma(k)-2k\sigma(k+1)+\sigma(k+2)+\cdots+\sigma(2k+1).
    \end{align*}
   Since $\sigma,\pi\in\Sym_{2k+1}$, 
   \[
   \sigma(1)+\cdots+\sigma(2k+1)=\frac{(2k+1)(2k+2)}{2}=\pi(1)+\cdots+\pi(2k+1).
   \]Therefore if $\sigma''=\pi''$, then
   \[
   \frac{(2k+1)(2k+2)}{2}-(2k+1)\sigma(k+1)=\frac{(2k+1)(2k+2)}{2}-(2k+1)\pi(k+1).
   \]
   Consequently, $\sigma(k+1)=\pi(k+1)$.
\end{proof}
More generally, we can extract similar information by looking at a weighted sum of second derivatives:
\begin{lem}\label{lem:weighted-sum-second-derivative} Let $\sigma,\pi\in \Sym_n$. If $\sigma''=\pi''$, then 
for each $1\leq k\leq n-1$,
\[
    A(k)\sigma(k)+B(k)\sigma(k+1)=A(k)\pi(k)+B(k)\pi(k+1),
\]
where $A(k)=\frac{n(n-1-2k)}{2}$ and $B(k)=-\frac{n(n+1-2k)}{2}$.
\end{lem}
\begin{proof}
    Given $1\leq k\leq n-1$, we compute directly
    \[
    \sum_{i=2}^k \binom{i}{2}\sigma''(i)+\sum_{i=k+1}^{n-1}\binom{n-i+1}{2}\sigma''(i)=A(k)\sigma(k)+B(k)\sigma(k+1)+\sum_{i=1}^n\sigma(i),
   \]
    where \[
    A(k)=\binom{k-1}{2}-2\binom{k}{2}+\binom{n-k}{2}-1=\frac{n(n-1-2k)}{2},
    \]  
    and
   \[
    B(k)=\binom{k}{2}-2\binom{n-k}{2}+\binom{n-k-1}{2}-1=-\frac{n(n+1-2k)}{2}.
    \]
       Since \[
\sum_{i=1}^n\sigma(i)=\sum_{i=1}^n \pi(i)=\frac{(n+1)n}{2},\] this yields \[A(k)\sigma(k)+B(k)\sigma(k+1)=A(k)\pi(k)+B(k)\pi(k+1)\] for $1\leq k\leq n-1$. 
\end{proof}

\begin{rem} We notice that
\[
\gcd(A(k),B(k))=\begin{cases} n/2, \text{ if } n \text{ is even}\\
n, \text{ if } n \text{ is odd.}
\end{cases}
\]
Indeed, we have $A(k)=\frac{n(n-2k-1)}{2}$ and $B(k)=-\frac{n(n-2k+1)}{2}$. Since $\gcd(a,b)=\gcd(a,b-a)$, $\gcd\left(\frac{n(n-2k-1)}{2},\frac{n(n-2k+1)}{2}\right)=\gcd\left(\frac{n(n-2k-1)}{2}, n\right)$.
If $n$ is odd, then $n-2k-1$ is even, and $n$ divides $A(k)$, giving $\gcd(A(k),B(k))=n$.

If $n$ is even, then $n-2k-1$ is odd. Writing $n=2\ell$, we have $\frac{n(n-2k-1)}{2}=\ell(n-2k-1)$ and $\text{gcd}(A(k),B(k))=\text{gcd}(\ell(n-2k-1),2\ell)=\ell\text{gcd}(n-2k-1,2)=\ell=n/2$.
    
\end{rem}

\begin{ex}\label{ex:n7}
Let $n=7$ and let $\sigma,\pi\in\Sym_7$ with $\sigma''=\pi''$. Here $A(k)=7(3-k)$ and $B(k)=-7(4-k)$, so
\[
\begin{array}{c|cccccc}
k & 1 & 2 & 3 & 4 & 5 & 6\\ \hline
A(k) & 14 & 7 & 0 & -7 & -14 & -21\\
B(k) & -21 & -14 & -7 & 0 & 7 & 14
\end{array}
\]
Dividing each equation of Lemma~\ref{lem:weighted-sum-second-derivative} by $\gcd(A(k),B(k))=7$ gives
\begin{align*}
2\sigma(1)-3\sigma(2)&=2\pi(1)-3\pi(2) && (k=1)\\
\sigma(2)-2\sigma(3)&=\pi(2)-2\pi(3) && (k=2)\\
\sigma(4)&=\pi(4) && (k=3 \text{ and } k=4)\\
\sigma(6)-2\sigma(5)&=\pi(6)-2\pi(5) && (k=5)\\
2\sigma(7)-3\sigma(6)&=2\pi(7)-3\pi(6) && (k=6).
\end{align*}
Writing $\delta(i)=\sigma(i)-\pi(i)$, this system reads
\[
2\delta(1)=3\delta(2),\qquad \delta(2)=2\delta(3),\qquad \delta(4)=0,\qquad \delta(6)=2\delta(5),\qquad 2\delta(7)=3\delta(6).
\]
 The solutions are
\[
(\delta(1),\ldots,\delta(7))=(3t,\,2t,\,t,\,0,\,s,\,2s,\,3s)
\]
for integers $s,t$. Since $\sum_{i=1}^7\delta(i)=0$, we get $s=-t$, that is, $\delta(i)=t(4-i)$. Finally, $|\delta(1)|\le 6$ forces $|t|\le 2$. The case $t=0$ corresponds to $\sigma=\pi$, and $t=\pm2$ gives $\{\sigma,\pi\}=\{e,e^r\}$. The case $t=\pm 1$ occurs for other pairs of permutations: for instance, $\sigma=4567123$ and $\pi=1357246$ satisfy $\sigma''=\pi''=(0,0,-7,7,0)$ and $\sigma-\pi=(3,2,1,0,-1,-2,-3)$.
\end{ex}

We start the analysis of the injectivity of the second derivative in $\Sym_n$ for $n$ even by noticing that in that case, if two permutations have the same second derivative and coincide at one value, they must be the same:

\begin{lem}\label{L:equalifeven}
    Let $\sigma,\pi\in \Sym_n$ be such that $\sigma''=\pi''$. If $n$ even and there exists $i\in[n]$ with $\sigma(i)=\pi(i)$, then $\sigma=\pi$.
\end{lem}

\begin{proof} By Lemma \ref{lem:weighted-sum-second-derivative}, for every $1\le k\le n-1$,
\[
A(k)(\sigma(k)-\pi(k))+B(k)(\sigma(k+1)-\pi(k+1))=0,
\]
  where $A(k)=\frac{n(n-1-2k)}{2}$ and $B(k)=-\frac{n(n+1-2k)}{2}$. Also notice that $A(k)=0$ if and only if $k=\frac{n-1}{2}$, and $B(k)=0$ if and only if $k=\frac{n+1}{2}$, neither of which happens when $n$ is even. Consequently, if $n$ is even, we can write for $1\le k\le n-1$,
  \[
  \sigma(k+1)-\pi(k+1)=-\frac{A(k)}{B(k)}(\sigma(k)-\pi(k)).
  \]
  Therefore $\sigma(j)-\pi(j)=0$ if and only if $\sigma(j+1)-\pi(j+1)=0$ for every $1\le j\le n-1$. This gives a chain of equalities showing that if $\sigma(i)=\pi(i)$ for some $i\in[n]$, then $\sigma(j)=\pi(j)$ for all $j\in[n]$.
\end{proof}

The following theorem proves when $n$ is even, the only two distinct permutations with the same second derivative are the identity and its reverse:

\begin{thm}\label{thm:injectivity-of-second-derivative}
    Let $\sigma,\pi$ be distinct permutations of $\Sym_n$, where $n$ is even. If $\sigma''=\pi''$, then $\sigma=\pi^r$ and $\sigma$ or $\pi$ is the identity. 
\end{thm}
\begin{proof}
    By Lemma \ref{lem:weighted-sum-second-derivative} we know that $
    (n-3)\sigma(1)-(n-1)\sigma(2)=(n-3)\pi(1)-(n-1)\pi(2)$. Rearranging terms yields
    \[\frac{n-3}{n-1}(\sigma(1)-\pi(1))=\sigma(2)-\pi(2)\]
    Since $n$ is even, $\gcd(n-3,n-1)=1$. Hence, $(n-1)\mid (\sigma(1)-\pi(1))$, where $-n<\sigma(1)-\pi(1)<n$. If $\sigma(1)-\pi(1)=0$, then by Lemma \ref{L:equalifeven} we would obtain $\sigma=\pi$, contradicting the assumption that they are distinct. Hence, $|\sigma(1)-\pi(1)|=n-1$ so that $\sigma(1)=1, \pi(1)=n$ or $\sigma(1)=n,\pi(1)=1$. Assume without loss of generality that $\sigma(1)=1$ and $\pi(1)=n$.  

    We now use strong induction to finish the proof.  Let $2\le i\le n$ and assume $\sigma(j)=j$ and $\pi(j)=n+1-j$ for all $1\le j<i$. Then $\sigma(i)\ge i$ and $\pi(i)\le n-i+1$. Therefore
    \begin{equation}\label{sigmapi}
    \sigma(i)-\pi(i)\ge i-(n-i+1)=2i-n-1,
    \end{equation}
    with equality forcing $\sigma(i)=i$ and $\pi(i)=n-i+1$. Lemma \ref{lem:weighted-sum-second-derivative} (with $k=i-1$) gives
    \[
    \frac{n(n-1-2(i-1))}{2}(i-1-(n+1-(i-1)))=-\frac{n(n+1-2(i-1))}{2}(\pi(i)-\sigma(i))
    \]
    that is, 
    \[
    \frac{n(n+1-2i)}{2}(2i-n-3)=-\frac{n(n+3-2i)}{2}(\pi(i)-\sigma(i))
    \]
    hence $    (n+1-2i)=\pi(i)-\sigma(i)$, from which \eqref{sigmapi} gives us $\sigma(i)=i$ and $\pi(i)=n-i+1$, as desired. In this case, $\sigma=e$ and $\pi=e^r$.
\end{proof}
When $n$ is odd, the second derivative is not injective, even when excluding the identity and its reverse: 
\begin{thm}\label{L:notinjodd}
Let $k\ge 1$. The second derivative is not injective on $\Sym_{2k+1}\setminus\{e,e^r\}$.
\end{thm}
\begin{proof}
    Let $n=2k+1$, and take $\sigma=135\cdots n24\cdots (n-1)$ and $\pi=(k+1)(k+2)\cdots n123\cdots k$. Then $\sigma''=0\cdots0(-n)n0\cdots 0=\pi''$.
\end{proof}

\begin{rem}
For odd $n=2m+1$, two distinct permutations $\sigma$, $\pi\in\Sym_n$ with $\sigma''=\pi''$ satisfy $\sigma(i)-\pi(i)=c\,(i-(m+1))$ for some $c\in\{\pm1,\pm2\}$, and $c=\pm2$ only for $\{e,e^r\}$. The pairs with $c=\pm1$ correspond to permutations $\gamma$ of $\{-m,\dots,m\}$ for which $j\mapsto \gamma(j)+j$ is a permutation. These are the additive permutations of Kotzig and Laufer \cite{kotzig-laufer}, counted by \cite{oeis-A002047}. Hence the number of unordered pairs of distinct permutations in $\Sym_{2k+1}$ with the same second derivative is $1+a(k)$, where $a(k)$ denotes the $k$th term of \cite{oeis-A002047} (so $a(1),\dots,a(4)=2,6,28,244$). For $n=3,5,7,9$ this gives $3,7,29,245$ pairs, respectively.
\end{rem}

\section{Admissible Inflection Sets}\label{S:Admissible-Inflections}

We say that if a set of pairs of integers $I$ is such that there exists a permutation $\sigma$ with inflection set $\mathcal{I}(\sigma)=I$, the set is \emph{witnessed} by it. The sets that are witnessed are split into two types, based on whether there are infinitely many $n$ such that there is a witness in $\Sym_n$. In this section, we provide examples to illustrate inadmissibility and each of the two types of admissibility.

\begin{defn}
    We say that a set of pairs of integers $I$ is \textit{admissible} if there exist $n$ and a permutation $w\in \Sym_n$ (called a \emph{witness} for $I$) such that $\mathcal{I}(w)=I$. We say that an admissible inflection set $I$ is \textit{non-vanishing} if there exist infinitely many $n$ so that there is a witness in $\Sym_n$. Otherwise, we say that $I$ is \textit{vanishing}.
\end{defn}

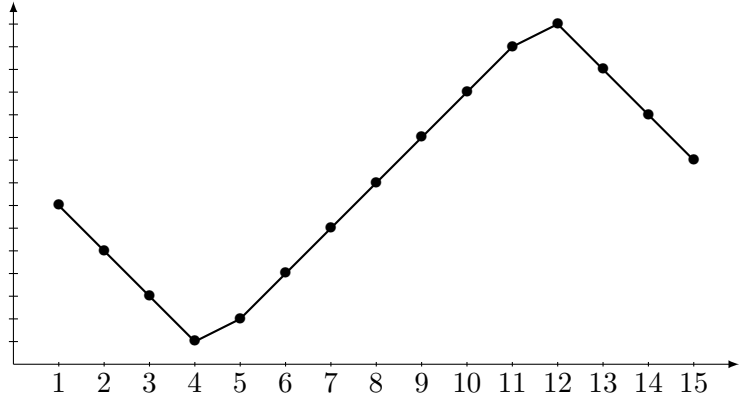
\begin{figure}[htb]
    \center
    \begin{tikzpicture}[yscale=.3, xscale=.6]

    \draw[-latex] (0,0) to (0,16);
    \draw[-latex] (0,0) to (16,0);

    \foreach \i in {1,2,...,15}{
        \draw[] (-0.1, \i) to (0.1,\i) ;
        \draw[] (\i,-0.1) to (\i,0.1);         
    }

    \def\yvalues{7,5,3,1,2,4,6,8,10,12,14,15,13,11,9}

    \foreach \y [
        count=\x,
        remember=\x as \previousx,
        remember=\y as \previousy
    ] in \yvalues{
        \node at (\x,\y) {\small$\bullet$};
        \ifnum \x>1
            \draw[thick] (\previousx,\previousy) to (\x,\y);
        \fi
    }

        \foreach \i in {1,2,...,15}{
        \node at (\i,-0.8){$\i$};
        }
    \end{tikzpicture} 
    \caption{The graph of a permutation witnessing inflection set $\{(5,11)\}.$}
    \label{fig:InflectionSet(5,11)}
\end{figure}

\begin{ex}[Admissible vanishing]
Consider the set $\{(5,11)\}$. The permutation 
\[\sigma=7\,\ 5\,\ 3\,\ 1\,\ 2\,\ 4\,\ 6\,\ 8\,\ 10\,\ 12\,\ 14\,\ 15\,\ 13\,\ 11\,\ 9\]
illustrated in Figure~\ref{fig:InflectionSet(5,11)} witnesses the inflection set $\{(5,11)\}$. In fact, there is a witness for $\{(5,11)\}$ in $\Sym_n$ if and only if $13\le n\le 19$, see Theorem~\ref{T:windows}; in particular, $\{(5,11)\}$ is vanishing.
\end{ex}

\begin{ex}[Admissible non-vanishing]
   The set $\{(2,3)\}$ is witnessed by $n1(n-1)(n-2)\dots2\in \Sym_n$ for every $n\ge 4$. More generally, every set $\{(k,m)\}$ with $2\le k<m\le 2k$ is non-vanishing, see Corollary~\ref{C:dichotomy}.
\end{ex}

However, it is not the case that every set of pairs of integers is admissible.

\begin{ex}[Inadmissible]
    The set $\{(3,6),(6,10)\}$ is not witnessed by any permutation. After relacing $\sigma$ by its complement if necessary, we may assume that $\sigma''(3)<0, \sigma''(6)>0$, and $\sigma''(10)<0$. Denote $\sigma'(3)=\sigma'(4)=\sigma'(5)=a$ and $\sigma'(6)=\sigma'(7)=\sigma'(8)=\sigma'(9)=b$,
    where $a<b$. Note that 
    \begin{equation}\label{eq:(3,6)(6,10)tail}
        \sigma'(1)\geq \sigma'(2)>a\quad \text{ and } \quad b>\sigma'(10)\geq \sigma'(11)\geq\dots\geq \sigma'(n-1).
    \end{equation} We refer to the consecutive subword $\sigma([11,n])$ as the suffix of $\sigma$, which must be concave. Furthermore, we denote the values below $\sigma(10)$ missed by the first 10 values of $\sigma$ by $C$. That is,
    \[C=[1,\sigma(10)]\setminus \{\sigma([1,10])\}=\{c_1,c_2,\dots,c_k\},\]
    where $c_1<c_2<\dots<c_k$. In particular, we must have that $\sigma(n-i+1)=c_{i}$ for $1\leq i\leq k$, and $\sigma'(n-i)=c_i-c_{i+1}$ for $1\leq i\leq k-1$. By (\ref{eq:(3,6)(6,10)tail}), if there is some consecutive $c_i,c_{i+1}$ then each of the values $c_{i+1},\dots,c_k$, must be consecutive. Equivalently, if $c_{i-1}=c_{i}-1$ for some $i<k$, but $[c_i,c_k]\setminus C$ is nonempty this contradicts (\ref{eq:(3,6)(6,10)tail}) saying that the suffix $\sigma$ is concave. We will use this condition to establish $a<0<b$ and then rule out all values of $b$.

    First, assume that $a>0$ to obtain a contradiction. 
    Then $\sigma([1,10])$ is increasing. If $\sigma'(1),\sigma'(2)$, or $b\geq 3$, then there are consecutive values $c_{i-1},c_{i}\in C$ for some $i<k$, but $[c_i,c_k]\setminus C$ is nonempty. Hence, $\sigma'(1)=\sigma'(2)=b=2$, and $a=1$ so that 
    \[
    \sigma([1,7])=\sigma(1),\sigma(1)+2,\sigma(1)+4,\sigma(1)+5,\sigma(1)+6,\sigma(1)+7,\sigma(1)+9,
    \]
    which implies that the $\sigma(1)+8,\sigma(1)+3,\sigma(1)+1$ must be a consecutive subword of the suffix of $\sigma$, again contradicting (\ref{eq:(3,6)(6,10)tail}). Therefore, $a<0$.

    Next, we will show that $b>0$. Assume $b<0$ to obtain a contradiction. Since $a<b<0$ we have $a\leq-2$ and $b\leq -1$ so that $\sigma(3)-\sigma(10)=-3a-4b\geq 10$. Also, $\sigma'(n-1)\leq \dots\leq \sigma'(10)<-1$ such that  $1\leq \sigma(n)<\sigma(n-1)<\dots<\sigma(10)$. Thus, all integers in the interval $[\sigma(10),\sigma(3)]$ must occur in the first 10 positions, but there are at least 11 integers in the interval $[\sigma(10),\sigma(3)]$, which is a contradiction. Therefore, we must have $a<0<b$.

    Consider the consecutive subword
    \[
    \sigma([3,10])=\sigma(6)-3a,\sigma(6)-2a,\sigma(6)-a,\sigma(6),\sigma(6)+b,\sigma(6)+2b,\sigma(6)+3b,\sigma(6)+4b.
    \]
    We split into the below cases and show each is impossible.
    
    If $b=1$, then $a\leq -5$ so that $\sigma(3)-\sigma(10)=-3a-4\geq 11$. However, the suffix of $\sigma$ must be decreasing. Hence, all of the integers in the interval $[\sigma(10),\sigma(3)]$ of which there are least 12, must occur in the first 10 positions, which is a contradiction.

    If $b=2$, then $\sigma(4)=\sigma(6)-2a\neq \sigma(7),\dots,\sigma(10)$ so that $a\leq -5$. Hence, $\sigma(4)>\sigma(10)$. Since $\sigma'(10)<2$ the suffix of $\sigma$ may have some portion of consecutive values, i.e. slope 1. However, it cannot reach any value larger than $\sigma(4)$. Therefore, the at least six integers in the interval $[\sigma(4),\sigma(3)]$ must be among $\sigma(1),\ldots,\sigma(4)$, which is a contradiction.

    If $b=3$, then $\sigma(3)=\sigma(6)-3a\neq \sigma(7),\dots,\sigma(10)$ so that $a\leq -5$. Concavity of the suffix of $\sigma$ implies that two consecutive integers below $\sigma(9)$ cannot both be missing from the first 10 entries, $\sigma([1,10])$. Thus each pair $\{\sigma(6)+1,\sigma(6)+2\}$, $\{\sigma(6)+4,\sigma(6)+5\}$, $\{\sigma(6)+7,\sigma(6)+8\}$ must contain at least one value occurring in $\sigma([1,5])$. Since $\sigma(3),\sigma(4)\geq \sigma(6)+10$, only $\sigma(1),\sigma(2)$, and $\sigma(5)$ can achieve these values. Hence, $\sigma(6)<\sigma(1),\sigma(2)<\sigma(6)+9=\sigma(9)$. Since $\sigma'(1)\geq\sigma'(2)$, $\sigma(2)-\sigma(1)\geq -3a+\sigma(6)-\sigma(2)$, i.e.
    \[
    15\leq -3a\leq 2\sigma(2)-\sigma(1)-\sigma(6)\leq 15.
    \]
    But then $a=-5$, $\sigma(1)=\sigma(6)+1$, and $\sigma(2)=\sigma(6)+8$, which implies that $\sigma(6)+7$, $\sigma(6)+4$, $\sigma(6)+2$ must be a consecutive subword of the suffix of $\sigma$, contradicting the concavity of the suffix.

    If $b=4$, then making sure $\{\sigma(6)-3a,\sigma(6)-2a,\sigma(6)-a\}\cap\{\sigma(6)+4,\sigma(6)+8,\sigma(6)+12,\sigma(6)+16\}=\emptyset$, implies $a\neq -2,-4,-6$. We consider the cases $a=-1$, $a=-3$, and $a\leq-5$ separately. \begin{itemize}
        \item If $a=-1$, then $\sigma(6)+5$ and $\sigma(6)+6$ are missing from $\sigma(1)\dots\sigma(10)$, which is a contradiction. 
        \item If $a=-3$ then the two intervals $[\sigma(6)+1,\sigma(5)-1]=[\sigma(6)+1,\sigma(6)+2]$  and $[\sigma(3)+1,\sigma(9)-1]=[\sigma(6)+10,\sigma(6)+11]$ must each contain one of $\sigma(1)$ or $\sigma(2)$. In the second interval, whichever value is not achieved by $\sigma(1)$ or $\sigma(2)$ is followed by $\sigma(6)+7$ and $\sigma(6)+5$ in the suffix, which contradicts concavity.
        \item Next, consider $a\leq -5$, note also that $a\neq -6$. Then the two intervals $[\sigma(6)+1,\sigma(6)+3]$ and $[\sigma(6)+5,\sigma(6)+7]$ each must contain a value occurring at $\sigma(1)$ or $\sigma(2)$ since only $\sigma(5)\geq \sigma(6)+5$ can have nonempty intersection with the second interval at $\sigma(6)+5$ or $\sigma(6)+7$. Hence, $\sigma(6)+1\leq \sigma(1),\sigma(2)\leq \sigma(6)+7$ and $\sigma'(1)\geq \sigma'(2)$ imply that 
        \[
        15\leq -3a\leq 2\sigma(2)-\sigma(1)-\sigma(6)\leq 13,
        \]
        which is a contradiction.
    \end{itemize}
        
    If $b\geq 5$, then the intervals $[\sigma(6)+1,\sigma(7)-1]$, $[\sigma(7)+1,\sigma(8)-1]$, $[\sigma(8)+1,\sigma(9)-1]$, contain at least 4 consecutive integers. As before, concavity of the suffix implies that two consecutive integers below $\sigma(9)$ cannot both be missing from the first 10 entries $\sigma(1)\dots\sigma(10)$. But since these intervals have empty intersection with $\sigma(6)\dots\sigma(10)$ then a pair of entries from each of the three intervals must be achieved by the values $\sigma(1)\dots\sigma(5)$, which is a contradiction.
    
\end{ex}

\section{Permutations with no inflections}\label{S:NoInflection}
The characterization of permutations with no inflections in Theorem~\ref{thm:no inflection points} is equivalent to the
classification of convex permutations by Brualdi and Dahl
\cite[Theorem~6]{brualdi-dahl}: a permutation has no inflection points if and
only if it or its complement is convex. We give a different,
inductive proof, and record the resulting count, which is novel.

\begin{thm}\label{thm:no inflection points}
    Fix $n \ge 4$. A permutation $w \in \Sym_n$ has no inflection points if and only if either $w$ or a fundamental symmetry of $w$ is one of the following permutations:
    \begin{enumerate}[label={\rm (\Alph*)}]
        \item $123\cdots n$, \label{0 inflection points: 12345}
        \item $1345\cdots n2$, \label{0 inflection points: 13452}
        \item $1357\cdots 8642$, or \label{0 inflection points: 13542}
        \item $1n(n-1)(n-2)\cdots 32$. \label{0 inflection points: 15432}
    \end{enumerate}
    Then the number of permutations in $\Sym_n$ with zero inflection points is 
    \begin{align*}
    \begin{cases}
        10 & \text{ if } \ n = 4, \text{ and}\\
        14 & \text{ if } \ n \ge 5.
          \end{cases}
    \end{align*}
\end{thm}

\begin{proof}
   The fact that the listed permutations and their fundamental symmetries all have zero inflection points is easy to check.   
   Fix $w \in \Sym_n$ and assume that $w$ has zero inflection points. We will prove that $w$ or one of its fundamental symmetries has one of the listed forms. The result is easy to check in small cases, so fix $n \ge 5$ and we assume inductively that the result holds in $\Sym_m$ whenever $m \in [4,n-1]$.
    
    Suppose, first, that $w''(i) = 0$ for all $i \in [2,n-1]$. Then, in particular, we must have $\{1,n\} = \{w(1), w(n)\}$, because otherwise $w''(w^{-1}(1)) > 0$ and/or $w''(w^{-1}(n)) < 0$. Suppose, without loss of generality, that $w(1) = 1$. Then, if $w(j) = 2$ for some $j > 2$, we will find that $w''(j) > 0$, which would be a contradiction. Proceeding in this way, we find that $w$ must be the identity permutation. Had we instead assumed that $w(1) = n$, then $w$ would have been the reverse of the identity permutation.

    Now suppose that $w''(i) \neq 0$ for some $i \in [2,n-1]$. If $w(a) = 1$ and $w(b) = n$ for $a,b \not\in \{1,n\}$, then we will have $w''(a) > 0 > w''(b)$. Thus we can assume, without loss of generality, that $w(1) = 1$. (The other three cases are fundamental symmetries of this choice.) Consider the permutation $v \in \Sym_{n-1}$ defined by $v(j) := w(j+1) - 1$. Then $v''(j) = w''(j+1)$ for all $j \in [2,n-2]$. Therefore, by the inductive hypothesis, $v$ (or a fundamental symmetry of it) must have one of the listed forms.

    The remainder of the proof follows from a straightforward case analysis. We consider each possibility of $v$ being a fundamental symmetry of one of the four listed forms, and analyze when $w$ (obtained by appending a smallest element to the left end of that $v$) is acceptable; i.e., when $w$ has zero inflection points: 
    \begin{itemize}
        \item When $v$ is of form~\ref{0 inflection points: 12345}, then $w$ is also of that form. When $v$ is the reverse of form~\ref{0 inflection points: 12345}, then $w$ has form~\ref{0 inflection points: 15432}.
        \item When $v$ is any fundamental symmetry of form~\ref{0 inflection points: 13452}, the resulting $w$ will not be acceptable. 
        \item When $v$ is the reverse of form~\ref{0 inflection points: 13542}, then $w$ has form~\ref{0 inflection points: 13542}. No other fundamental symmetries of this form will produce an acceptable $w$.
        \item When $v$ is the reverse of form~\ref{0 inflection points: 15432}, then $w$ has form~\ref{0 inflection points: 13452}. No other fundamental symmetries of this form will produce an acceptable $w$.
    \end{itemize}

Form~\ref{0 inflection points: 12345} permits two permutations (the complement of the identity always coincides with the reverse of the identity), form~\ref{0 inflection points: 13452} permits four permutations, form~\ref{0 inflection points: 13542} (which coincides with form~\ref{0 inflection points: 13452} when $n=4$) produces four permutations, and form~\ref{0 inflection points: 15432} produces four permutations. This gives a total of $10$ inflection point-free permutations in $\Sym_4$.  When $n \ge 5$, none of the listed forms coincide, and there will always be $2+4+4+4 = 14$ inflection point-free permutations in $\Sym_n$.
\end{proof}

\section{Permutations with one inflection}\label{S:Perms-with-one-inflection}
In this section, we address permutations of $[n]$ with one inflection. We determine, for every $2\le k<m$, all the values of $n$ for which there exists $\sigma\in\Sym_n$ with $\mathcal{I}(\sigma)=\{(k,m)\}$: see Theorem \ref{T:windows} for $m\ge 2k+1$, and Theorem \ref{T:short} for $k<m\le 2k$. In particular, $I=\{(k,m)\}$ is always an admissible inflection set, and $I$ is vanishing if and only if $m\ge 2k+1$, see Corollary \ref{C:dichotomy}.

\begin{lem}\label{L:AP}
    Let $\sigma\in \Sym_n$ with $(k,m)\in \mathcal{I}(\sigma)$. Then $\sigma(k), \ \sigma(k+1),\ldots, \sigma(m)$ must form an arithmetic progression.
    \end{lem}

    \begin{proof}
Since $\sigma''(t)=0$ for all $k<t<m$, $\sigma(t-1)+\sigma(t+1)-2\sigma(t)=0,$ hence
\[
\sigma(t+1)-\sigma(t)=\sigma(t)-\sigma(t-1).
\]
Let $d=\sigma(k+1)-\sigma(k)=\sigma(t+1)-\sigma(t)$ for $k< t< m$, so for $k\le j\le m$,
\[
\sigma(j)=\sigma(k)+(j-k)d.
\]
\end{proof}

\begin{lem}\label{L:sign}  Let $\sigma\in\Sym_n$ with $\mathcal{I}(\sigma)=\{(k,m)\}$. By Lemma \ref{L:AP}, we can say $\sigma(j)=\sigma(k)+(j-k)d$ for $k\le j\le m$. 
The complement satisfies $(\overline{\sigma})''=-\sigma''$, hence
$\mathcal I(\overline{\sigma})=\mathcal I(\sigma)$, while the slope of the arithmetic progression changes sign.
We may therefore assume $d>0$. Under this assumption, we have $\sigma''(k)>0$ and $\sigma''(m)<0$.
\end{lem}

\begin{proof}
Given $\sigma\in\Sym_n$, it is clear that $(\overline{\sigma})''=-\sigma''$ and $\mathcal I(\overline{\sigma})=\mathcal I(\sigma)$. Assume $\mathcal{I}(\sigma)=\{(k,m)\}$ and $d>0$.
Assume by contradiction that there exists $\ell_0$ with $2\le \ell_0<k$ and
$\sigma''(\ell_0)\sigma''(k)<0$. Let \[
i=\max\{\ell \in [2,k) \ : \ \sigma''(\ell)\sigma''(k)<0\}\]
and
\[
j=\min\{\ell\in (i,k] \ : \ \sigma''(\ell)\neq0\}\le k.
\]
Then $\sigma''(i)\sigma''(j)<0$ and
$\sigma''\equiv0$ on $(i,j)$, so $(i,j)\in\mathcal I(\sigma)$ with $j\le k<m$, a
contradiction. Hence $\sigma''(\ell)\sigma''(k)\ge0$ on $[2,k]$, and similarly
$\sigma''(\ell)\sigma''(m)\ge0$ on $[m,n-1]$.

Assume by contradiction that $\sigma''(k)<0$ (hence $\sigma''(m)>0$). Then $\sigma'$ is non-increasing on
$[1,k-1]$ and $\sigma'(k-1)>\sigma'(k)=d$, so $\sigma'>0$ on $[1,k-1]$. Similarly, $\sigma'$ is
 non-decreasing on $[m,n-1]$ with $\sigma'(m)>\sigma'(m-1)=d$, so $\sigma'>0$ on $[m,n-1]$.
As $\sigma'\equiv d>0$ on $[k,m-1]$, we conclude that $\sigma$ is strictly increasing, which forces $\sigma=e$. This implies $\sigma''\equiv0$ and $\mathcal I(\sigma)=\emptyset$, a contradiction.
\end{proof}

The next lemma fully characterizes when $\mathcal{I}(\sigma)=\{(k,m)\}$. It will be used both to describe the permutations with a single inflection and to create permutations with a prescribed single inflection.

\begin{lem}\label{L:itisinflection} Let $\sigma\in \Sym_n$ and $2\le k<m\le n-1$. Let $d=\sigma(k+1)-\sigma(k)$ and assume $d>0$. Then $\mathcal I(\sigma)=\{(k,m)\}$ if and only if the following three conditions hold: 
\begin{gather}
\sigma'(1)\le\sigma'(2)\le\ldots\le \sigma'(k-1)<d=\sigma'(k), \label{mono1}\\
   \sigma'(j)=d \qquad \text{for }
    \ k\le j\le m-1,\label{derg}\\
   d=\sigma'(m-1)>\sigma'(m)\ge \sigma'(m+1)\ge\ldots\ge \sigma'(n-1).\label{mono2}
   \end{gather}
\end{lem}

\begin{proof} Assume first $\mathcal I(\sigma)=\{(k,m)\}$. 
By Lemma \ref{L:AP}, we have 
$\sigma(j)=\sigma(k)+(j-k)d$ for $k\le j\le m$, that is, \eqref{derg} holds. Since $d>0$, Lemma \ref{L:sign} gives $\sigma''(k)>0$ and $\sigma''(m)<0$. Moreover, the proof of Lemma \ref{L:sign} shows that $\sigma''(\ell)\ge 0$ for $2\le\ell\le k$ and $\sigma''(\ell)\le 0$ for $m\le\ell\le n-1$. Since $\sigma''(\ell)=\sigma'(\ell)-\sigma'(\ell-1)$ and $\sigma'(k)=\sigma'(m-1)=d$, this means that
\begin{gather*}
\sigma'(1)\le \sigma'(2)\le \ldots\le\sigma'(k-1)<\sigma' (k)=d,\\
d=\sigma'(m-1)>\sigma'(m)\ge\sigma'(m+1)\ge\ldots\ge\sigma'(n-1),
\end{gather*}
that is, \eqref{mono1} and \eqref{mono2} hold.

Reciprocally, assume \eqref{mono1}--\eqref{mono2} and let us prove that $\mathcal{I}(\sigma)=\{(k,m)\}$. Since $\sigma''(\ell)=\sigma'(\ell)-\sigma'(\ell-1)$, \eqref{mono1}--\eqref{mono2} yield \begin{itemize}
    \item $\sigma''(\ell)\ge 0$ for $2\le \ell<k$, 
    \item $\sigma''(k)>0$,
    \item $\sigma''(\ell)=0$ for $k<\ell<m$, 
    \item $\sigma''(m)<0$,
    \item $\sigma''(\ell)\le 0$ for $m<\ell\le n-1.$
\end{itemize} 
In particular, $(k,m)\in \mathcal{I}(\sigma).$ Let now $(i,j)\in \mathcal{I}(\sigma)$, so that $i<j$ and $\sigma''(i)\sigma''(j)<0$. Since $\sigma''$ is positive only at indices in $[2,k]$ and negative only at indices in $[m,n-1]$, and since $i<j$, we must have $i\le k$ and $j\ge m$.
If $i<k$, then $i<k<j$, but $\sigma''(k)\neq 0$, contradicting $\sigma''(\ell)=0$ for $i<\ell<j$. Therefore, $i=k$. If $j>m$, then $i<m<j$ but $\sigma''(m)\neq 0$, leading to the same contradiction. Therefore, $j=m$ as well.
\end{proof}

Let $2\le k<m$, $n\ge m+1$ and $\sigma\in\Sym_n$ with $\mathcal{I}(\sigma)=\{(k,m)\}$. By Lemma \ref{L:AP}, $\sigma(j)=\sigma(k)+(j-k)d$ for $k\le j\le m$, where $d=\sigma(k+1)-\sigma(k)$. By Lemma \ref{L:sign}, we may assume $d>0$, in which case $\sigma''(k)>0$ and $\sigma''(m)<0$. Using Lemma \ref{L:itisinflection} and the fact that $\sigma'$ takes integer values, we conclude that $\sigma'(j)\le d-1$ for $1\le j\le k-1$ and for $m\le j\le n-1$.

The next lemma describes the shape of a permutation with a single inflection.

\begin{lem}\label{L:shape}
Let $2\le k<m$ and $n\ge m+1$. Let $\sigma\in\Sym_n$ with $\mathcal{I}(\sigma)=\{(k,m)\}$ and $d=\sigma(k+1)-\sigma(k)>0$. Then there exist $k^*\in[1,k]$ and $m^*\in[m,n]$ so that
\begin{enumerate}
\item[(a)] \begin{itemize}
    \item $\sigma$ is strictly decreasing on $[1,k^*]$, and
    \item $\sigma$ is strictly increasing on $[k^*,m^*]$, and
    \item $\sigma$ is strictly decreasing on $[m^*,n]$,
    \end{itemize}
    with $0<\sigma'(j)\le d-1$ for $j\in[k^*,k-1]\cup[m,m^*-1]$.
\item[(b)] 
\begin{itemize}
    \item $|\sigma'|$ is  non-increasing on $[1,k^*-1]$ and 
    \item $|\sigma'|$ is  non-decreasing on $[m^*,n-1]$.
\end{itemize}
\item[(c)] $\{1,\ldots,\sigma(k^*)-1\}\subseteq \sigma([m^*+1,n])$ and $\{\sigma(m^*)+1,\ldots, n\}\subseteq \sigma([1,k^*-1])$.
\item[(d)] Moreover, $\sigma(k^*)\le 2$ and $\sigma(m^*)\ge n-1$.
\end{enumerate}
\end{lem}

\begin{proof}
By \eqref{mono1}, either there exists $1\le k^*\le k-1$ so that $\{j\in [1,k-1] \ : \ \sigma'(j)>0\}=[k^*,k-1]$, or $\{j\in [1,k-1] \ : \ \sigma'(j)>0\}=\emptyset$ (in this case, we will write $k^*=k$). Similarly, by \eqref{mono2}, either there exists $m^*\in[m+1,n]$ so that $\{j\in [m,n-1] \ : \ \sigma'(j)>0\}=[m,m^*-1]$ or $\{j\in [m,n-1] \ : \ \sigma'(j)>0\}=\emptyset$ (in this case, we will write $m^*=m$). This implies that
\begin{itemize}
    \item $\sigma$ is strictly decreasing on $[1,k^*]$,
    \item $\sigma$ is strictly increasing on $[k^*,m^*]$, and
    \item $\sigma$ is strictly decreasing on $[m^*,n]$,
\end{itemize} with $0<\sigma'(j)\le d-1$ for $j\in[k^*,k-1]\cup[m,m^*-1]$. This proves (a).

Given $j\in [m^*,n-1]$, $\sigma'(j)<0$ and, by \eqref{mono2}, $\sigma'$ is non-increasing in $[m^*,n-1]$. Similarly, given $j\in [1,k^*-1]$, $\sigma'(j)<0$, and $\sigma'$ is non-decreasing by \eqref{mono1}. We conclude that
\begin{itemize}
    \item $|\sigma'|$ is non-decreasing on $[m^*,n-1]$,
    \item $|\sigma'|$ is non-increasing on $[1,k^*-1]$,
\end{itemize}
which is (b).

Given $j\in [k^*,m^*]$, $\sigma(j)\ge \sigma(k^*)$, hence
$\{1,\ldots,\sigma(k^*)-1\}\cap\{\sigma(k^*),\ldots,\sigma(m^*)\}=\emptyset$. Moreover, $\{1,\ldots,\sigma(k^*)-1\}\cap\{\sigma(1),\ldots,\sigma(k^*-1)\}=\emptyset$, as $\sigma$ is strictly decreasing on $[1,k^*]$. This implies that
$\{1,\ldots,\sigma(k^*)-1\}\subseteq \{\sigma(m^*+1),\ldots,\sigma(n)\}$, and $\{1,\ldots,\sigma(k^*)-1\}$ are the smallest values in this last third of the permutation.

Assume by contradiction that $\sigma(k^*)\ge 3$. Since $\{1,\ldots,\sigma(k^*)-1\}\subseteq \{\sigma(m^*+1),\ldots,\sigma(n)\}$ and $\sigma$ is strictly decreasing on $[m^*,n]$, $|\sigma'(n-1)|=1$. Since $|\sigma'|$ is non-decreasing on $[m^*,n-1]$, this shows that $\{\sigma(m^*),\ldots,\sigma(n)\}=\{1,\ldots,\sigma(m^*)\}$. Now, since $1\le \sigma(k^*)<\sigma(m^*)$, $\sigma(k^*)\in\{\sigma(m^*),\ldots,\sigma(n)\}$. However, $k^*\le k<m\le m^*$, contradicting the fact that $\sigma$ is a permutation (hence injective). We conclude that $\sigma(k^*)\le 2.$

Now, given $j\in [m^*+1,n]$, $\sigma(j)<\sigma(m^*)$, hence $\{\sigma(m^*)+1,\ldots, n\}\cap\{\sigma(m^*+1),\ldots,\sigma(n)\}=\emptyset$. Moreover, $\{\sigma(m^*)+1,\ldots,n\}\cap\{\sigma(k^*),\ldots,\sigma(m^*)\}=\emptyset$, as $\sigma$ is strictly increasing on $[k^*,m^*]$. Hence
 $\{\sigma(m^*)+1,\ldots, n\}\subseteq \{\sigma(1),\ldots,\sigma(k^*-1)\}$, and $\{\sigma(m^*)+1,\ldots, n\}$ are the largest values of this first third of the permutation.

Assume by contradiction that $\sigma(m^*)\le n-2$.
Then $\{n-1,n\}\subseteq \{\sigma(1),\ldots,\sigma(k^*-1)\}$. Since $\sigma$ is strictly decreasing on $[1,k^*]$, we conclude  that $\sigma(1)=n$ and $\sigma(2)=n-1$. Since $|\sigma'|$ is non-increasing on  $[1,k^*-1]$, we conclude that $\{\sigma(1),\ldots,\sigma(k^*-1),\sigma(k^*)\}=\{\sigma(k^*),\ldots,n\}$. We also have $\sigma(m^*)\in \{\sigma(k^*),\ldots, n\}$, but $m^*\notin[1,k^*]$, contradicting the injectivity of $\sigma$. Therefore $\sigma(m^*)\ge n-1$. This proves (c) and (d).
\end{proof}

From now on, given $\sigma$ as in Lemma \ref{L:shape}, we fix $k^*$ and $m^*$ as in that lemma, and we call
\[
r=k-k^*\ge 0,\qquad p=m^*-m\ge 0 .
\]
We will refer to $\sigma([1,k^*-1])$ as the first third of $\sigma$ and to $\sigma([m^*+1,n])$ as the last third of $\sigma$. Notice that they have $k-1-r$ and $n-m-p$ elements, respectively. The set of values that $\sigma$ misses within the inflection is
\[
H_{infl}=\{\sigma(k)+id+j \ : \ 0\le i\le m-k-1, \ 1\le j\le d-1\}.
\]
Since $\sigma$ is strictly increasing on $[k^*,m^*]$ and, for each $0\le i\le m-k-1$, $\sigma(k)+id=\sigma(k+i)$ and $\sigma(k)+(i+1)d=\sigma(k+i+1)$, we have $H_{infl}\cap\sigma([k^*,m^*])=\emptyset$. Therefore
\begin{equation}\label{Hsplit}
H_{infl}\subseteq \sigma([1,k^*-1])\cup\sigma([m^*+1,n]).
\end{equation}
We now analyze the possible values of $d$ separately.

\begin{figure}[ht]
\centering
\begin{tikzpicture}[x=0.8cm,y=0.42cm]
  \draw[gray!60] (0.5,0.5) -- (8.6,0.5);
  \draw[gray!60] (0.5,0.5) -- (0.5,8.8);
  \foreach \i in {1,...,8} \node[below,font=\scriptsize,gray!60!black] at (\i,0.5) {\i};
  \foreach \v in {1,...,8} \node[left,font=\scriptsize,gray!60!black] at (0.5,\v) {\v};
  \draw[densely dotted,gray] (2,0.5) -- (2,9.6);
  \draw[densely dotted,gray] (7,0.5) -- (7,9.6);
  \node[font=\scriptsize,above] at (2,9.6) {$k=2$};
  \node[font=\scriptsize,above] at (7,9.6) {$m$};
  \draw[thick] (1,8)--(2,2)--(3,3)--(4,4)--(5,5)--(6,6)--(7,7)--(8,1);
  \foreach \p/\v in {1/8,2/2,3/3,4/4,5/5,6/6,7/7,8/1} \fill (\p,\v) circle (2.2pt);
  \node[font=\scriptsize,anchor=west,fill=white,inner sep=1.5pt] at (1.25,8.1) {$\sigma(1)=m+1$};
  \node[font=\scriptsize,anchor=west] at (2.3,1.8) {$\sigma(k)=2$};
  \node[font=\scriptsize,anchor=east] at (6.75,7.2) {$\sigma(m)=m$};
  \node[font=\scriptsize,anchor=west] at (8.2,1) {$\sigma(n)=1$};
  \draw[decorate,decoration={brace,mirror,amplitude=3pt}] (1.05,-1.0) -- (1.95,-1.0)
       node[midway,below=4pt,font=\scriptsize,align=center] {one step\\ of $1-m$};
  \draw[decorate,decoration={brace,mirror,amplitude=3pt}] (2.05,-1.0) -- (6.95,-1.0)
       node[midway,below=4pt,font=\scriptsize,align=center] {$m-2$ steps\\ of $1$};
  \draw[decorate,decoration={brace,mirror,amplitude=3pt}] (7.05,-1.0) -- (7.95,-1.0)
       node[midway,below=4pt,font=\scriptsize,align=center] {one step\\ of $1-m$};
\end{tikzpicture}
\caption{The permutation $\sigma=(m+1)\,2\,3\ldots m\,1$ of Lemma \ref{L:d1} when $m=7$, that is, $\sigma=8\;2\;3\;4\;5\;6\;7\;1$. Since $d=1$, the progression $\sigma(k),\ldots,\sigma(m)$ skips no values, so the only entries outside of it are $n$, in the first position, and $1$, in the last one.}
\label{F:d1}
\end{figure}
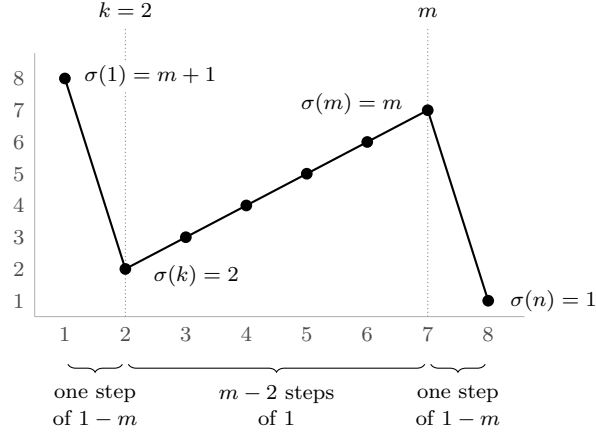

\begin{lem}\label{L:d1}
Let $\sigma$ be as in Lemma \ref{L:shape}. If $d=1$, then $k=2$, $n=m+1$ and $\sigma=(m+1)\,2\,3\ldots m\,1$, see Figure \ref{F:d1}.
\end{lem}

\begin{proof}
Since $d=1$, Lemma \ref{L:shape}(a) shows that there is no $j$ in $[k^*,k-1]\cup[m,m^*-1]$, hence $k^*=k$ and $m^*=m$, and $\sigma([k,m])=\{\sigma(k),\sigma(k)+1,\ldots,\sigma(m)\}$.

Since $k\ge 2$ and $\sigma'(k-1)<0$, we have $\sigma(k-1)>\sigma(k)$. As $\sigma(k-1)\notin\sigma([k,m])$, this gives $\sigma(k-1)\ge\sigma(m)+1$. By Lemma \ref{L:shape}(d), $\sigma(m)\ge n-1$, so $\sigma(k-1)=n$. If $k\ge 3$, then $\sigma(k-2)>\sigma(k-1)=n$, as $\sigma$ is strictly decreasing on $[1,k]$, which is a contradiction. Hence $k=2$.

Similarly, since $n\ge m+1$ and $\sigma'(m)<0$, $\sigma(m+1)<\sigma(m)$, and $\sigma(m+1)\notin\sigma([k,m])$ gives $\sigma(m+1)\le\sigma(k)-1$. By Lemma \ref{L:shape}(d), $\sigma(k)\le 2$, so $\sigma(m+1)=1$ and $\sigma(k)=2$. If $n\ge m+2$, then $\sigma(m+2)<\sigma(m+1)=1$, which is a contradiction. Hence $n=m+1$.

We conclude that $\sigma(1)=n=m+1$, $\sigma(j)=j$ for $2\le j\le m$, and $\sigma(m+1)=1$.
\end{proof}

\begin{lem}\label{L:d3}
Let $\sigma$ be as in Lemma \ref{L:shape}. If $d\ge 3$, then $m\le 2k$ and $n\ge 2m-k-1$.
\end{lem}

\begin{proof}
For each $0\le i\le m-k-1$, let
\[
H_i=\{\sigma(k)+id+1,\ldots,\sigma(k)+id+(d-1)\},
\]
so that $H_{infl}$ is the disjoint union of $H_0,\ldots,H_{m-k-1}$, see Figure \ref{F:d3}.

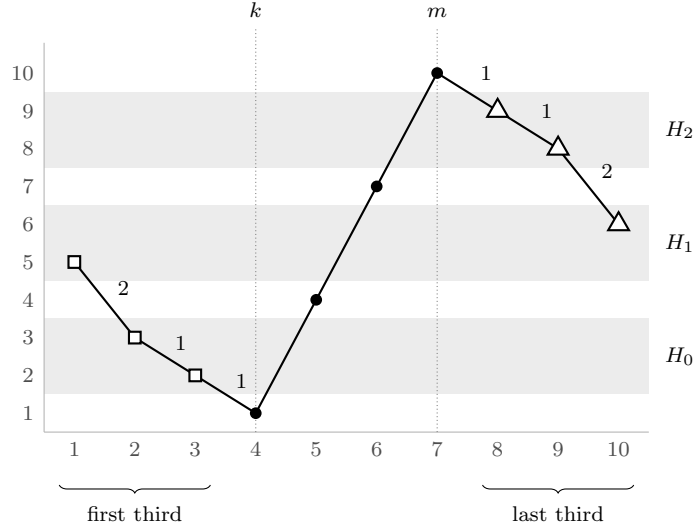
\begin{figure}[ht]
\centering
\begin{tikzpicture}[x=0.8cm,y=0.5cm]
  \foreach \lo/\hi/\i in {1.5/3.5/0, 4.5/6.5/1, 7.5/9.5/2}
     { \fill[gray!15] (0.5,\lo) rectangle (10.5,\hi);
       \node[font=\scriptsize,anchor=west] at (10.6,{(\lo+\hi)/2}) {$H_{\i}$}; }
  \draw[gray!60] (0.5,0.5) -- (10.5,0.5);
  \draw[gray!60] (0.5,0.5) -- (0.5,10.8);
  \foreach \i in {1,...,10} \node[below,font=\scriptsize,gray!60!black] at (\i,0.5) {\i};
  \foreach \v in {1,...,10} \node[left,font=\scriptsize,gray!60!black] at (0.5,\v) {\v};
  \draw[densely dotted,gray] (4,0.5) -- (4,11.2);
  \draw[densely dotted,gray] (7,0.5) -- (7,11.2);
  \node[font=\scriptsize,above] at (4,11.2) {$k$};
  \node[font=\scriptsize,above] at (7,11.2) {$m$};
  \draw[thick] (1,5)--(2,3)--(3,2)--(4,1)--(5,4)--(6,7)--(7,10)--(8,9)--(9,8)--(10,6);
  \foreach \p/\v in {4/1,5/4,6/7,7/10} \fill (\p,\v) circle (2.2pt);
  \foreach \p/\v in {1/5,2/3,3/2} \node[draw,fill=white,inner sep=2.2pt,thick] at (\p,\v) {};
  \foreach \p/\v in {8/9,9/8,10/6} \node[draw,fill=white,regular polygon,regular polygon sides=3,inner sep=1.6pt,thick] at (\p,\v) {};
  \node[font=\scriptsize,anchor=west] at (1.55,4.3) {$2$};
  \node[font=\scriptsize,anchor=west] at (2.5,2.85) {$1$};
  \node[font=\scriptsize,anchor=west] at (3.5,1.85) {$1$};
  \node[font=\scriptsize,anchor=west] at (7.55,10.0) {$1$};
  \node[font=\scriptsize,anchor=west] at (8.55,9.0) {$1$};
  \node[font=\scriptsize,anchor=west] at (9.55,7.4) {$2$};
  \draw[decorate,decoration={brace,mirror,amplitude=3pt}] (0.75,-0.9) -- (3.25,-0.9)
       node[midway,below=4pt,font=\scriptsize] {first third};
  \draw[decorate,decoration={brace,mirror,amplitude=3pt}] (7.75,-0.9) -- (10.25,-0.9)
       node[midway,below=4pt,font=\scriptsize] {last third};
\end{tikzpicture}
\caption{Lemma \ref{L:d3} with $\sigma=5\;3\;2\;1\;4\;7\;10\;9\;8\;6$, where $(k,m,n)=(4,7,10)$ and $d=3$. The shaded bands are the sets $H_0=\{2,3\}$, $H_1=\{5,6\}$ and $H_2=\{8,9\}$, of values skipped by the progression. The numbers next to the graph are the values of $|\sigma'|$ in the first and last thirds.}
\label{F:d3}
\end{figure}

Since $d\ge 3$, each $H_i$ consists of $d-1\ge 2$ consecutive integers. 

We first claim that if $H_i\subseteq \sigma([1,k^*-1])$, then $i=0$. Indeed, let $b,b+1\in H_i$. Since $\sigma$ is strictly decreasing on $[1,k^*]$, there exists $j\in[1,k^*-2]$ so that $\sigma(j)=b+1$ and $\sigma(j+1)=b$, that is, $|\sigma'(j)|=1$. Since $|\sigma'|$ is non-increasing on $[1,k^*-1]$, we conclude that $\{\sigma(j),\ldots,\sigma(k^*)\}=\{\sigma(k^*),\ldots,b+1\}$. Now, $\sigma(k^*)\le\sigma(k)\le\sigma(k)+id<b$, hence $\sigma(k+i)=\sigma(k)+id\in\{\sigma(j),\ldots,\sigma(k^*)\}$. By injectivity of $\sigma$, $k+i\le k^*\le k$, that is, $i=0$.

Next we claim that if $H_i\subseteq \sigma([m^*+1,n])$, then $i=m-k-1$. Indeed, let $b,b+1\in H_i$. Since $\sigma$ is strictly decreasing on $[m^*,n]$, there exists $j\in[m^*+1,n-1]$ so that $\sigma(j)=b+1$ and $\sigma(j+1)=b$, that is, $|\sigma'(j)|=1$. Since $|\sigma'|$ is non-decreasing on $[m^*,n-1]$, we conclude that $\{\sigma(m^*),\ldots,\sigma(j+1)\}=\{b,\ldots,\sigma(m^*)\}$. Now, $b<\sigma(k)+(i+1)d\le\sigma(m)\le\sigma(m^*)$, hence $\sigma(k+i+1)=\sigma(k)+(i+1)d\in\{\sigma(m^*),\ldots,\sigma(j+1)\}$. By injectivity of $\sigma$, $k+i+1\ge m^*\ge m$, that is, $i=m-k-1$.

By \eqref{Hsplit} and the first claim, each of the $m-k-1$ disjoint sets $H_1,\ldots,H_{m-k-1}$ contains an element of $\sigma([m^*+1,n])$. Hence $n-m^*\ge m-k-1$, and so $n\ge m^*+m-k-1\ge 2m-k-1$. By \eqref{Hsplit} and the second claim, each of the $m-k-1$ disjoint sets $H_0,\ldots,H_{m-k-2}$ contains an element of $\sigma([1,k^*-1])$. Hence $m-k-1\le k^*-1\le k-1$, that is, $m\le 2k$.
\end{proof}

We are left with the case $d=2$. In this case, the progression $\sigma(k),\ldots,\sigma(m)$ skips exactly one value between any two consecutive entries, so it skips $\ell=m-k$ values $h_1<\ldots<h_\ell$. By \eqref{Hsplit}, each of them appears in the first third or in the last third of $\sigma$. We mark this choice in two sets of indices, $A$ and $B$. We will see that $\sigma$ is completely described by $A$ and $B$ together with four numbers: the lengths $r$ and $p$ of the two parts next to the progression where $\sigma$ increases by $1$, and two numbers $a,c\in\{0,1\}$, which tell us whether $n$ is the first entry of $\sigma$ and whether $1$ is the last one (see Figure \ref{F:build} for an example). Counting entries gives two relations between these parameters, $k$ and $n$, namely \eqref{count1} and \eqref{count2}, and \eqref{dist} allows us to compute $|\sigma'(j)|$ for $j\in[1,k^*-1]\cup[m^*,n-1]$, in terms of the parameters. By Lemma \ref{L:shape}(b), these values are non-increasing on $[1,k^*-1]$ and non-decreasing on $[m^*,n-1]$, and this is what restricts the parameters in the proofs below.

So, let $\sigma$ be as in Lemma \ref{L:shape}, assume $d=2$, and write $\ell=m-k$. By Lemma \ref{L:shape}(a), $\sigma'(j)=1$ for $j\in[k^*,k-1]\cup[m,m^*-1]$, hence
\begin{equation}\label{ends}
\sigma(k^*)=\sigma(k)-r,\qquad \sigma(m^*)=\sigma(m)+p=\sigma(k)+2\ell+p,
\end{equation}
and
\[
\sigma([k^*,m^*])=\{\sigma(k)-r,\ldots,\sigma(k)\}\cup\{\sigma(k)+2i \ : \ 0\le i\le \ell\}\cup\{\sigma(m),\ldots,\sigma(m)+p\}.
\]
In this case $H_{infl}=\{h_1,\ldots,h_\ell\}$, where
\[
h_i=\sigma(k)+2i-1,\qquad 1\le i\le \ell,
\]
and the set of values that $\sigma$ misses on $[k^*,m^*]$ is $\{1,\ldots,\sigma(k^*)-1\}\cup H_{infl}\cup\{\sigma(m^*)+1,\ldots,n\}$. Let
\[
A=\{i\in[1,\ell] \ : \ h_i\in\sigma([1,k^*-1])\},\qquad B=\{i\in[1,\ell] \ : \ h_i\in\sigma([m^*+1,n])\}.
\]
By \eqref{Hsplit}, $[1,\ell]$ is the disjoint union of $A$ and $B$. Recall that the values in $\sigma([1,k^*-1])$ are larger than $\sigma(k^*)$, and the values in $\sigma([m^*+1,n])$ are smaller than $\sigma(m^*)$. Together with Lemma \ref{L:shape}(c), this gives
\begin{equation}\label{thirds}
\begin{gathered}
\sigma([1,k^*-1])=\{h_i \ : \ i\in A\}\cup\{\sigma(m^*)+1,\ldots,n\},\\
\sigma([m^*+1,n])=\{h_i \ : \ i\in B\}\cup\{1,\ldots,\sigma(k^*)-1\}.
\end{gathered}
\end{equation}
By Lemma \ref{L:shape}(d), the numbers
\[
a=n-\sigma(m^*)\qquad\text{and}\qquad c=\sigma(k^*)-1
\]
belong to $\{0,1\}$. Let $q=|A|$ and $t=|B|=\ell-q$, and write $A=\{i_1>i_2>\ldots>i_q\}$ and $B=\{j_1>j_2>\ldots>j_t\}$. Counting the elements of both sets in \eqref{thirds}, we obtain $k^*-1=q+a$ and $n-m^*=t+c$, that is,
\begin{equation}\label{count1}
q+a+r=k-1,
\end{equation}
\begin{equation}\label{count2}
n=m+p+\ell-q+c=2\ell+k+p+c-q .
\end{equation}
Since $\sigma$ is strictly decreasing on $[1,k^*]$ and on $[m^*,n]$, \eqref{thirds} determines
the order of the entries in the first and last thirds of $\sigma$. Recall that $k^*-1=q+a$
and $n-m^*=t+c$.

In the first third, the values $h_{i_1}>h_{i_2}>\ldots>h_{i_q}$ appear in this order, and they are preceded by $n$ when $a=1$. That is,
\begin{equation}\label{order}
\sigma(a+s)=h_{i_s}\quad\text{for } 1\le s\le q,\qquad\text{and}\qquad \sigma(1)=n\ \text{ if } a=1.
\end{equation}
In the last third, the values $h_{j_1}>h_{j_2}>\ldots>h_{j_t}$ appear in this order, and they are followed by $1$ when $c=1$. That is,
\begin{equation}\label{order2}
\sigma(m^*+s)=h_{j_s}\quad\text{for } 1\le s\le t,\qquad\text{and}\qquad \sigma(n)=1\ \text{ if } c=1.
\end{equation}
In particular, $\sigma(k^*-1)=h_{i_q}$ if $q\ge 1$, and $\sigma(m^*+1)=h_{j_1}$ if $t\ge 1$.
Finally, by \eqref{ends} and the definitions of $a$ and $c$, for $i,i'\in[1,\ell]$ we have
\begin{equation}\label{dist}
\begin{gathered}
h_i-h_{i'}=2(i-i'),\qquad h_i-\sigma(k^*)=2i-1+r,\qquad \sigma(m^*)-h_i=2(\ell-i)+p+1,\\
n-h_i=2(\ell-i)+p+1+a,\qquad h_i-1=2i+r+c-1 .
\end{gathered}
\end{equation}

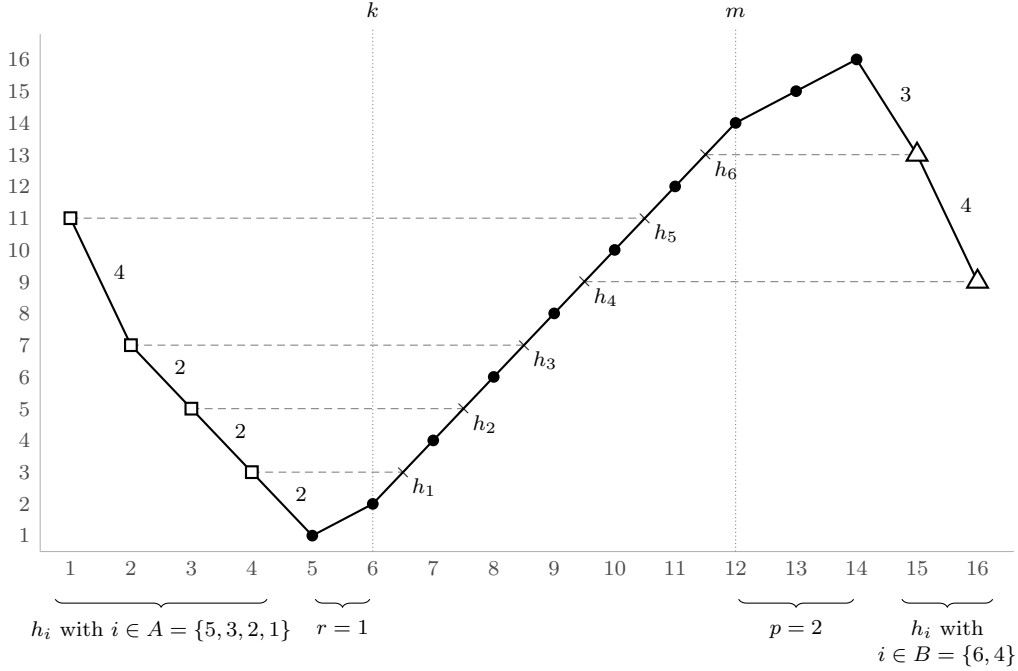
\begin{figure}[ht]
\centering
\begin{tikzpicture}[x=0.8cm,y=0.42cm]
  \draw[gray!60] (0.5,0.5) -- (16.6,0.5);
  \draw[gray!60] (0.5,0.5) -- (0.5,16.8);
  \foreach \i in {1,...,16} \node[below,font=\scriptsize,gray!60!black] at (\i,0.5) {\i};
  \foreach \v in {1,...,16} \node[left,font=\scriptsize,gray!60!black] at (0.5,\v) {\v};
  \draw[densely dotted,gray] (6,0.5) -- (6,17);
  \draw[densely dotted,gray] (12,0.5) -- (12,17);
  \node[font=\scriptsize,above] at (6,17) {$k$};
  \node[font=\scriptsize,above] at (12,17) {$m$};
  \draw[densely dashed,gray] (1,11) -- (10.5,11);
  \draw[densely dashed,gray] (2,7)  -- (8.5,7);
  \draw[densely dashed,gray] (3,5)  -- (7.5,5);
  \draw[densely dashed,gray] (4,3)  -- (6.5,3);
  \draw[densely dashed,gray] (11.5,13) -- (15,13);
  \draw[densely dashed,gray] (9.5,9)   -- (16,9);
  \draw[thick] (1,11)--(2,7)--(3,5)--(4,3)--(5,1)--(6,2)--(7,4)--(8,6)--(9,8)--(10,10)--(11,12)--(12,14)--(13,15)--(14,16)--(15,13)--(16,9);
  \foreach \p/\v in {5/1,6/2,7/4,8/6,9/8,10/10,11/12,12/14,13/15,14/16} \fill (\p,\v) circle (2.2pt);
  \foreach \p/\v in {1/11,2/7,3/5,4/3} \node[draw,fill=white,inner sep=2.2pt,thick] at (\p,\v) {};
  \foreach \p/\v in {15/13,16/9} \node[draw,fill=white,regular polygon,regular polygon sides=3,inner sep=1.6pt,thick] at (\p,\v) {};
  \foreach \q/\v/\n in {6.5/3/1,7.5/5/2,8.5/7/3,9.5/9/4,10.5/11/5,11.5/13/6}
     { \node[font=\tiny,inner sep=0.5pt] at (\q,\v) {$\times$};
       \node[font=\scriptsize,anchor=north west,inner sep=1pt] at (\q+0.1,\v-0.1) {$h_{\n}$}; }
  \node[font=\scriptsize,anchor=west] at (1.55,9.3) {$4$};
  \node[font=\scriptsize,anchor=west] at (2.55,6.3) {$2$};
  \node[font=\scriptsize,anchor=west] at (3.55,4.3) {$2$};
  \node[font=\scriptsize,anchor=west] at (4.55,2.3) {$2$};
  \node[font=\scriptsize,anchor=west] at (14.55,14.9) {$3$};
  \node[font=\scriptsize,anchor=west] at (15.55,11.4) {$4$};
  \draw[decorate,decoration={brace,mirror,amplitude=3pt}] (0.75,-1.0) -- (4.25,-1.0)
       node[midway,below=4pt,font=\scriptsize] {$h_i$ with $i\in A=\{5,3,2,1\}$};
  \draw[decorate,decoration={brace,mirror,amplitude=3pt}] (5.05,-1.0) -- (5.95,-1.0)
       node[midway,below=4pt,font=\scriptsize] {$r=1$};
  \draw[decorate,decoration={brace,mirror,amplitude=3pt}] (12.05,-1.0) -- (13.95,-1.0)
       node[midway,below=4pt,font=\scriptsize] {$p=2$};
  \draw[decorate,decoration={brace,mirror,amplitude=3pt}] (14.75,-1.0) -- (16.25,-1.0)
       node[midway,below=4pt,font=\scriptsize,align=center] {$h_i$ with\\ $i\in B=\{6,4\}$};
\end{tikzpicture}
\caption{Lemma \ref{L:d2} with $\sigma=11\;7\;5\;3\;1\;2\;4\;6\;8\;10\;12\;14\;15\;16\;13\;9$, where $(k,m,n)=(6,12,16)$, $\ell=6$ and $h_i=2i+1$. The numbers next to the graph are the values of $|\sigma'|$ in the first and last thirds. Since $3,2\in A$, the values $h_3$ and $h_2$ are consecutive entries and $|\sigma'(2)|=2$. From then on, $|\sigma'|\le 2$ in the first third, which forces $h_1$ to be the next entry, as in (a). Any two consecutive elements of $A$, and of $B$, differ by at most $2$, as in (c). The last value of $|\sigma'|$ in the first third is $2i_q-1+r=2$, and the first one in the last third is $2(\ell-j_1)+p+1=3$, as in (d) and (e). Furthermore, $p=2$ and $B=\{\ell,\ell-2\}$.}
\label{F:d2}
\end{figure}

The next lemma shows that the sets $A$ and $B$ are very rigid, see Figure \ref{F:d2}. If two neighbouring skipped values $h_{i-1}$ and $h_i$ are both used in the first third, they are consecutive entries of $\sigma$ and differ by $2$. Since $|\sigma'|$ is non-increasing there, from then on $\sigma$ can only decrease by $2$ at a time, which forces all the smaller skipped values to be used in the first third as well; this is (a), and (b) is the symmetric statement for the last third. As a consequence, the elements of $A$, and those of $B$, cannot be far apart, which is (c). Finally, (d) and (e) compare the difference between the two entries next to the increasing part of $\sigma$, which grows with $r$ (respectively, with $p$), with the neighbouring difference, which is at most $4$ by (c). This bounds $r$ and $p$.

\begin{lem}\label{L:d2}
Let $\sigma$ be as in Lemma \ref{L:shape} with $d=2$, and let $\ell$, $h_i$, $A$, $B$, $q$, $t$, $a$, $c$ be as above.
\begin{enumerate}
\item[(a)] If $i\ge 2$ and $i,i-1\in A$, then $\{1,\ldots,i\}\subseteq A$.
\item[(b)] If $j\ge 2$ and $j,j-1\in B$, then $\{j-1,\ldots,\ell\}\subseteq B$.
\item[(c)] Any two consecutive elements of $A$ differ by at most $2$, and any two consecutive elements of $B$ differ by at most $2$.
\item[(d)] If $q\ge 2$, then $2i_q-1+r\le 2(i_{q-1}-i_q)\le 4$. In particular, $r\le 3$, and if $r\ge 2$, then $i_q=1$ and $i_{q-1}=3$.
\item[(e)] If $t\ge 2$, then $2(\ell-j_1)+p+1\le 2(j_1-j_2)\le 4$. In particular, $p\le 3$. If $j_1-j_2=1$, then $p\le 1$. If $p\ge 2$, then $B=\{\ell,\ell-2,\ell-4,\ldots,\ell-2(t-1)\}$.
\end{enumerate}
\end{lem}

\begin{proof}
Throughout the proof we will use the following three facts: first, by \eqref{Hsplit}, no $h_i$ belongs to $\sigma([k^*,m^*])$. In particular, $h_i\ne\sigma(k^*)$ and $h_i\ne\sigma(m^*)$. Second, for $2\le i\le\ell$, the only integer strictly between $h_{i-1}$ and $h_i$ is
\begin{equation}\label{between}
h_i-1=h_{i-1}+1=\sigma(k)+2(i-1)=\sigma(k+i-1),\qquad\text{where } k+1\le k+i-1\le m-1 .
\end{equation}
Third, by Lemma \ref{L:shape}(b), $|\sigma'|$ is non-increasing on $[1,k^*-1]$ and non-decreasing on $[m^*,n-1]$.

\medskip
\noindent Proof of (a): Let $i\ge 2$ with $i,i-1\in A$, that is, $h_i,h_{i-1}\in\sigma([1,k^*-1])$.

We first show that $h_i$ and $h_{i-1}$ appear in consecutive positions. Indeed, $\sigma$ is strictly decreasing on $[1,k^*]$, and by \eqref{between} the only integer strictly between $h_{i-1}$ and $h_i$ is $\sigma(k+i-1)$, which is not in $\sigma([1,k^*-1])$ because $k+i-1>k^*$. Hence there exists $s\in[1,k^*-2]$ so that
\[
\sigma(s)=h_i\qquad\text{and}\qquad\sigma(s+1)=h_{i-1}.
\]
In particular, $|\sigma'(s)|=h_i-h_{i-1}=2$. Since $|\sigma'|$ is non-increasing on $[1,k^*-1]$, we conclude that
\begin{equation}\label{smallsteps1}
|\sigma'(u)|\le 2\qquad\text{for all } s\le u\le k^*-1.
\end{equation}

Next we claim the following: if $u\in[s,k^*-1]$ and $\sigma(u)=h_{i'}$ for some $i'\ge 2$, then $\sigma(u+1)=h_{i'-1}$ (and therefore $u+1\le k^*-1$, since $h_{i'-1}\ne\sigma(k^*)$). Indeed, $\sigma(u+1)<\sigma(u)$ and, by \eqref{smallsteps1}, $\sigma(u)-\sigma(u+1)\le 2$, so
\[
\sigma(u+1)\in\{h_{i'}-1,\ h_{i'}-2\}.
\]
By \eqref{between}, $h_{i'}-1=\sigma(k+i'-1)$ with $k+i'-1\ge k+1>k^*\ge u+1$, so the injectivity of $\sigma$ gives $\sigma(u+1)\ne h_{i'}-1$. Therefore $\sigma(u+1)=h_{i'}-2=h_{i'-1}$, which proves the claim.

Applying the claim repeatedly, starting from $\sigma(s+1)=h_{i-1}$, we obtain
\[
\sigma(s)=h_i,\quad\sigma(s+1)=h_{i-1},\quad\sigma(s+2)=h_{i-2},\quad\ldots,\quad\sigma(s+i-1)=h_1,
\]
with $s+i-1\le k^*-1$. Hence $h_1,\ldots,h_i\in\sigma([1,k^*-1])$, that is, $\{1,\ldots,i\}\subseteq A$.

\medskip
\noindent Proof of (b): The argument is symmetric. Let $j\ge 2$ with $j,j-1\in B$, that is, $h_j,h_{j-1}\in\sigma([m^*+1,n])$.

Since $\sigma$ is strictly decreasing on $[m^*,n]$, and by \eqref{between} the only integer strictly between $h_{j-1}$ and $h_j$ is $\sigma(k+j-1)$, which is not in $\sigma([m^*+1,n])$ because $k+j-1<m^*$, there exists $s\in[m^*+1,n-1]$ so that
\[
\sigma(s)=h_j\qquad\text{and}\qquad\sigma(s+1)=h_{j-1}.
\]
In particular, $|\sigma'(s)|=2$. Since $|\sigma'|$ is non-decreasing on $[m^*,n-1]$, we conclude that
\begin{equation}\label{smallsteps2}
|\sigma'(u)|\le 2\qquad\text{for all } m^*\le u\le s.
\end{equation}

Next we claim the following: if $u\in[m^*+1,s]$ and $\sigma(u)=h_{j'}$ for some $j'\le\ell-1$, then $\sigma(u-1)=h_{j'+1}$ (and therefore $u-1\ge m^*+1$, since $h_{j'+1}\ne\sigma(m^*)$). Indeed, $\sigma(u-1)>\sigma(u)$ and, by \eqref{smallsteps2}, $\sigma(u-1)-\sigma(u)\le 2$, so
\[
\sigma(u-1)\in\{h_{j'}+1,\ h_{j'}+2\}.
\]
By \eqref{between}, $h_{j'}+1=\sigma(k+j')$ with $k+j'\le m-1<m^*\le u-1$, so the injectivity of $\sigma$ gives $\sigma(u-1)\ne h_{j'}+1$. Therefore $\sigma(u-1)=h_{j'}+2=h_{j'+1}$, which proves the claim.

Applying the claim repeatedly, starting from $\sigma(s)=h_j$, we obtain
\[
\sigma(s+1)=h_{j-1},\quad\sigma(s)=h_j,\quad\sigma(s-1)=h_{j+1},\quad\ldots,\quad\sigma(s-(\ell-j))=h_\ell,
\]
with $s-(\ell-j)\ge m^*+1$. Hence $h_{j-1},h_j,\ldots,h_\ell\in\sigma([m^*+1,n])$, that is, $\{j-1,\ldots,\ell\}\subseteq B$.

\medskip
\noindent Proof of (c): Recall that $[1,\ell]$ is the disjoint union of $A$ and $B$.

Assume by contradiction that $i>i'$ are consecutive elements of $A$ with $i-i'\ge 3$. Then $i-1$ and $i-2$ are not in $A$, so $i-1,i-2\in B$. By (b) (applied with $j=i-1$), $\{i-2,\ldots,\ell\}\subseteq B$. In particular $i\in B$, contradicting $i\in A$.

Similarly, assume by contradiction that $j>j'$ are consecutive elements of $B$ with $j-j'\ge 3$. Then $j-1$ and $j-2$ are not in $B$, so $j-1,j-2\in A$. By (a) (applied with $i=j-1$), $\{1,\ldots,j-1\}\subseteq A$. In particular $j'\in A$, contradicting $j'\in B$.

\medskip
\noindent Proof of (d): Assume $q\ge 2$. Since $a+q=k^*-1$, \eqref{order} (with $s=q-1$ and $s=q$) gives
\[
\sigma(k^*-2)=h_{i_{q-1}}\qquad\text{and}\qquad\sigma(k^*-1)=h_{i_q}.
\]
By \eqref{dist},
\[
|\sigma'(k^*-2)|=h_{i_{q-1}}-h_{i_q}=2(i_{q-1}-i_q)\qquad\text{and}\qquad|\sigma'(k^*-1)|=h_{i_q}-\sigma(k^*)=2i_q-1+r.
\]
Since $|\sigma'|$ is non-increasing on $[1,k^*-1]$, $|\sigma'(k^*-1)|\le|\sigma'(k^*-2)|$. Moreover, $i_{q-1}-i_q\le 2$ by (c). Therefore
\begin{equation}\label{ivineq}
2i_q-1+r\le 2(i_{q-1}-i_q)\le 4 .
\end{equation}
Since $i_q\ge 1$, \eqref{ivineq} gives $r\le 5-2i_q\le 3$.

Assume now $r\ge 2$. Then \eqref{ivineq} gives $2i_q\le 5-r\le 3$, hence $i_q=1$. Using $i_q=1$ in \eqref{ivineq}, we obtain $2(i_{q-1}-1)\ge 1+r\ge 3$, hence $i_{q-1}\ge 3$. Since $i_{q-1}-i_q\le 2$, we conclude that $i_{q-1}=3$.

\medskip
\noindent Proof of (e): assume $t\ge 2$. By \eqref{order2} (with $s=1$ and $s=2$),
\[
\sigma(m^*+1)=h_{j_1}\qquad\text{and}\qquad\sigma(m^*+2)=h_{j_2}.
\]
By \eqref{dist},
\[
|\sigma'(m^*)|=\sigma(m^*)-h_{j_1}=2(\ell-j_1)+p+1\qquad\text{and}\qquad|\sigma'(m^*+1)|=h_{j_1}-h_{j_2}=2(j_1-j_2).
\]
Since $|\sigma'|$ is non-decreasing on $[m^*,n-1]$, $|\sigma'(m^*)|\le|\sigma'(m^*+1)|$. Moreover, $j_1-j_2\le 2$ by (c). Therefore
\begin{equation}\label{vineq}
2(\ell-j_1)+p+1\le 2(j_1-j_2)\le 4 .
\end{equation}
Since $j_1\le\ell$, \eqref{vineq} gives $p+1\le 4$, that is, $p\le 3$. If $j_1-j_2=1$, \eqref{vineq} gives $p+1\le 2$, that is, $p\le 1$.

Assume now $p\ge 2$. Then the left-hand side of \eqref{vineq} is at least $3$, so $2(j_1-j_2)\ge 3$, hence $j_1-j_2=2$. Using this in \eqref{vineq}, we obtain $2(\ell-j_1)\le 3-p\le 1$, hence $j_1=\ell$. It remains to show that any two consecutive elements of $B$ differ by exactly $2$. Since $|\sigma'|$ is non-decreasing on $[m^*,n-1]$,
\[
|\sigma'(u)|\ge|\sigma'(m^*+1)|=2(j_1-j_2)=4\qquad\text{for all } m^*+1\le u\le n-1.
\]
For $1\le v\le t-1$, \eqref{order2} gives $\sigma(m^*+v)=h_{j_v}$ and $\sigma(m^*+v+1)=h_{j_{v+1}}$, so $2(j_v-j_{v+1})=|\sigma'(m^*+v)|\ge 4$, that is, $j_v-j_{v+1}\ge 2$. By (c), $j_v-j_{v+1}=2$. Since $j_1=\ell$, we conclude that $B=\{\ell,\ell-2,\ell-4,\ldots,\ell-2(t-1)\}$.
\end{proof}

To construct permutations with slope $d=2$ we will use Lemma \ref{L:build} below. The idea is to first write the increasing part of $\sigma$ (blocks (B3)-(B5)): a progression with $d=2$, preceded by $r$ steps of $1$ and followed by $p$ steps of $1$. The $\ell$ values $h_1,\ldots,h_\ell$ skipped by the progression are then split into two sets, which are written in decreasing order before and after the increasing part (blocks (B2) and (B6)), possibly together with $n$ in the first position and $1$ in the last one (blocks (B1) and (B7)). By Lemma \ref{L:itisinflection}, the only condition left to check is that the drops $\alpha_s$ at the beginning get smaller, and the drops $\beta_s$ at the end get larger, see Figure \ref{F:build}.

\begin{lem}\label{L:build}
Let $k\ge 2$, $\ell\ge 1$ and $m=k+\ell$. Let $a,c\in\{0,1\}$ and $r,p\ge 0$ be integers, and let $[1,\ell]$ be the disjoint union of $A=\{i_1>\ldots>i_q\}$ and $B=\{j_1>\ldots>j_t\}$, where $t\ge 1$ and $q+a+r=k-1$. Let
\[
x=c+r+1,\qquad h_i=x+2i-1\ \text{ for } 1\le i\le\ell,\qquad n=x+2\ell+p+a .
\]
Let $\sigma$ be obtained by listing, one after the other, the following seven blocks:
\begin{enumerate}
\item[(B1)] the single entry $n$ if $a=1$, and nothing if $a=0$,
\item[(B2)] $h_{i_1}\,h_{i_2}\ldots h_{i_q}$, which decreases,
\item[(B3)] $(c+1)(c+2)\ldots x$, which has $r+1$ entries and increases by $1$,
\item[(B4)] $(x+2)(x+4)\ldots(x+2\ell)$, which has $\ell$ entries and increases by $2$,
\item[(B5)] $(x+2\ell+1)(x+2\ell+2)\ldots(x+2\ell+p)$, which has $p$ entries and increases by $1$,
\item[(B6)] $h_{j_1}\,h_{j_2}\ldots h_{j_t}$, which decreases,
\item[(B7)] the single entry $1$ if $c=1$, and nothing if $c=0$.
\end{enumerate}
Then $\sigma\in\Sym_n$, with $\sigma(k)=x$ and $\sigma(m)=x+2\ell$, and $n=2\ell+r+a+p+c+1$. Moreover, the numbers
\begin{align*}
\alpha_s&=\sigma(s)-\sigma(s+1) && \text{for } 1\le s\le q+a,\\
\beta_s&=\sigma(m+p+s-1)-\sigma(m+p+s) && \text{for } 1\le s\le t+c,
\end{align*}
are positive, and 
\begin{equation}\label{buildderiv}
\sigma'=\underbrace{(-\alpha_1)\ldots(-\alpha_{q+a})}_{q+a}\;\underbrace{1\ldots1}_{r}\;\underbrace{2\ldots2}_{\ell}\;\underbrace{1\ldots1}_{p}\;\underbrace{(-\beta_1)\ldots(-\beta_{t+c})}_{t+c}.
\end{equation}
If
\begin{equation}\label{need}
\alpha_1\ge\alpha_2\ge\ldots\ge\alpha_{q+a}\qquad\text{and}\qquad\beta_1\le\beta_2\le\ldots\le\beta_{t+c},
\end{equation}
then $\mathcal{I}(\sigma)=\{(k,m)\}$.
\end{lem}

\begin{ex}\label{e:lemma511}
Let $\ell=5$, $a=c=1$, $r=2$, $p=1$, $A=\{3,1\}$ and $B=\{5,4,2\}$, so that $k=6$, $m=11$, $x=4$, $h_i=2i+3$ and $n=16$, see Figure \ref{F:build}. Then
\[
\sigma=\underbrace{16}_{\text{(B1)}}\;\;\underbrace{9\;\,5}_{\text{(B2)}}\;\;\underbrace{2\;\,3\;\,4}_{\text{(B3)}}\;\;\underbrace{6\;\,8\;\,10\;\,12\;\,14}_{\text{(B4)}}\;\;\underbrace{15}_{\text{(B5)}}\;\;\underbrace{13\;\,11\;\,7}_{\text{(B6)}}\;\;\underbrace{1}_{\text{(B7)}},
\]
with $(\alpha_1,\alpha_2,\alpha_3)=(7,4,3)$ and $(\beta_1,\beta_2,\beta_3,\beta_4)=(2,2,4,6)$. Since \eqref{need} holds, $\mathcal{I}(\sigma)=\{(6,11)\}$.

\begin{figure}[ht]
\centering
\begin{tikzpicture}[x=0.72cm,y=0.42cm]
  \draw[gray!60] (0.5,0.5) -- (16.6,0.5);
  \draw[gray!60] (0.5,0.5) -- (0.5,16.8);
  \foreach \i in {1,...,16} \node[below,font=\scriptsize,gray!60!black] at (\i,0.5) {\i};
  \foreach \v in {1,2,4,14,15,16} \node[left,font=\scriptsize,gray!60!black] at (0.5,\v) {\v};
  \draw[densely dotted,gray] (6,0.5) -- (6,17);
  \draw[densely dotted,gray] (11,0.5) -- (11,17);
  \node[font=\scriptsize,above] at (6,17) {$k=6$};
  \node[font=\scriptsize,above] at (11,17) {$m=11$};
  \draw[decorate,decoration={brace,amplitude=3pt}] (1.05,19.3) -- (3.95,19.3)
       node[midway,above=4pt,font=\scriptsize,align=center] {$1,\ldots,a+q+1$\\ $\sigma'=-\alpha_s$};
  \draw[decorate,decoration={brace,amplitude=3pt}] (4.05,19.3) -- (11.95,19.3)
       node[midway,above=4pt,font=\scriptsize,align=center] {$a+q+1,\ldots,m+p$\\ $\sigma'$ is $1\ldots1$, then $2\ldots2$, then $1\ldots1$};
  \draw[decorate,decoration={brace,amplitude=3pt}] (12.05,19.3) -- (15.95,19.3)
       node[midway,above=4pt,font=\scriptsize,align=center] {$m+p,\ldots,n$\\ $\sigma'=-\beta_s$};
  \draw[densely dashed,gray] (2,9)  -- (8.5,9);
  \draw[densely dashed,gray] (3,5)  -- (6.5,5);
  \draw[densely dashed,gray] (10.5,13) -- (13,13);
  \draw[densely dashed,gray] (9.5,11)  -- (14,11);
  \draw[densely dashed,gray] (7.5,7)   -- (15,7);
  \draw[thick] (1,16)--(2,9)--(3,5)--(4,2)--(5,3)--(6,4)--(7,6)--(8,8)--(9,10)--(10,12)--(11,14)--(12,15)--(13,13)--(14,11)--(15,7)--(16,1);
  \foreach \p/\v in {4/2,5/3,6/4,7/6,8/8,9/10,10/12,11/14,12/15} \fill (\p,\v) circle (2.2pt);
  \foreach \p/\v in {2/9,3/5} \node[draw,fill=white,inner sep=2.2pt,thick] at (\p,\v) {};
  \foreach \p/\v in {13/13,14/11,15/7} \node[draw,fill=white,regular polygon,regular polygon sides=3,inner sep=1.6pt,thick] at (\p,\v) {};
  \foreach \p/\v in {1/16,16/1} \node[draw,fill=white,circle,inner sep=2pt,thick] at (\p,\v) {};
  \foreach \q/\v/\n in {6.5/5/1,7.5/7/2,8.5/9/3,9.5/11/4,10.5/13/5}
     { \node[font=\tiny,inner sep=0.5pt] at (\q,\v) {$\times$};
       \node[font=\scriptsize,anchor=south east,inner sep=1pt] at (\q-0.05,\v+0.1) {$h_{\n}$}; }
  \node[font=\scriptsize,anchor=west] at (1.55,12.6) {$\alpha_1=7$};
  \node[font=\scriptsize,anchor=west] at (2.55,7.1)  {$\alpha_2=4$};
  \node[font=\scriptsize,anchor=east] at (3.45,3.1)  {$\alpha_3=3$};
  \node[font=\scriptsize,anchor=west] at (12.55,14.3) {$\beta_1=2$};
  \node[font=\scriptsize,anchor=west] at (13.55,12.3) {$\beta_2=2$};
  \node[font=\scriptsize,anchor=west] at (14.55,9.2)  {$\beta_3=4$};
  \node[font=\scriptsize,anchor=east] at (15.45,3.6)  {$\beta_4=6$};
  \foreach \l/\r/\n in {0.75/1.25/1, 1.75/3.25/2, 3.75/6.25/3, 6.75/11.25/4, 11.75/12.25/5, 12.75/15.25/6, 15.75/16.25/7}
     \draw[decorate,decoration={brace,mirror,amplitude=3pt}] (\l,-1.0) -- (\r,-1.0)
          node[midway,below=4pt,font=\scriptsize] {(B\n)};
\end{tikzpicture}
\caption{The permutation of Example \ref{e:lemma511}, $\sigma=16\;9\;5\;2\;3\;4\;6\;8\;10\;12\;14\;15\;13\;11\;7\;1$.}
\label{F:build}
\end{figure}
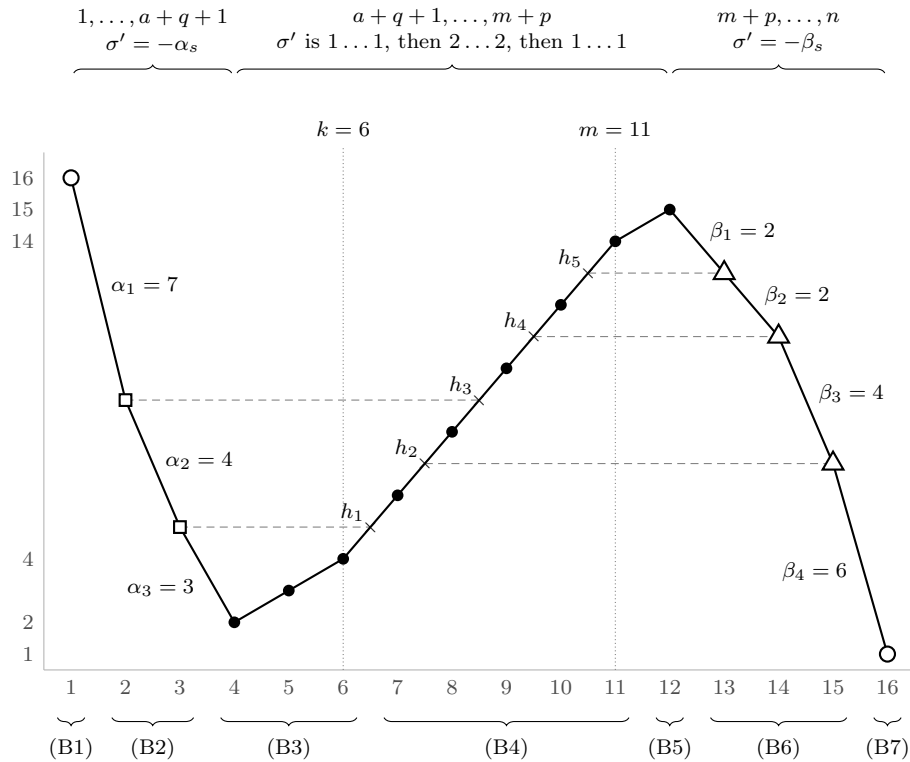

\end{ex}

\begin{proof}
We first show that $\sigma\in\Sym_n$. The entries of (B3) are the elements of $\{c+1,\ldots,x\}$, and the entries of (B5) are the elements of $\{x+2\ell+1,\ldots,x+2\ell+p\}$. The entries of (B4), (B2) and (B6) are the numbers $x+2i$ and $h_i=x+2i-1$ for $1\le i\le\ell$, which together with $x$ are the elements of $\{x,\ldots,x+2\ell\}$. Finally, (B7) contributes $1$ exactly when $c=1$, that is, when $c+1>1$, and (B1) contributes $n$ exactly when $a=1$, that is, when $n>x+2\ell+p$. Hence every element of $[1,n]$ appears exactly once in $\sigma$.

Next we locate the blocks. Blocks (B1) and (B2) have $a+q$ entries, so (B3) occupies the positions $a+q+1,\ldots,a+q+r+1=k$, and (B4) occupies  positions $k+1,\ldots,k+\ell=m$. In particular, $\sigma(k)=x$ and $\sigma(m)=x+2\ell$. Block (B5) occupies the positions $m+1,\ldots,m+p$, block (B6) occupies the positions $m+p+1,\ldots,m+p+t$, and $n=m+p+t+c$. Since $t\ge 1$, we have $m\le n-1$. The formula $n=2\ell+r+a+p+c+1$ follows from $n=x+2\ell+p+a$ and $x=c+r+1$.

We now compute $\sigma'$ and prove \eqref{buildderiv}. We split the positions of $\sigma$ into three ranges.

\begin{itemize}
\item Positions $1,\ldots,a+q+1$: here, the entries of $\sigma$ are, in order, $n$ (only if $a=1$), $h_{i_1},\,h_{i_2},\ldots,h_{i_q}$ (only if $q\ge 1$), and $c+1$.

This list is strictly decreasing. Indeed, if $a=1$ then $n=x+2\ell+p+1>h_\ell$, so $n$ is larger than every $h_i$ and than $c+1$. Moreover, $h_{i_1}>h_{i_2}>\ldots>h_{i_q}$ because $i_1>i_2>\ldots>i_q$. Finally, every $h_i$ satisfies $h_i\ge h_1=x+1>c+1$. (If $q=a=0$, this range is the single position $1$ and there are no $\alpha_s$.) Therefore, for $1\le s\le q+a$,
\[
\alpha_s=\sigma(s)-\sigma(s+1)>0\qquad\text{and}\qquad\sigma'(s)=-\alpha_s .
\]
\item Positions $a+q+1,\ldots,m+p$. Here, the entries of $\sigma$ are those of the blocks
(B3), (B4) and (B5), that is, in order,
\[
c+1,\,c+2,\ldots,x,\qquad x+2,\,x+4,\ldots,x+2\ell,\qquad x+2\ell+1,\,x+2\ell+2,\ldots,x+2\ell+p .
\]
Two consecutive entries of this list differ by $1$ (the first $r$ times), then by $2$ (the next
$\ell$ times), then by $1$ (the last $p$ times). Since $a+q+r+1=k$ and $k+\ell=m$, this means that
\[
\sigma'(s)=
\begin{cases}
1 & \text{for } a+q+1\le s\le k-1,\\
2 & \text{for } k\le s\le m-1,\\
1 & \text{for } m\le s\le m+p-1 .
\end{cases}
\]

\item Positions $m+p,\ldots,n$. Here, the entries of $\sigma$ are, in order,
\[
x+2\ell+p,\qquad h_{j_1},\,h_{j_2},\ldots,h_{j_t},\qquad 1\ \ (\text{only if } c=1).
\]
This list is strictly decreasing: $x+2\ell+p>h_\ell\ge h_{j_1}$, $h_{j_1}>h_{j_2}>\ldots>h_{j_t}$
because $j_1>j_2>\ldots>j_t$, and $h_{j_t}\ge h_1=x+1>1$. Therefore, for $1\le s\le t+c$,
\[
\beta_s=\sigma(m+p+s-1)-\sigma(m+p+s)>0\qquad\text{and}\qquad\sigma'(m+p+s-1)=-\beta_s .
\]
\end{itemize}
Reading the three ranges in order, we obtain \eqref{buildderiv}.

Finally, assume \eqref{need}. We check the three conditions of Lemma \ref{L:itisinflection}. By \eqref{buildderiv}, $\sigma'(j)=2$ for $k\le j\le m-1$, so that $d=\sigma(k+1)-\sigma(k)=2>0$ and \eqref{derg} holds. On $[1,k-1]$, $\sigma'$ is $-\alpha_1\le-\alpha_2\le\ldots\le-\alpha_{q+a}<0$, followed by $r$ entries equal to $1$, so \eqref{mono1} holds. On $[m,n-1]$, $\sigma'$ is $p$ entries equal to $1$ followed by $0>-\beta_1\ge-\beta_2\ge\ldots\ge-\beta_{t+c}$, so \eqref{mono2} holds. By Lemma \ref{L:itisinflection}, $\mathcal{I}(\sigma)=\{(k,m)\}$.
\end{proof}

\begin{figure}[ht]
\centering
    \begin{tikzpicture}[x=0.8cm,y=0.42cm]
  \fill[gray!15] (0.5,3.5) rectangle (15.6,13.5);
  \node[font=\scriptsize,anchor=west] at (2.4,12.4) {the $2\ell$ values $\sigma(k)+1,\ldots,\sigma(m)$};
  \draw[gray!60] (0.5,0.5) -- (15.6,0.5);
  \draw[gray!60] (0.5,0.5) -- (0.5,15.8);
  \foreach \i in {1,...,15} \node[below,font=\scriptsize,gray!60!black] at (\i,0.5) {\i};
  \foreach \v in {1,...,15} \node[left,font=\scriptsize,gray!60!black] at (0.5,\v) {\v};
  \draw[densely dotted,gray] (4,0.5) -- (4,16);
  \draw[densely dotted,gray] (9,0.5) -- (9,16);
  \node[font=\scriptsize,above] at (4,16) {$k$};
  \node[font=\scriptsize,above] at (9,16) {$m$};
  \draw[densely dashed,gray] (2,4)   -- (4.5,4);
  \draw[densely dashed,gray] (5.5,6)  -- (14,6);
  \draw[densely dashed,gray] (6.5,8)  -- (13,8);
  \draw[densely dashed,gray] (7.5,10) -- (12,10);
  \draw[densely dashed,gray] (8.5,12) -- (11,12);
  \draw[thick] (1,15)--(2,4)--(3,2)--(4,3)--(5,5)--(6,7)--(7,9)--(8,11)--(9,13)--(10,14)--(11,12)--(12,10)--(13,8)--(14,6)--(15,1);
  \foreach \p/\v in {5/5,6/7,7/9,8/11,9/13} \fill (\p,\v) circle (2.2pt);
  \foreach \q/\v in {4.5/4,5.5/6,6.5/8,7.5/10,8.5/12} \node[font=\tiny,inner sep=0.5pt] at (\q,\v) {$\times$};
  \node[draw,fill=white,inner sep=2.2pt,thick] at (2,4) {};
  \foreach \p/\v in {11/12,12/10,13/8,14/6} \node[draw,fill=white,regular polygon,regular polygon sides=3,inner sep=1.6pt,thick] at (\p,\v) {};
  \foreach \p/\v in {1/15,3/2,4/3,10/14,15/1}
     { \fill (\p,\v) circle (2.2pt); \draw[thick] (\p,\v) circle (4.6pt); }
  \node[font=\scriptsize,anchor=west] at (1.3,15) {$n$};
  \node[font=\scriptsize,anchor=north west] at (4.25,2.9) {$\sigma(k)$};
  \node[font=\scriptsize,anchor=south east] at (8.85,13.2) {$\sigma(m)$};
  \node[font=\scriptsize,anchor=east] at (14.7,1) {$1$};
  \foreach \l/\r/\t in {0.75/1.25/a, 1.75/2.25/q, 2.75/4.25/r+1, 4.75/9.25/\ell, 9.75/10.25/p, 10.75/14.25/\ell-q, 14.75/15.25/c}
     \draw[decorate,decoration={brace,mirror,amplitude=3pt}] (\l,-1.0) -- (\r,-1.0)
          node[midway,below=4pt,font=\scriptsize] {$\t$};
\end{tikzpicture}
\caption{The idea of the proof of Theorem \ref{T:windows}, for $\sigma=15\;4\;2\;3\;5\;7\;9\;11\;13\;14\;12\;10\;8\;6\;1$, where $(k,m,n)=(4,9,15)$ and $\ell=5$. The shaded band contains $2\ell$ values: the $\ell$ entries $\sigma(k+1),\ldots,\sigma(m)$ of the progression and the $\ell$ values that it skips. Each skipped value is used before the progression ($\square$) or after it ($\triangle$). The remaining $N=n-2\ell=r+a+p+c+1$ entries are circled. Here $q=1$, $r=a=p=c=1$ and $N=5$.}
\label{F:idea}
\end{figure}
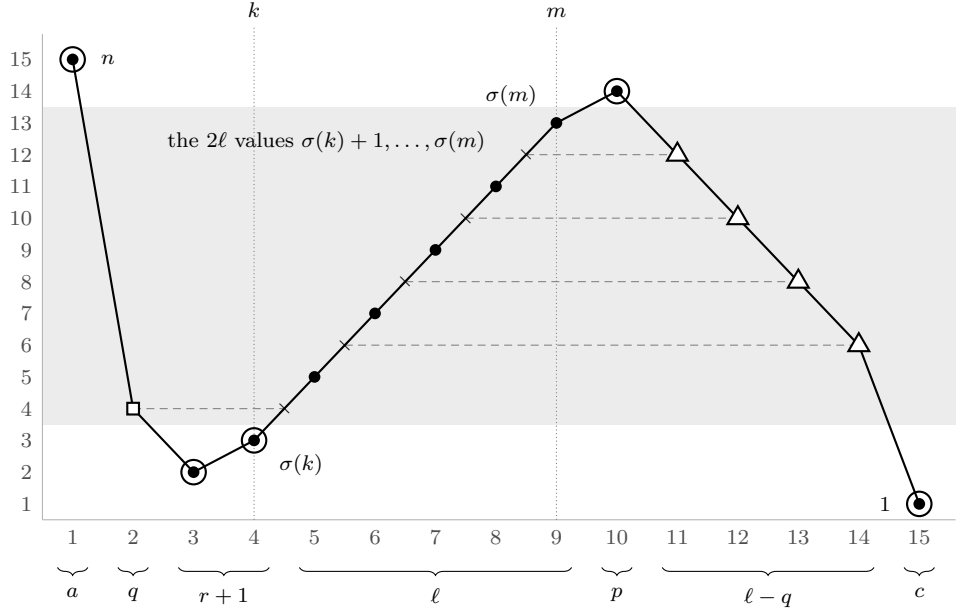

We can now state and prove the first main result of this section, which fully characterizes, for $m\ge 2k+1$, when $\{\sigma\in \Sym_n \mid \mathcal{I}(\sigma)=\{(k,m)\}\}\neq\emptyset$. The idea of the proof is as follows (see Figure \ref{F:idea}): outside of case (5) below, the progression $\sigma(k),\ldots,\sigma(m)$ has $d=2$, so it skips $\ell=m-k$ values, and each of them must be used before or after the progression. Hence $n=2\ell+N$, where $N$ counts the entries that are neither in the progression after $\sigma(k)$ nor skipped by it, and the monotonicity of $\sigma'$ before $k$ and after $m$ leaves space for only a few such entries. 

\begin{thm}\label{T:windows}
Let $k\ge 2$, $m\ge 2k+1$ and $n\ge m+1$, and write $n=2(m-k)+N$. Then $|\{\sigma\in \Sym_n \mid \mathcal{I}(\sigma)=\{(k,m)\}\}|>0$ if and only if one of the following holds:
\begin{enumerate}
\item[(1)] $1\le N\le \min\{7,k+2\}$;
\item[(2)] $k-1\le N\le k+2$;
\item[(3)] $N=8$ and $m\in\{3k-10,\,3k-8\}$;
\item[(4)] $N=9$ and $m=3k-10$;
\item[(5)] $k=2$ and $n=m+1$.
\end{enumerate}
\end{thm}

\begin{proof}
Write $\ell=m-k$, so that $\ell\ge k+1\ge 3$.

We first prove that the conditions are necessary. Let $\sigma\in\Sym_n$ with $\mathcal{I}(\sigma)=\{(k,m)\}$. By Lemma \ref{L:sign}, we may assume $d=\sigma(k+1)-\sigma(k)>0$, so that $\sigma$ is as in Lemma \ref{L:shape}. Since $m>2k$, Lemma \ref{L:d3} gives $d\le 2$. If $d=1$, Lemma \ref{L:d1} gives (5). Assume from now on that $d=2$, and let $h_i$, $A$, $B$, $q$, $t$, $a$, $c$, $r$, $p$ be as in Lemma \ref{L:d2}. By \eqref{count1} and \eqref{count2},
\begin{equation}\label{Nformula}
N=n-2\ell=k-q+p+c=r+a+p+c+1 .
\end{equation}
In particular, $N\ge 1$. Moreover, \eqref{count1} gives $q\le k-1$, hence $t=\ell-q\ge 2$, and Lemma \ref{L:d2}(e) gives $p\le 3$. We divide the proof in cases, depending on $q$.
\begin{enumerate}
\item[\parbox{1.2cm}{Case \\ $q=0$:}] In this case $B=[1,\ell]$, so $j_1-j_2=1$ and Lemma \ref{L:d2}(e) gives $p\le 1$. By \eqref{Nformula}, $N=k+p+c$, hence $k\le N\le k+2$ and (2) holds.

\item[\parbox{1.2cm}{Case \\ $q=1$:}] By \eqref{Nformula}, $N=k-1+p+c$, hence $k-1\le N\le k+3$, and $N=k+3$ only if $p=3$ and $c=1$. Assume by contradiction that $p=3$ and $c=1$. By Lemma \ref{L:d2}(e), $B=\{\ell,\ell-2,\ldots,\ell-2(t-1)\}$ with $t=\ell-1$. Since $\ell-2(t-1)\ge 1$, we obtain $\ell\le 3$, hence $\ell=3$, $k=2$ and $B=\{3,1\}$. By \eqref{count1}, $r=a=0$. By \eqref{order2}, $\sigma(n-2)=h_3$, $\sigma(n-1)=h_1$ and $\sigma(n)=1$, and by \eqref{dist}, $|\sigma'(n-2)|=4$ and $|\sigma'(n-1)|=h_1-1=2$. This contradicts the fact that $|\sigma'|$ is non-decreasing on $[m^*,n-1]$. Therefore $k-1\le N\le k+2$ and (2) holds.

\item[\parbox{1.2cm}{Case \\ $q\ge 2$:}] By \eqref{Nformula}, $N=k-q+p+c\le k-2+3+1=k+2$. By Lemma \ref{L:d2}(d), $r\le 3$, so \eqref{Nformula} also gives $N=r+a+p+c+1\le 9$. If $N\le 7$, then (1) holds. Assume now $N\ge 8$. Then $r+p\ge 5$, and since $r\le 3$ and $p\le 3$, we must have $r\ge 2$ and $p\ge 2$. By Lemma \ref{L:d2}(d), $1\in A$, $3\in A$ and $2\notin A$, that is, $2\in B$. By Lemma \ref{L:d2}(e), $B=\{\ell,\ell-2,\ldots,\ell-2(t-1)\}$. Since $2\in B$, $\ell$ is even and $B=\{2,4,\ldots,\ell\}$, so $q=\ell/2$. By \eqref{count1},
\[
\ell=2q=2(k-1-r-a).
\]
If $N=9$, then $r=3$, $a=1$, $p=3$ and $c=1$, hence $\ell=2k-10$, that is, $m=3k-10$, and (4) holds. If $N=8$, then $(r,a,p,c)$ is one of $(3,1,3,0)$, $(3,1,2,1)$, $(3,0,3,1)$, $(2,1,3,1)$, hence $\ell\in\{2k-10,2k-8\}$, that is, $m\in\{3k-10,3k-8\}$, and (3) holds.
\end{enumerate}
This proves that one of (1)--(5) holds.

We now prove that the conditions are sufficient. If (5) holds, then $\sigma=(m+1)\,2\,3\ldots m\,1\in\Sym_{m+1}$ has first derivative
\[
(1-m)\underbrace{1\ldots 1}_{m-2}(1-m),
\]
so that $\mathcal{I}(\sigma)=\{(2,m)\}$ by Lemma \ref{L:itisinflection} (with $d=1$). 

For the remaining cases, we use Lemma \ref{L:build}, with three different choices of sets $A$ and $B$. For each choice, we say which missing values go in the blocks (B2) and (B6), and we compute the numbers $\alpha_s$ and $\beta_s$ of that lemma, to find out when \eqref{need} holds. Throughout, $x=c+r+1$, $h_i=x+2i-1$ and $n=x+2\ell+p+a$, so that
\[
h_i-h_{i'}=2(i-i'),\qquad h_i-(c+1)=2i-1+r,\qquad (x+2\ell+p)-h_i=2(\ell-i)+p+1,
\]
and moreover $n-h_i=2(\ell-i)+p+2$ if $a=1$, and $h_i-1=2i+r$ if $c=1$.

\medskip
\noindent First choice: $A=\{1,\ldots,q\}$ and $B=\{q+1,\ldots,\ell\}$, with $0\le q\le k-1$.

The block (B2) is $h_q\,h_{q-1}\ldots h_1$, and the block (B6) is $h_\ell\,h_{\ell-1}\ldots h_{q+1}$. If $q\ge 1$, then
\[
(\alpha_1,\ldots,\alpha_{q+a})=
\begin{cases}
\big(\underbrace{2,\ldots,2}_{q-1},\ 1+r\big) & \text{if } a=0,\\[14pt]
\big(2(\ell-q)+p+2,\ \underbrace{2,\ldots,2}_{q-1},\ 1+r\big) & \text{if } a=1.
\end{cases}
\]
If $q=0$, there are no $\alpha_s$ when $a=0$, and there is only $\alpha_1=n-(c+1)$ when $a=1$. In all cases,
\[
(\beta_1,\ldots,\beta_{t+c})=
\begin{cases}
\big(p+1,\ \underbrace{2,\ldots,2}_{\ell-q-1}\big) & \text{if } c=0,\\[14pt]
\big(p+1,\ \underbrace{2,\ldots,2}_{\ell-q-1},\ 2q+2+r\big) & \text{if } c=1.
\end{cases}
\]
Notice that $\ell-q-1\ge 1$, since $q\le k-1\le\ell-2$. Hence the second part of \eqref{need} holds if and only if $p\le 1$. The first part of \eqref{need} holds in each of the following situations: if $q=0$, or if $q=1$ (when $a=1$ it requires $1+r\le 2(\ell-1)+p+2$, which holds since $r\le k-2\le\ell-3$), and if $q\ge 2$ and $r\le 1$. In summary, with this choice, \eqref{need} holds if
\[
p\le 1,\qquad\text{and}\qquad r\le 1\ \text{ in case } q\ge 2 .
\]

\medskip
\noindent Second choice: $A=\{1,3,5,\ldots,2q-1\}$ and $B=\{2,4,\ldots,2q\}\cup\{2q+1,\ldots,\ell\}$, with $q\ge 2$ and $\ell\ge 2q$.

The block (B2) is $h_{2q-1}\,h_{2q-3}\ldots h_3\,h_1$, and the block (B6) is $h_\ell\,h_{\ell-1}\ldots h_{2q+1}\,h_{2q}\,h_{2q-2}\ldots h_4\,h_2$ (Example \ref{e:lemma511} is of this type, with $q=2$). We have
\[
(\alpha_1,\ldots,\alpha_{q+a})=
\begin{cases}
\big(\underbrace{4,\ldots,4}_{q-1},\ 1+r\big) & \text{if } a=0,\\[14pt]
\big(2(\ell-2q)+p+4,\ \underbrace{4,\ldots,4}_{q-1},\ 1+r\big) & \text{if } a=1,
\end{cases}
\]
and
\[
(\beta_1,\ldots,\beta_{t+c})=
\begin{cases}
\big(p+1,\ \underbrace{2,\ldots,2}_{\ell-2q},\ \underbrace{4,\ldots,4}_{q-1}\big) & \text{if } c=0,\\[14pt]
\big(p+1,\ \underbrace{2,\ldots,2}_{\ell-2q},\ \underbrace{4,\ldots,4}_{q-1},\ 4+r\big) & \text{if } c=1.
\end{cases}
\]
The first part of \eqref{need} holds if and only if $1+r\le 4$. The second part holds if and only if $p+1\le 2$ (when $\ell>2q$) or $p+1\le 4$ (when $\ell=2q$). In summary, with this choice, \eqref{need} holds if
\[
r\le 3,\qquad\text{and}\qquad p\le 1,\ \text{ or }\ p\le 3\ \text{ and }\ \ell=2q .
\]

\medskip
\noindent Third choice: $B=\{j,j+2,j+4,\ldots,\ell\}$ and $A=\{1,\ldots,j-1\}\cup\{j+1,j+3,\ldots,\ell-1\}$, where $j=2q+2-\ell$ and $q+2\le\ell\le 2q$.

Notice that $j\ge 2$, that $B$ has $\ell-q$ elements, and that $A$ has $q$ elements. The block (B2) is $h_{\ell-1}\,h_{\ell-3}\ldots h_{j+1}\,h_{j-1}\,h_{j-2}\ldots h_1$, and the block (B6) is $h_\ell\,h_{\ell-2}\ldots h_{j+2}\,h_j$. We have
\[
(\alpha_1,\ldots,\alpha_{q+a})=
\begin{cases}
\big(\underbrace{4,\ldots,4}_{\ell-q-1},\ \underbrace{2,\ldots,2}_{j-2},\ 1+r\big) & \text{if } a=0,\\[14pt]
\big(p+4,\ \underbrace{4,\ldots,4}_{\ell-q-1},\ \underbrace{2,\ldots,2}_{j-2},\ 1+r\big) & \text{if } a=1,
\end{cases}
\]
and
\[
(\beta_1,\ldots,\beta_{t+c})=
\begin{cases}
\big(p+1,\ \underbrace{4,\ldots,4}_{\ell-q-1}\big) & \text{if } c=0,\\[14pt]
\big(p+1,\ \underbrace{4,\ldots,4}_{\ell-q-1},\ 2j+r\big) & \text{if } c=1.
\end{cases}
\]
Notice that $\ell-q-1\ge 1$ and $2j+r\ge 4$. In summary, with this choice \eqref{need} holds if
\[
r\le 1\qquad\text{and}\qquad p\le 3 .
\]

\medskip
In all three choices, Lemma \ref{L:build} gives a permutation $\sigma\in\Sym_n$ with $\mathcal{I}(\sigma)=\{(k,m)\}$ and $n=2\ell+N$, where $N=r+a+p+c+1$ and $q=k-1-r-a$. We now go through the values of $N$ in (1)-(4).
\begin{itemize}
\item If $N=k-1$, we use the first choice with $q=1$, $a=0$, $r=k-2$ and $p=c=0$.
\item If $N\in\{k,k+1,k+2\}$, we use the first choice with $q=0$, $a=0$, $r=k-1$, and $p,c\in\{0,1\}$ with $p+c=N-k$.
\item If $1\le N\le\min\{5,k+2\}$, we use the first choice with $r,a,p,c\in\{0,1\}$ and $r+a+p+c=N-1$, where we take $r=a=0$ if $N\le 3$, and $r=1$, $a=0$ if $N=4$. Then $q=k-1-r-a\ge 0$ (notice that if $N=5$, then $k\ge 3$).
For example, if $N=1$, then $q=k-1$ and $r=a=p=c=0$, and Lemma \ref{L:build} gives $\sigma=(2k-2)(2k-4)\ldots 2\,1\,3\,5\ldots(2\ell+1)(2\ell)(2\ell-2)\ldots(2k)\in\Sym_{2\ell+1}$, with 
\[
\sigma'=\underbrace{(-2)\ldots(-2)}_{k-2}(-1)\underbrace{2\ldots 2}_{\ell}(-1)\underbrace{(-2)\ldots(-2)}_{m-2k}.
\]
\item Let $N\in\{6,7\}$. If $N\ge k-1$, we are in one of the first two items, so assume $N\le k-2$, hence $k\ge 8$. If $\ell\ge 2k-10$, we use the second choice with $r=3$, $a=1$, $p=1$, $c=N-6$ and $q=k-5$. If $\ell\le 2k-11$, we use the third choice with $r=1$, $a=1$, $p=3$, $c=N-6$ and $q=k-3$ (notice that $q+2\le\ell\le 2q$ since $\ell\ge k+1$).
\item If $N\in\{8,9\}$ and $m=3k-10$, then $\ell=2k-10$, and we use the second choice with $q=k-5=\ell/2$, $r=3$, $a=1$, $p=3$ and $c=N-8$.
\item If $N=8$ and $m=3k-8$, then $\ell=2k-8$, and we use the second choice with $q=k-4=\ell/2$, $r=3$, $a=0$, $p=3$ and $c=1$.
\end{itemize}
In the last two items, $q=\ell/2\ge 2$ since $\ell\ge 3$.
\end{proof}

\begin{cor}\label{C:windows}
Let $k\ge 2$ and $m\ge 2k+1$.
\begin{enumerate}
\item[(a)] If $\sigma\in\Sym_n$ and $\mathcal{I}(\sigma)=\{(k,m)\}$, then $2m-2k+1\le n\le 2m-k+2$, or $k=2$ and $n=m+1$. In particular, $|\{\sigma\in \Sym_n \mid \mathcal{I}(\sigma)=\{(k,m)\}\}|=0$ for all $n\ge 2m-k+3$, and this is sharp.
\item[(b)] There exists $\sigma\in\Sym_n$ with $\mathcal{I}(\sigma)=\{(k,m)\}$ for every $2m-2k+1\le n\le 2m-k+2$ if and only if $k\le 9$ or $(k,m)\in\{(10,22),(11,23)\}$.
\item[(c)] If $k\ge 10$ and $n=2m-2k+N$ with $8\le N\le k-2$, there exists no $\sigma\in\Sym_n$ with $\mathcal{I}(\sigma)=\{(k,m)\}$, except when $N=8$ and $m\in\{3k-10,3k-8\}$, or $N=9$ and $m=3k-10$. For instance, there is no $\sigma\in\Sym_{30}$ with $\mathcal{I}(\sigma)=\{(10,21)\}$, although such permutations exist in $\Sym_{29}$ and in $\Sym_{31}$.
\end{enumerate}
\end{cor}

\begin{proof}
Writing $n=2(m-k)+N$, the inequalities $2m-2k+1\le n\le 2m-k+2$ are equivalent to $1\le N\le k+2$. Hence (a) and (c) follow from Theorem \ref{T:windows}. Notice that (a) is sharp since $N=k+2$ satisfies Theorem \ref{T:windows}(2). 

For (b), conditions (1) and (2) from Theorem \ref{T:windows} cover all $1\le N\le k+2$ exactly when $k\le 9$. 
If $k\ge 10$, the values $8\le N\le k-2$ must satisfy (3) or (4) from Theorem \ref{T:windows}. If $k\ge 12$, this fails for $N=10$. If $k=11$, $N=9$ forces $m=3k-10=23$, and then $N=8$ satisfies (3) from Theorem \ref{T:windows}. If $k=10$, $N=8$ forces $m\in\{20,22\}$, and $m\ge 2k+1$ gives $m=22$. Finally, $(k,m,n)=(10,21,30)$ has $N=8$ and $m\notin\{20,22\}$, while $n=29$ and $n=31$ correspond to $N=7$ and $N=9=k-1$.
\end{proof}

\begin{ex}
    Consider the set $\{(3,7)\}$. By Corollary \ref{C:windows}, there exists $\sigma\in\Sym_n$ with $\mathcal{I}(\sigma)=\{(3,7)\}$ if and only if $9\le n\le 13$.
\end{ex}

We now turn to the case $k<m\le 2k$. In Theorem \ref{T:witnesses}, we construct permutations with the single inflection $(k,m)$. The construction itself depends on the relationship between $k,m$ and $n$.

In (a), $m$ is close to $2k$, so the number $\ell=m-k$ of values skipped by a progression with $d=2$ is close to $k-1$, the number of positions before $k$. We fill these positions with skipped values in decreasing order, together with at most two other entries, and after the progression, $\sigma$ simply increases by $1$ for as long as needed. This is why every $n\ge m+1$ works, see Figures \ref{F:m2k}-\ref{F:m2km3}. 

In (b), we use a progression with $d=3$, which skips two values between any two consecutive entries. We use one of them before the progression and the other one after it (blocks $P$ and $U$), so that $\sigma$ decreases by $3$ at the beginning and at the end. The values that are smaller than $\sigma(k)$ or larger than $\sigma(m)$ are placed next to the progression (blocks $Q$ and $S$), in an order that goes down by $2$ and then up by $2$ (respectively, up by $2$ and then down by $2$), which keeps $\sigma'$ monotone. Since $S$ can have any number of entries, this works for every $n\ge 2m-k$, see Figure \ref{F:PQRSU}. 

In (c), we again use $d=2$, but now only one skipped value is used before the progression, and all the others are used after it.

\begin{thm}\label{T:witnesses}
Let $2\le k<m\le 2k$ and write $\ell=m-k$.
\begin{enumerate}
\item[(a)] If $m\ge 2k-3$ and $n\ge m+1$, there exists $\sigma\in\Sym_n$ with $\mathcal{I}(\sigma)=\{(k,m)\}$.
\item[(b)] If $m\le 2k-1$ and $n\ge 2m-k$, there exists $\sigma\in\Sym_n$ with $\mathcal{I}(\sigma)=\{(k,m)\}$.
\item[(c)] If $m\ge k+2$ and $n=2m-k-1$, there exists $\sigma\in\Sym_n$ with $\mathcal{I}(\sigma)=\{(k,m)\}$.
\end{enumerate}
\end{thm}

\begin{proof}
(a) We break into four cases depending on the value of $m$. In all of them, the arithmetic progression $\sigma(k),\ldots,\sigma(m)$ has $d=2$, the first derivative is non-decreasing and at most $1$ on $[1,k-1]$, and non-increasing and at most $1$ on $[m,n-1]$. Hence \eqref{mono1}--\eqref{mono2} hold with $d=2$, and Lemma \ref{L:itisinflection} gives $\mathcal{I}(\sigma)=\{(k,m)\}$.

When $m=2k$, let $\sigma=(2k-2)(2k-4)\ldots 2\,1\,3\,5\ldots(2k+1)(2k+2)\ldots n\,(2k)$ (see Figure \ref{F:m2k}), where $\sigma(k)=1$, $\sigma(m)=2k+1$, and 
\[
\sigma'=\underbrace{(-2)\ldots(-2)}_{k-2}(-1)\underbrace{2\ldots2}_{k}\underbrace{1\ldots1}_{n-2k-1}(2k-n).
\]
Notice that $2k-n\le -1$ since $n\ge m+1=2k+1$.

\begin{figure}[ht]
\centering
\begin{tikzpicture}[x=0.8cm,y=0.42cm]
  \draw[gray!60] (0.5,0.5) -- (12.6,0.5);
  \draw[gray!60] (0.5,0.5) -- (0.5,12.8);
  \foreach \i in {1,...,12} \node[below,font=\scriptsize,gray!60!black] at (\i,0.5) {\i};
  \foreach \v in {1,...,12} \node[left,font=\scriptsize,gray!60!black] at (0.5,\v) {\v};
  \draw[densely dotted,gray] (4,0.5) -- (4,13);
  \draw[densely dotted,gray] (8,0.5) -- (8,13);
  \node[font=\scriptsize,above] at (4,13) {$k$};
  \node[font=\scriptsize,above] at (8,13) {$m=2k$};
  \draw[thick] (1,6)--(2,4)--(3,2)--(4,1)--(5,3)--(6,5)--(7,7)--(8,9)--(9,10)--(10,11)--(11,12)--(12,8);
  \foreach \p/\v in {1/6,2/4,3/2,4/1,5/3,6/5,7/7,8/9,9/10,10/11,11/12,12/8} \fill (\p,\v) circle (2.2pt);
  \node[font=\scriptsize,anchor=south west] at (1.05,6.1) {$2k-2$};
  \node[font=\scriptsize,anchor=west] at (4.25,1.1) {$\sigma(k)=1$};
  \node[font=\scriptsize,anchor=south east] at (7.9,9.1) {$\sigma(m)=2k+1$};
  \node[font=\scriptsize,anchor=south] at (11,12.2) {$n$};
  \node[font=\scriptsize,anchor=west] at (12.15,8) {$2k$};
  \draw[decorate,decoration={brace,mirror,amplitude=3pt}] (1.05,-1.0) -- (2.95,-1.0)
       node[midway,below=4pt,font=\scriptsize,align=center] {$k-2$ steps\\ of $-2$};
  \draw[decorate,decoration={brace,mirror,amplitude=3pt}] (3.05,-3.6) -- (3.95,-3.6)
       node[midway,below=4pt,font=\scriptsize,align=center] {one step\\ of $-1$};
  \draw[decorate,decoration={brace,mirror,amplitude=3pt}] (4.05,-1.0) -- (7.95,-1.0)
       node[midway,below=4pt,font=\scriptsize,align=center] {$k$ steps\\ of $2$};
  \draw[decorate,decoration={brace,mirror,amplitude=3pt}] (8.05,-1.0) -- (10.95,-1.0)
       node[midway,below=4pt,font=\scriptsize,align=center] {$n-2k-1$ steps\\ of $1$};
  \draw[decorate,decoration={brace,mirror,amplitude=3pt}] (11.05,-3.6) -- (11.95,-3.6)
       node[midway,below=4pt,font=\scriptsize,align=center] {one step\\ of $2k-n$};
\end{tikzpicture}
\caption{$\sigma=(2k-2)(2k-4)\ldots 2\,1\,3\,5\ldots(2k+1)(2k+2)\ldots n\,(2k)$ with $k=4$ and $n=12$, that is, $\sigma=6\;4\;2\;1\;3\;5\;7\;9\;10\;11\;12\;8$.}
\label{F:m2k}
\end{figure}
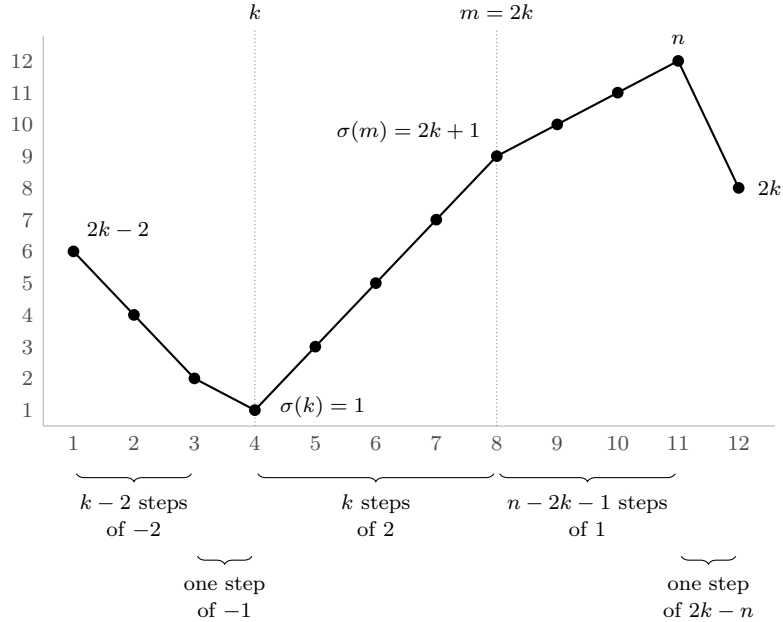

When $m=2k-1$, let $\sigma=(2k-2)(2k-4)\ldots 2\,1\,3\,5\ldots(2k-1)(2k)(2k+1)\ldots n$, see Figure \ref{F:m2km1}. Then $\sigma\in\Sym_n$, with $\sigma(k)=1$ and $\sigma(m)=2k-1$, and $\sigma$ has first derivative
\[
\underbrace{(-2)\ldots(-2)}_{k-2}(-1)\underbrace{2\ldots2}_{k-1}\underbrace{1\ldots1}_{n-2k+1}.
\]

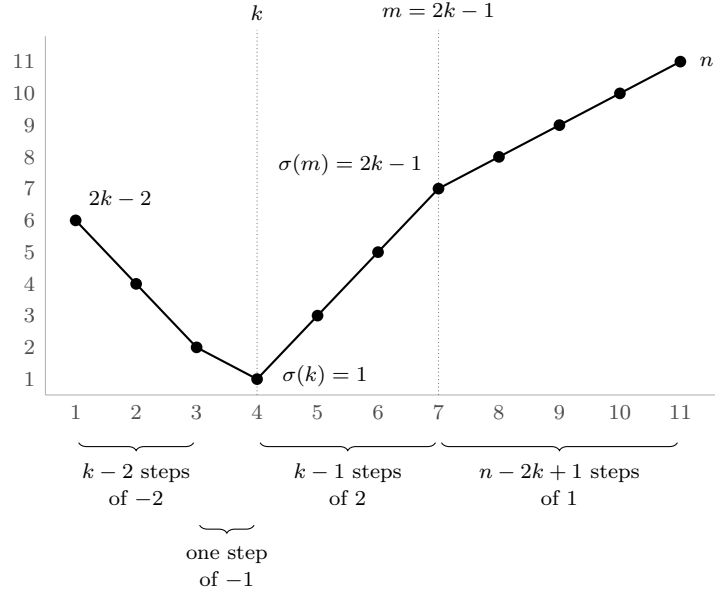
\begin{figure}[ht]
\centering
\begin{tikzpicture}[x=0.8cm,y=0.42cm]
  \draw[gray!60] (0.5,0.5) -- (11.6,0.5);
  \draw[gray!60] (0.5,0.5) -- (0.5,11.8);
  \foreach \i in {1,...,11} \node[below,font=\scriptsize,gray!60!black] at (\i,0.5) {\i};
  \foreach \v in {1,...,11} \node[left,font=\scriptsize,gray!60!black] at (0.5,\v) {\v};
  \draw[densely dotted,gray] (4,0.5) -- (4,12);
  \draw[densely dotted,gray] (7,0.5) -- (7,12);
  \node[font=\scriptsize,above] at (4,12) {$k$};
  \node[font=\scriptsize,above] at (7,12) {$m=2k-1$};
  \draw[thick] (1,6)--(2,4)--(3,2)--(4,1)--(5,3)--(6,5)--(7,7)--(8,8)--(9,9)--(10,10)--(11,11);
  \foreach \p/\v in {1/6,2/4,3/2,4/1,5/3,6/5,7/7,8/8,9/9,10/10,11/11} \fill (\p,\v) circle (2.2pt);
  \node[font=\scriptsize,anchor=south west] at (1.05,6.1) {$2k-2$};
  \node[font=\scriptsize,anchor=west] at (4.25,1.1) {$\sigma(k)=1$};
  \node[font=\scriptsize,anchor=south east] at (6.9,7.1) {$\sigma(m)=2k-1$};
  \node[font=\scriptsize,anchor=west] at (11.15,11) {$n$};
  \draw[decorate,decoration={brace,mirror,amplitude=3pt}] (1.05,-1.0) -- (2.95,-1.0)
       node[midway,below=4pt,font=\scriptsize,align=center] {$k-2$ steps\\ of $-2$};
  \draw[decorate,decoration={brace,mirror,amplitude=3pt}] (3.05,-3.6) -- (3.95,-3.6)
       node[midway,below=4pt,font=\scriptsize,align=center] {one step\\ of $-1$};
  \draw[decorate,decoration={brace,mirror,amplitude=3pt}] (4.05,-1.0) -- (6.95,-1.0)
       node[midway,below=4pt,font=\scriptsize,align=center] {$k-1$ steps\\ of $2$};
  \draw[decorate,decoration={brace,mirror,amplitude=3pt}] (7.05,-1.0) -- (10.95,-1.0)
       node[midway,below=4pt,font=\scriptsize,align=center] {$n-2k+1$ steps\\ of $1$};
\end{tikzpicture}
\caption{The permutation $\sigma=(2k-2)(2k-4)\ldots 2\,1\,3\,5\ldots(2k-1)(2k)(2k+1)\ldots n$ with $k=4$ and $n=11$, that is, $\sigma=6\;4\;2\;1\;3\;5\;7\;8\;9\;10\;11$.}
\label{F:m2km1}
\end{figure}

When $m=2k-2$, let $\sigma=(2k-3)(2k-5)\ldots 3\,1\,2\,4\ldots(2k-2)(2k-1)(2k)\ldots n$, see Figure \ref{F:m2km2}. Then $\sigma\in\Sym_n$, with $\sigma(k)=2$ and $\sigma(m)=2k-2$, and $\sigma$ has first derivative
\[
\underbrace{(-2)\ldots(-2)}_{k-2}\,1\,\underbrace{2\ldots2}_{k-2}\underbrace{1\ldots1}_{n-2k+2}.
\]

\begin{figure}[ht]
\centering
\begin{tikzpicture}[x=0.8cm,y=0.42cm]
  \draw[gray!60] (0.5,0.5) -- (11.6,0.5);
  \draw[gray!60] (0.5,0.5) -- (0.5,11.8);
  \foreach \i in {1,...,11} \node[below,font=\scriptsize,gray!60!black] at (\i,0.5) {\i};
  \foreach \v in {1,...,11} \node[left,font=\scriptsize,gray!60!black] at (0.5,\v) {\v};
  \draw[densely dotted,gray] (4,0.5) -- (4,12);
  \draw[densely dotted,gray] (6,0.5) -- (6,12);
  \node[font=\scriptsize,above] at (4,12) {$k$};
  \node[font=\scriptsize,above] at (6,12) {$m=2k-2$};
  \draw[thick] (1,5)--(2,3)--(3,1)--(4,2)--(5,4)--(6,6)--(7,7)--(8,8)--(9,9)--(10,10)--(11,11);
  \foreach \p/\v in {1/5,2/3,3/1,4/2,5/4,6/6,7/7,8/8,9/9,10/10,11/11} \fill (\p,\v) circle (2.2pt);
  \node[font=\scriptsize,anchor=south west] at (1.05,5.1) {$2k-3$};
  \node[font=\scriptsize,anchor=west] at (4.25,1.8) {$\sigma(k)=2$};
  \node[font=\scriptsize,anchor=north west] at (6.2,5.7) {$\sigma(m)=2k-2$};
  \node[font=\scriptsize,anchor=west] at (11.15,11) {$n$};
  \draw[decorate,decoration={brace,mirror,amplitude=3pt}] (1.05,-1.0) -- (2.95,-1.0)
       node[midway,below=4pt,font=\scriptsize,align=center] {$k-2$ steps\\ of $-2$};
  \draw[decorate,decoration={brace,mirror,amplitude=3pt}] (3.05,-3.6) -- (3.95,-3.6)
       node[midway,below=4pt,font=\scriptsize,align=center] {one step\\ of $1$};
  \draw[decorate,decoration={brace,mirror,amplitude=3pt}] (4.05,-1.0) -- (5.95,-1.0)
       node[midway,below=4pt,font=\scriptsize,align=center] {$k-2$ steps\\ of $2$};
  \draw[decorate,decoration={brace,mirror,amplitude=3pt}] (6.05,-1.0) -- (10.95,-1.0)
       node[midway,below=4pt,font=\scriptsize,align=center] {$n-2k+2$ steps\\ of $1$};
\end{tikzpicture}
\caption{$\sigma=(2k-3)(2k-5)\ldots 3\,1\,2\,4\ldots(2k-2)(2k-1)(2k)\ldots n$ with $k=4$ and $n=11$, that is, $\sigma=5\;3\;1\;2\;4\;6\;7\;8\;9\;10\;11$.}
\label{F:m2km2}
\end{figure}
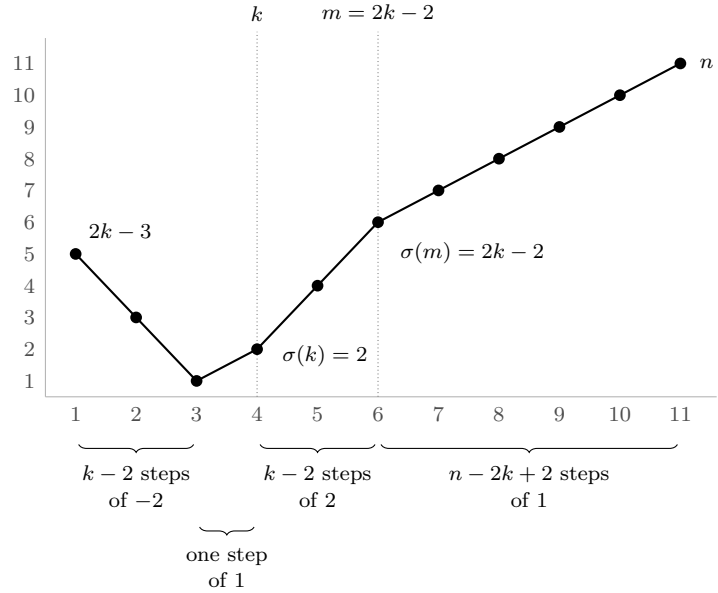

When $m=2k-3$ (so that $k\ge 4$, since $m>k$), let $\sigma=n\,(2k-5)(2k-7)\ldots 3\,1\,2\,4\ldots(2k-4)(2k-3)(2k-2)\ldots(n-1)$, see Figure \ref{F:m2km3}. Then $\sigma\in\Sym_n$, with $\sigma(k)=2$ and $\sigma(m)=2k-4$, and $\sigma$ has first derivative
\[
(2k-5-n)\underbrace{(-2)\ldots(-2)}_{k-3}\,1\,\underbrace{2\ldots2}_{k-3}\underbrace{1\ldots1}_{n-2k+3},
\]
where $2k-5-n\le -3$ since $n\ge m+1=2k-2$.

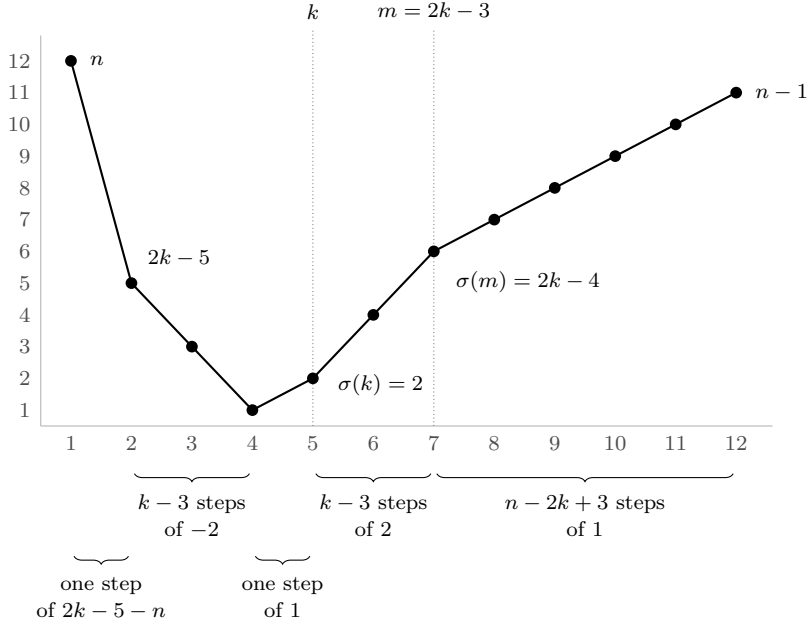
\begin{figure}[ht]
\centering
\begin{tikzpicture}[x=0.8cm,y=0.42cm]
  \draw[gray!60] (0.5,0.5) -- (12.6,0.5);
  \draw[gray!60] (0.5,0.5) -- (0.5,12.8);
  \foreach \i in {1,...,12} \node[below,font=\scriptsize,gray!60!black] at (\i,0.5) {\i};
  \foreach \v in {1,...,12} \node[left,font=\scriptsize,gray!60!black] at (0.5,\v) {\v};
  \draw[densely dotted,gray] (5,0.5) -- (5,13);
  \draw[densely dotted,gray] (7,0.5) -- (7,13);
  \node[font=\scriptsize,above] at (5,13) {$k$};
  \node[font=\scriptsize,above] at (7,13) {$m=2k-3$};
  \draw[thick] (1,12)--(2,5)--(3,3)--(4,1)--(5,2)--(6,4)--(7,6)--(8,7)--(9,8)--(10,9)--(11,10)--(12,11);
  \foreach \p/\v in {1/12,2/5,3/3,4/1,5/2,6/4,7/6,8/7,9/8,10/9,11/10,12/11} \fill (\p,\v) circle (2.2pt);
  \node[font=\scriptsize,anchor=west] at (1.15,12) {$n$};
  \node[font=\scriptsize,anchor=south west] at (2.1,5.2) {$2k-5$};
  \node[font=\scriptsize,anchor=west] at (5.25,1.8) {$\sigma(k)=2$};
  \node[font=\scriptsize,anchor=north west] at (7.2,5.7) {$\sigma(m)=2k-4$};
  \node[font=\scriptsize,anchor=west] at (12.15,11) {$n-1$};
  \draw[decorate,decoration={brace,mirror,amplitude=3pt}] (1.05,-3.6) -- (1.95,-3.6)
       node[midway,below=4pt,font=\scriptsize,align=center] {one step\\ of $2k-5-n$};
  \draw[decorate,decoration={brace,mirror,amplitude=3pt}] (2.05,-1.0) -- (3.95,-1.0)
       node[midway,below=4pt,font=\scriptsize,align=center] {$k-3$ steps\\ of $-2$};
  \draw[decorate,decoration={brace,mirror,amplitude=3pt}] (4.05,-3.6) -- (4.95,-3.6)
       node[midway,below=4pt,font=\scriptsize,align=center] {one step\\ of $1$};
  \draw[decorate,decoration={brace,mirror,amplitude=3pt}] (5.05,-1.0) -- (6.95,-1.0)
       node[midway,below=4pt,font=\scriptsize,align=center] {$k-3$ steps\\ of $2$};
  \draw[decorate,decoration={brace,mirror,amplitude=3pt}] (7.05,-1.0) -- (11.95,-1.0)
       node[midway,below=4pt,font=\scriptsize,align=center] {$n-2k+3$ steps\\ of $1$};
\end{tikzpicture}
\caption{The permutation $\sigma=n\,(2k-5)(2k-7)\ldots 3\,1\,2\,4\ldots(2k-4)(2k-3)(2k-2)\ldots(n-1)$ with $k=5$ and $n=12$, that is, $\sigma=12\;5\;3\;1\;2\;4\;6\;7\;8\;9\;10\;11$.}
\label{F:m2km3}
\end{figure}

\medskip
(b) Let $b=2k-m$ and $T=2m-k$. Since $m\le 2k-1$ and $n\ge 2m-k$, we have $b\ge 1$ and $n\ge T$. Notice also that $T=b+3\ell$. We define $\sigma$ by listing, one after the other, the entries of five blocks,
\[
\sigma=P\;Q\;R\;S\;U,
\]
which are defined as follows.
\begin{itemize}
\item $P=(T-2)(T-5)(T-8)\ldots(b+4)(b+1)$. It has $\ell$ entries, and it decreases by $3$.
\item $Q$ contains all the elements of $\{1,\ldots,b-1\}$, listed in two halves: first $(b-1)(b-3)(b-5)\ldots$, decreasing by $2$ until it reaches $1$ or $2$, then the remaining elements, increasing by $2$ until they reach $b-2$. It has $b-1$ entries (so $Q$ is empty if $b=1$, and $Q=1$ if $b=2$).
\item $R=b\,(b+3)(b+6)\ldots(b+3\ell)$. It has $\ell+1$ entries, it increases by $3$, and its last entry is $T$.
\item $S$ contains all the elements of $\{T+1,\ldots,n\}$, listed in two halves: first $(T+2)(T+4)(T+6)\ldots$, increasing by $2$ until it reaches $n-1$ or $n$; then the remaining elements, decreasing by $2$ until they reach $T+1$. It has $n-T$ entries (so $S$ is empty if $n=T$, and $S=T+1$ if $n=T+1$).
\item $U=(T-1)(T-4)(T-7)\ldots(b+5)(b+2)$. It has $\ell$ entries, and it decreases by $3$.
\end{itemize}
For example, if $(k,m,n)=(7,10,16)$, then $\ell=3$, $b=4$, $T=13$, and
\[
\sigma=\underbrace{11\;\,8\;\,5}_{P}\;\;\underbrace{3\;\,1\;\,2}_{Q}\;\;\underbrace{4\;\,7\;\,10\;\,13}_{R}\;\;\underbrace{15\;\,16\;\,14}_{S}\;\;\underbrace{12\;\,9\;\,6}_{U},
\]
with $\sigma'=(-3)(-3)(-2)(-2)\,1\,2\,3\,3\,3\,2\,1\,(-2)(-2)(-3)(-3)$ (see Figure \ref{F:PQRSU}).

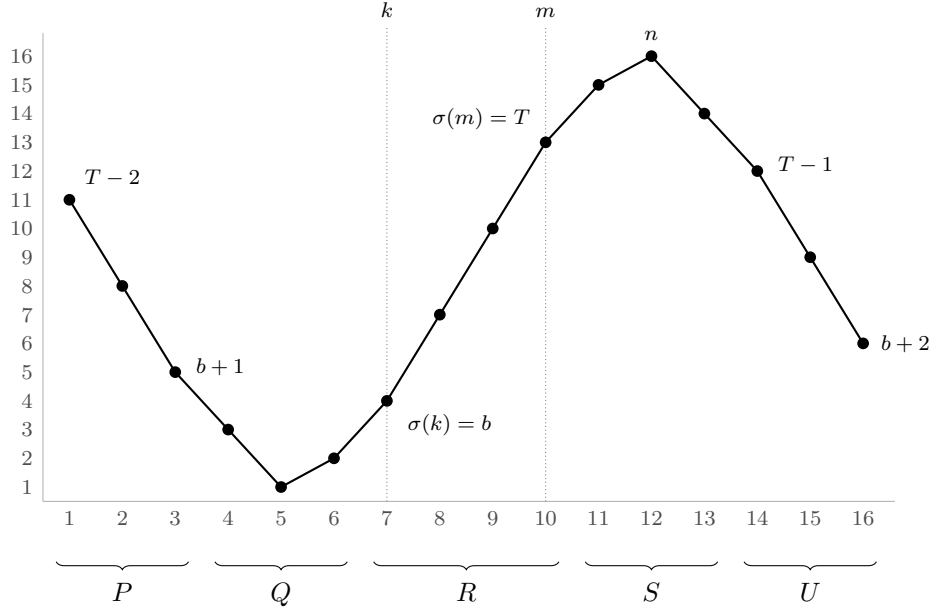
\begin{figure}[ht]
\centering
\begin{tikzpicture}[x=0.7cm,y=0.38cm]
  \draw[gray!60] (0.5,0.5) -- (16.6,0.5);
  \draw[gray!60] (0.5,0.5) -- (0.5,16.8);
  \foreach \i in {1,...,16} \node[below,font=\scriptsize,gray!60!black] at (\i,0.5) {\i};
  \foreach \v in {1,...,16} \node[left,font=\scriptsize,gray!60!black] at (0.5,\v) {\v};
  \draw[densely dotted,gray] (7,0.5) -- (7,17);
  \draw[densely dotted,gray] (10,0.5) -- (10,17);
  \node[font=\scriptsize,above] at (7,17) {$k$};
  \node[font=\scriptsize,above] at (10,17) {$m$};
  \draw[thick] (1,11)--(2,8)--(3,5)--(4,3)--(5,1)--(6,2)--(7,4)--(8,7)--(9,10)--(10,13)--(11,15)--(12,16)--(13,14)--(14,12)--(15,9)--(16,6);
  \foreach \p/\v in {1/11,2/8,3/5,4/3,5/1,6/2,7/4,8/7,9/10,10/13,11/15,12/16,13/14,14/12,15/9,16/6} \fill (\p,\v) circle (2.2pt);
  \node[font=\scriptsize,anchor=south west] at (1.1,11.1) {$T-2$};
  \node[font=\scriptsize,anchor=west] at (3.2,5.2) {$b+1$};
  \node[font=\scriptsize,anchor=north west] at (7.2,3.9) {$\sigma(k)=b$};
  \node[font=\scriptsize,anchor=south east] at (9.9,13.1) {$\sigma(m)=T$};
  \node[font=\scriptsize,anchor=south] at (12,16.2) {$n$};
  \node[font=\scriptsize,anchor=west] at (14.2,12.2) {$T-1$};
  \node[font=\scriptsize,anchor=west] at (16.15,6) {$b+2$};
  \foreach \l/\r/\n in {1/3/P, 4/6/Q, 7/10/R, 11/13/S, 14/16/U}
     \draw[decorate,decoration={brace,mirror,amplitude=3pt}] (\l-0.25,-1.6) -- (\r+0.25,-1.6)
          node[midway,below=4pt,font=\small] {$\n$};
\end{tikzpicture}
\caption{The permutation $\sigma=P\,Q\,R\,S\,U$ with $(k,m,n)=(7,10,16)$, so that $\sigma=11\;8\;5\;\,3\;1\;2\;\,4\;7\;10\;13\;\,15\;16\;14\;\,12\;9\;6$.}
\label{F:PQRSU}
\end{figure}

We first show that $\sigma\in\Sym_n$, and that $R$ occupies the positions $k,\ldots,m$. The entries of $Q$ are the elements of $\{1,\ldots,b-1\}$, and the entries of $S$ are the elements of $\{T+1,\ldots,n\}$. Among the elements of $\{b,\ldots,T\}$, those congruent to $b$, $b+1$ and $b+2$ modulo $3$ are exactly the entries of $R$, $P$ and $U$, respectively. Hence, every element of $[1,n]$ appears exactly once in $\sigma$. Moreover, $P$ and $Q$ together have $\ell+(b-1)=k-1$ entries, so the $\ell+1$ entries of $R$ occupy the positions $k,\ldots,k+\ell=m$. In particular, $\sigma(k)=b$, $\sigma(m)=T$, and $\sigma'(j)=3$ for $k\le j\le m-1$. Finally, $m\le n-1$, since $U$ is not empty.

We now compute $\sigma'(j)=\sigma(j+1)-\sigma(j)$ on $[1,k-1]$. If $\sigma(j)$ and $\sigma(j+1)$ are both in $P$, then $\sigma'(j)=-3$, if they are both in the first half of $Q$, then $\sigma'(j)=-2$, and if they are both in the second half of $Q$, then $\sigma'(j)=2$. It remains to compute $\sigma'(j)$ when $\sigma(j)$ is the last entry of a block (or of the first half of $Q$) and $\sigma(j+1)$ is the first entry of the next one.
\begin{itemize}
\item If $b\ge 2$, the last entry of $P$ is $b+1$ and the first entry of $Q$ is $b-1$, so here $\sigma'(j)=-2$.
\item If $b\ge 3$, the first half of $Q$ ends at $1$ and the second half starts at $2$ (when $b$ is even), or the first half ends at $2$ and the second half starts at $1$ (when $b$ is odd). So here $\sigma'(j)=(-1)^b$.
\item If $b\ge 3$, the last entry of $Q$ is $b-2$ and the first entry of $R$ is $b$, so here $\sigma'(j)=2$. If $b=2$, then $Q=1$ and the first entry of $R$ is $2$, so here $\sigma'(j)=1=(-1)^b$. If $b=1$, then $Q$ is empty, the last entry of $P$ is $2$ and the first entry of $R$ is $1$, so here $\sigma'(j)=-1=(-1)^b$.
\end{itemize}
Counting the entries of each kind, we conclude that, in all cases, $\sigma'$ restricted to $[1,k-1]$ is
\begin{equation}\label{prefixPQ}
\underbrace{(-3)\ldots(-3)}_{\ell-1}\;\underbrace{(-2)\ldots(-2)}_{\lceil (b-1)/2\rceil}\;(-1)^b\;\underbrace{2\ldots2}_{\lfloor (b-1)/2\rfloor}.
\end{equation}

We compute $\sigma'$ on $[m,n-1]$ in the same way. If $\sigma(j)$ and $\sigma(j+1)$ are both in the first half of $S$, then $\sigma'(j)=2$, if they are both in the second half of $S$, then $\sigma'(j)=-2$, and if they are both in $U$, then $\sigma'(j)=-3$. For the remaining positions:
\begin{itemize}
\item If $n\ge T+2$, the last entry of $R$ is $T$ and the first entry of $S$ is $T+2$, so here $\sigma'(j)=2$.
\item If $n\ge T+2$, the first half of $S$ ends at $n-1$ and the second half starts at $n$ (when $n-T$ is odd), or the first half ends at $n$ and the second half starts at $n-1$ (when $n-T$ is even). So here $\sigma'(j)=(-1)^{n-T+1}$.
\item If $n\ge T+2$, the last entry of $S$ is $T+1$ and the first entry of $U$ is $T-1$, so here $\sigma'(j)=-2$. If $n=T+1$, then $S=T+1$, so $\sigma'(m)=1=(-1)^{n-T+1}$ and $\sigma'(m+1)=-2$. If $n=T$, then $S$ is empty, the last entry of $R$ is $T$ and the first entry of $U$ is $T-1$, so $\sigma'(m)=-1=(-1)^{n-T+1}$.
\end{itemize}
Therefore, in all cases, $\sigma'$ restricted to $[m,n-1]$ is
\begin{equation}\label{suffixSU}
\underbrace{2\ldots2}_{\lfloor (n-T)/2\rfloor}\;(-1)^{n-T+1}\;\underbrace{(-2)\ldots(-2)}_{\lceil (n-T)/2\rceil}\;\underbrace{(-3)\ldots(-3)}_{\ell-1}.
\end{equation}

Putting everything together, 
\[
\sigma'=\underbrace{(-3)\ldots(-3)}_{\ell-1}\underbrace{(-2)\ldots(-2)}_{\lceil (b-1)/2\rceil}\,(\pm1)\,\underbrace{2\ldots2}_{\lfloor (b-1)/2\rfloor}\;\underbrace{3\ldots3}_{\ell}\;\underbrace{2\ldots2}_{\lfloor (n-T)/2\rfloor}\,(\pm1)\,\underbrace{(-2)\ldots(-2)}_{\lceil (n-T)/2\rceil}\underbrace{(-3)\ldots(-3)}_{\ell-1}.
\]
By \eqref{prefixPQ}, $\sigma'$ is non-decreasing on $[1,k-1]$ and $\sigma'(k-1)\le 2<3$, so \eqref{mono1} holds with $d=3$. By \eqref{suffixSU}, $\sigma'$ is non-increasing on $[m,n-1]$ and $\sigma'(m)\le 2<3$, so \eqref{mono2} holds. Since $\sigma'(j)=3$ for $k\le j\le m-1$, Lemma \ref{L:itisinflection} gives $\mathcal{I}(\sigma)=\{(k,m)\}$.

\medskip
(c) Since $m\ge k+2$, we have $\ell\ge 2$ and $n=2m-k-1=k+2\ell-1\ge m+1$. Let 
\[
\sigma=k\,1\,2\ldots(k-1)(k+1)(k+3)\ldots(k+2\ell-1)(k+2\ell-2)(k+2\ell-4)\ldots(k+2).
\]

See Figure \ref{F:partc}. 
\begin{figure}[ht]
\centering
\begin{tikzpicture}[x=0.8cm,y=0.42cm]
  \draw[gray!60] (0.5,0.5) -- (11.6,0.5);
  \draw[gray!60] (0.5,0.5) -- (0.5,11.8);
  \foreach \i in {1,...,11} \node[below,font=\scriptsize,gray!60!black] at (\i,0.5) {\i};
  \foreach \v in {1,...,11} \node[left,font=\scriptsize,gray!60!black] at (0.5,\v) {\v};
  \draw[densely dotted,gray] (4,0.5) -- (4,12.6);
  \draw[densely dotted,gray] (8,0.5) -- (8,12.6);
  \node[font=\scriptsize,above] at (4,12.6) {$k$};
  \node[font=\scriptsize,above] at (8,12.6) {$m$};
  \draw[thick] (1,4)--(2,1)--(3,2)--(4,3)--(5,5)--(6,7)--(7,9)--(8,11)--(9,10)--(10,8)--(11,6);
  \foreach \p/\v in {1/4,2/1,3/2,4/3,5/5,6/7,7/9,8/11,9/10,10/8,11/6} \fill (\p,\v) circle (2.2pt);
  \node[font=\scriptsize,anchor=south west] at (1.05,4.2) {$\sigma(1)=k$};
  \node[font=\scriptsize,anchor=north west] at (4.2,2.9) {$\sigma(k)=k-1$};
  \node[font=\scriptsize,anchor=south east] at (7.85,11.1) {$\sigma(m)=n$};
  \node[font=\scriptsize,anchor=west] at (9.2,10.2) {$k+2\ell-2$};
  \node[font=\scriptsize,anchor=west] at (11.15,6) {$k+2$};
  \draw[decorate,decoration={brace,mirror,amplitude=3pt}] (1.05,-3.6) -- (1.95,-3.6)
       node[midway,below=4pt,font=\scriptsize,align=center] {one step\\ of $1-k$};
  \draw[decorate,decoration={brace,mirror,amplitude=3pt}] (2.05,-1.0) -- (3.95,-1.0)
       node[midway,below=4pt,font=\scriptsize,align=center] {$k-2$ steps\\ of $1$};
  \draw[decorate,decoration={brace,mirror,amplitude=3pt}] (4.05,-1.0) -- (7.95,-1.0)
       node[midway,below=4pt,font=\scriptsize,align=center] {$\ell$ steps\\ of $2$};
  \draw[decorate,decoration={brace,mirror,amplitude=3pt}] (8.05,-3.6) -- (8.95,-3.6)
       node[midway,below=4pt,font=\scriptsize,align=center] {one step\\ of $-1$};
  \draw[decorate,decoration={brace,mirror,amplitude=3pt}] (9.05,-1.0) -- (10.95,-1.0)
       node[midway,below=4pt,font=\scriptsize,align=center] {$\ell-2$ steps\\ of $-2$};
\end{tikzpicture}
\caption{$\sigma=k\,1\,2\ldots(k-1)(k+1)(k+3)\ldots(k+2\ell-1)(k+2\ell-2)(k+2\ell-4)\ldots(k+2)$ with $k=4$ and $m=8$, so that $\ell=4$ and $n=11$, that is, $\sigma=4\;1\;2\;3\;5\;7\;9\;11\;10\;8\;6$.}
\label{F:partc}
\end{figure}
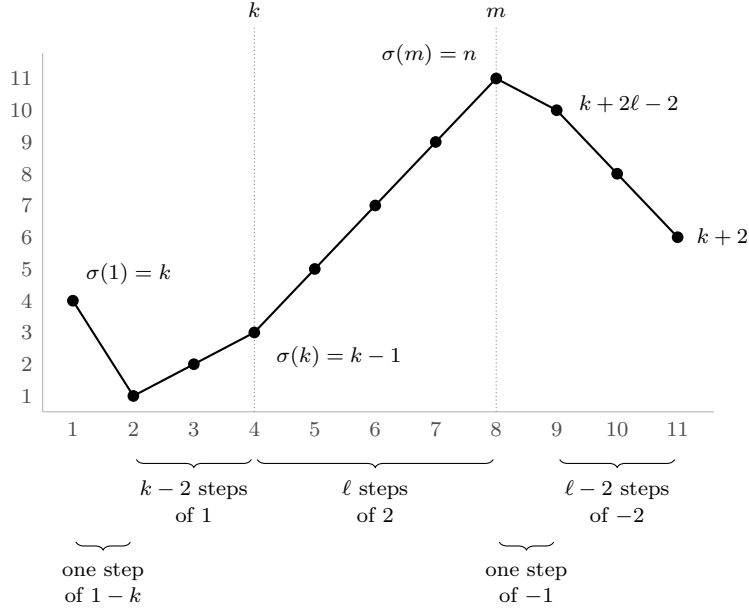

Then $\sigma\in\Sym_n$, with $\sigma(k)=k-1$ and $\sigma(m)=k+2\ell-1=n$, and 
\[
\sigma'=(1-k)\underbrace{1\ldots 1}_{k-2}\underbrace{2\ldots2}_{\ell}(-1)\underbrace{(-2)\ldots(-2)}_{\ell-2}.
\]
In particular, $\sigma'(j)=2$ for $k\le j\le m-1$, the first derivative is non-decreasing and at most $1$ on $[1,k-1]$, and non-increasing and at most $-1$ on $[m,n-1]$. Hence \eqref{mono1}--\eqref{mono2} hold with $d=2$, and Lemma \ref{L:itisinflection} gives $\mathcal{I}(\sigma)=\{(k,m)\}$.
\end{proof}

When $m=k+1$, Theorem \ref{T:witnesses}(b) gives a witness in $\Sym_n$ for every $n\ge m+1$. In this case there is also a simpler one: $\sigma=1\,2\ldots k\,n\,(n-1)(n-2)\ldots(k+1)\in\Sym_n$ has 
\[
\sigma'=\underbrace{1\ldots 1}_{k-1}(n-k)\underbrace{(-1)\ldots(-1)}_{n-k-1},
\]
so that $\mathcal{I}(\sigma)=\{(k,k+1)\}$ by Lemma \ref{L:itisinflection} (with $d=n-k\ge 2$). In particular, when $m=k+1$ the progression can have any $d\ge 2$, while $d\le 2$ when $m\ge 2k+1$ by Lemma \ref{L:d3}.

We are ready to fully characterize when $\{\sigma\in \Sym_n \mid \mathcal{I}(\sigma)=\{(k,m)\}\}\neq\emptyset$, in the case $k<m\le 2k$:

\begin{thm}\label{T:short}
Let $k\ge 2$, $k<m\le 2k$ and $n\ge m+1$. Then $|\{\sigma\in \Sym_n \mid \mathcal{I}(\sigma)=\{(k,m)\}\}|>0$ if and only if
\[
m\ge 2k-3\qquad\text{or}\qquad n\ge 2m-k-1 .
\]
\end{thm}

In other words, when $k<m\le 2k$ there is a witness in $\Sym_n$ for every large $n$, and small values of $n$ are excluded only when $m\le 2k-4$: in that case, the progression skips too few values to fill the positions before $k$, and we will show that at least $\ell-1$ entries of $\sigma$ must come after position $m$, that is, $n\ge m+\ell-1=2m-k-1$.

\begin{proof}
Write $\ell=m-k$.

\medskip
\noindent Sufficiency: we need to find $\sigma\in\Sym_n$ with $\mathcal{I}(\sigma)=\{(k,m)\}$ in each of the following cases.
\begin{itemize}
\item If $m\ge 2k-3$, this is Theorem \ref{T:witnesses}(a).
\item If $m\le 2k-4$ and $n\ge 2m-k$, this is Theorem \ref{T:witnesses}(b).
\item If $m\le 2k-4$ and $n=2m-k-1$, then $n\ge m+1$ gives $m\ge k+2$, and this is Theorem \ref{T:witnesses}(c).
\end{itemize}

\medskip
\noindent Necessity: assume $m\le 2k-4$, and let $\sigma\in\Sym_n$ with $\mathcal{I}(\sigma)=\{(k,m)\}$. We need to show that
\[
n\ge 2m-k-1 .
\]
Notice that $k<m\le 2k-4$ gives $k\ge 5$. By Lemma \ref{L:sign}, we may assume $d=\sigma(k+1)-\sigma(k)>0$, so that $\sigma$ is as in Lemma \ref{L:shape}. We divide the proof in cases, depending on $d$.
\begin{enumerate}
\item[\parbox{1.2cm}{Case \\ $d=1$:}] By Lemma \ref{L:d1}, this would give $k=2$, contradicting $k\ge 5$. So this case does not occur.

\item[\parbox{1.2cm}{Case \\ $d\ge 3$:}] By Lemma \ref{L:d3}, $n\ge 2m-k-1$.

\item[\parbox{1.2cm}{Case \\ $d=2$:}] Let $h_i$, $A$, $B$, $q$, $t$, $a$, $c$, $r$, $p$ be as in Lemma \ref{L:d2}. Since $m+\ell=2m-k$, \eqref{count2} gives
\[
n=m+p+\ell-q+c=2m-k+(p+c-q).
\]
Hence $n\ge 2m-k-1$ is equivalent to
\begin{equation}\label{goalq}
q\le p+c+1 .
\end{equation}
If $q\le 1$, then \eqref{goalq} holds, since $p\ge 0$ and $c\ge 0$. Assume from now on that $q\ge 2$, and assume by contradiction that \eqref{goalq} fails, that is,
\begin{equation}\label{contra}
q\ge p+c+2 .
\end{equation}
Since $A\subseteq[1,\ell]$, we have $q\le\ell=m-k\le k-4$. By \eqref{count1}, $r+a=k-1-q\ge 3$. Since $a\le 1$, we conclude that $r\ge 2$.
Now, since $q\ge 2$ and $r\ge 2$, Lemma \ref{L:d2}(d) gives $i_q=1$ and $i_{q-1}=3$. These are the two smallest elements of $A$, so $1\in A$, $3\in A$ and $2\notin A$.
Since $2\notin A$, we have $q\le\ell-1\le k-5$, so \eqref{count1} gives $r+a=k-1-q\ge 4$. On the other hand, $r\le 3$ by Lemma \ref{L:d2}(d), and $a\le 1$. Hence $r=3$, $a=1$ and $q=k-5=\ell-1$. Since $A$ has $\ell-1$ elements and $2\notin A$, we conclude that $A=[1,\ell]\setminus\{2\}$.

Since $3\in A$, we have $\ell\ge 3$. If $\ell\ge 4$, then $3,4\in A$ by what we have just shown, and Lemma \ref{L:d2}(a) (applied with $i=4$) gives $\{1,2,3,4\}\subseteq A$, contradicting the fact that $2\notin A$. Hence $\ell=3$, $A=\{1,3\}$ and $q=2$.
Finally, since $a=1$ and $A=\{1,3\}$, \eqref{order} gives
\[
\sigma(1)=n,\qquad\sigma(2)=h_3,\qquad\sigma(3)=h_1,
\]
and $k^*-1=q+a=3$. By \eqref{dist} (with $\ell=3$ and $a=1$),
\[
|\sigma'(1)|=n-h_3=p+2\qquad\text{and}\qquad|\sigma'(2)|=h_3-h_1=4 .
\]
Since $|\sigma'|$ is non-increasing on $[1,k^*-1]=[1,3]$, we have $p+2\ge 4$, that is, $p\ge 2$. However, \eqref{contra} and $q=2$ give $p+c\le 0$, that is, $p=0$. This is a contradiction, so \eqref{goalq} holds.
\end{enumerate}
In all cases, $n\ge 2m-k-1$.
\end{proof}

We summarize Theorems \ref{T:windows} and \ref{T:short} in the following corollary.

\begin{cor}\label{C:dichotomy}
Let $2\le k<m$ and $I=\{(k,m)\}$. Then $I$ is an admissible inflection set. Moreover, $I$ is vanishing if and only if $m\ge 2k+1$.
\end{cor}

\begin{proof}
If $m\ge 2k+1$, Theorem \ref{T:windows}(1) (with $N=1$) gives a witness in $\Sym_n$ for $n=2(m-k)+1$, and Corollary \ref{C:windows}(a) shows that there is no witness in $\Sym_n$ for $n\ge 2m-k+3$. Hence $I$ is admissible and vanishing. If $k<m\le 2k$, Theorem \ref{T:short} gives a witness in $\Sym_n$ for every $n\ge\max\{m+1,\,2m-k-1\}$. Hence $I$ is admissible and it is not vanishing.
\end{proof}

\section{Maximally inflected permutations}\label{S:max}

In this section, we study the size of $\mathcal M_n$, the set of maximally inflected permutations in $\Sym_n$. Our last main result, Theorem~\ref{T:estimatemax}, 
proves the explicit bound $|\mathcal M_n|/n!\le(2/3)^{\lfloor n/3\rfloor-1}$. To define a maximally inflected permutation, we start with a simple lemma:
\begin{lem}\label{L:maxnumber}
    If $\sigma\in \Sym_n$, then $|\mathcal{I}(\sigma)|\le n-3.$
\end{lem}

\begin{proof}
If $(i,j)\in\mathcal{I}(\sigma)$, then $j$ is the smallest index larger than $i$ with $\sigma''(j)\neq0$, so $j$ is determined by $i$. Hence $(i,j)\mapsto i$ is injective on $\mathcal{I}(\sigma)$, with image contained in $[2,n-2]$, and $|\mathcal{I}(\sigma)|\le n-3$.
\end{proof}

\begin{defn}
    For $n\geq 4$, if $\sigma\in\Sym_n$ and $|\mathcal{I}(\sigma)|=n-3$, we say $\sigma$ is \emph{maximally inflected}. We will write $\mathcal{M}_n=\{\sigma\in \Sym_n \ : \ \sigma \text{ is maximally inflected}\}$.
\end{defn}

\begin{defn}
A finite list $t_1,\dots,t_N$ is called \emph{zigzag} if
$(t_{j+1}-t_j)(t_{j+2}-t_{j+1})<0$ for $1\le j\le N-2$, that is, it is strictly
alternately increasing and decreasing.
\end{defn}

\begin{thm}\label{T:char}
For $\sigma\in\Sym_n$ the following are equivalent:
\begin{enumerate}
\item[(a)] $\sigma\in\mathcal{M}_n$.
\item[(b)] $\mathcal I(\sigma)=\{(2,3),(3,4),\dots,(n-2,n-1)\}$.
\item[(c)] $\sigma''(\ell)\neq0$ for all $\ell$ and $\sigma''(\ell)\sigma''(\ell+1)<0$
for $2\le\ell\le n-2$.
\item[(d)] The first derivative $D\sigma=(s_1,\dots,s_{n-1})$ is zigzag.
\end{enumerate}
\end{thm}

\begin{proof}
The equivalence between (a), (b) and (c) is immediate in view of the definition of inflection and Lemma \ref{L:maxnumber}.
To prove that (c) is equivalent to (d), notice that $\sigma''(\ell)=s_\ell-s_{\ell-1}$, hence (c) is exactly (d).
\end{proof}

\begin{cor}\label{C:alt}
Every alternating permutation is maximally inflected, so
\[
|\mathcal{M}_n| \ge\ 2E_n,
\]
where $E_n$ is the $n$-th Euler zigzag number.
\end{cor}

\begin{proof}
If $\sigma$ is alternating, then the signs of $s_1,\dots,s_{n-1}$ alternate and 
$D\sigma$ is zigzag. By Theorem \ref{T:char}(d), the result follows.
\end{proof}

Using code, we exhaustively enumerated the maximally inflected permutations for small $n$:
    \[
    \begin{tabular}{c|cccccccc}
      $n$  & 4&5&6&7&8&9&10 & 11 \\ \hline
       $\left|\mathcal{M}_n\right|$ & 14 & 44&190&904&5028& 30884 & 214598 & 1613896
    \end{tabular}\]
While a closed formula for $|\mathcal{M}_n|$ for general $n$ is still an open question, one can easily obtain initial bounds:
    
\begin{rem}
 We have
 \[
4\left(\frac{2}{\pi}\right)^{n+1}n!\lesssim |\mathcal{M}_n|\lesssim e^{-1/2}n!.
\]
In fact, one has $E_n\sim  2\left(\frac{2}{\pi}\right)^{n+1}n!$ (see, for example, \cite{stanley2010survey}), hence Corollary \ref{C:alt} gives the lower bound. For the upper bound, let \[
   A(n)=\{\sigma\in \Sym_n \ : \ \sigma \text{ contains three consecutive terms in arithmetic progression}\}.
   \]
   By the method of moments (see \cite{alon-spencer}, and also the discussion in \cite{mse5132077}), $|A(n)|\sim (1-e^{-1/2})n!$ as $n\rightarrow \infty$. Since $\mathcal{M}_n\subseteq \Sym_n\setminus A(n)$, we must have
    $|\mathcal{M}_n|\lesssim n!- \left(1-e^{-1/2}\right)n!,$ hence
\[
4\left(\frac{2}{\pi}\right)^{n+1}n!\lesssim |\mathcal{M}_n|\lesssim e^{-1/2}n!
\]
In particular, from the lower bound and Stirling's formula, $\lim\limits_{n\rightarrow\infty}|\mathcal{M}_n|^{1/n}=\infty$. 
\end{rem}

To obtain more precise bounds for $|\mathcal{M}_n|$ we first introduce the \emph{descent word of a permutation} (see \cite{stanley2010survey}, and \cite{bender-helton-richmond}, where descent words and their run words are used to count permutations with a prescribed pattern of rises and falls).

\begin{defn}\label{d:descent}
    Given $\sigma\in \Sym_n$, define the \emph{descent word} of $\sigma$ to be $w(\sigma)=w_1w_2\ldots w_{n-1}\in \{u,d\}^{n-1}$ given by
    \[
    w_j=\begin{cases} u \text{ if } \sigma(j)<\sigma(j+1)\\
    d, \text{ if } \sigma(j)>\sigma(j+1).
    \end{cases}
    \]
    Equivalently, $w_j=u$ if and only if $s_j>0$, where $D\sigma=(s_1,\ldots,s_{n-1})$ and $s_j=\sigma(j+1)-\sigma(j)$.
\end{defn}

\begin{defn} A \emph{run} of $w(\sigma)$ is a maximal block of consecutive letters, that is, a set of consecutive indices $R=\{\ell,\ell+1,\ldots, m\}\subseteq\{1,\ldots,n-1\}$ with $w_{\ell}=w_{\ell+1}=\ldots=w_m$ and $w_{\ell-1}\neq w_{\ell}$ (or $\ell=1$), and $w_{m+1}\neq w_m$ (or $m=n-1$). We define the \emph{length of the run} to be $|R|=m-\ell+1$. 

Notice that the runs partition $\{1,\ldots, n-1\}$, so we can list them in order as $R_1,\ldots, R_q$. We say that a run $R_j$ is \emph{interior} if $2\le j\le q-1$.    
\end{defn}

\begin{defn}
    Given $n\ge 1$, define
    \[
    \mathcal{O}_n=\{\sigma\in \Sym_n \ : \ \text{ every interior run of } w(\sigma) \text{ has odd length}\}.
    \]
    Let $O(n)=|\mathcal{O}_n|$ and $p_n=\frac{O(n)}{n!}$.
\end{defn}

\begin{thm}\label{T:odd} Given $n\ge 1$, we have $\mathcal{M}_n\subseteq \mathcal{O}_n$.
\end{thm}

\begin{proof}
   Let $\sigma\in\mathcal{M}_n$, and let $R=\{i,\ldots, j\}$ be an interior run of $w(\sigma)$ (so $2\le i\le j\le n-2$, $w_{i-1}\neq w_i$ and $w_{j+1}\neq w_j$). Then $\sigma''(i)$ and $\sigma''(j+1)$ must have opposite signs. Since $\sigma$ is maximally inflected, the sign of $\sigma''$ must alternate on $\{i,i+1,\ldots, j, j+1\}$, hence $\text{sign}(\sigma''(j+1))=(-1)^{j+1-i}\text{sign}(\sigma''(i))$. This implies $(-1)^{j+1-i}=-1$, that is, the length of $R$ is odd.
\end{proof}

We now follow \cite{bona2022combinatorics} and \cite{kitaev2011patterns}  to define the pattern of a permutation:
\begin{defn} Given a list $x=(r_1,\ldots, r_n)$ of $n$ distinct real numbers, its \emph{pattern/standardisation} is the permutation $st(x)\in \Sym_n$ defined by
\[
\text{st}(x)(i)=|\{j\in [n] \ : \ r_j\le r_i\}|.
\]
Given $\sigma\in \Sym_n$ and $L=\{p,p+1,\ldots,q\}$, we write
\[
\pi_L(\sigma)=\text{st}(\sigma(p),\sigma(p+1),\ldots,\sigma(q)).
\]
\end{defn}

\begin{ex} To fix ideas, let $n=5$, $\sigma=35412$, and $L=\{3,4,5\}$. In this case, $\pi_L(\sigma)=\text{st}(4,1,2)=312$.
\end{ex}
\begin{thm}\label{T:firstboundsMgood}
    For every $k,n \ge 2$, we have
\begin{equation}\label{Mnbound1}
|\mathcal{M}_n|\le (n!) p_k^{\lfloor n/k\rfloor}.
\end{equation}
    Additionally, $p_{k+\ell}\le p_kp_{\ell}$, hence 
    \begin{equation}\label{Mnbound2}
\limsup_{n\rightarrow\infty}\left(\frac{|\mathcal{M}_n|}{n!}\right)^{1/n}\le \lim_{n\rightarrow\infty} p_n^{1/n}.
    \end{equation}
\end{thm}

\begin{proof} Fix $k, n\ge 2$. Let $\sigma\in\mathcal{M}_n$ and $L=\{p,p+1,\ldots, p+k-1\}$ be $k$ consecutive natural numbers. For simplicity of notation, write $\pi_L(\sigma)=\pi$, that is, $\pi\in \Sym_k$ and $\text{st}(\sigma(p),\sigma(p+1),\ldots,\sigma(p+k-1))=\pi$. Following Definition \ref{d:descent}, write $w(\pi)=(w_1,\ldots,w_{k-1})$.
    We have $w(\pi)_l=w(\sigma)_{p+l-1}.$
  
    Let $R=\{i,\ldots, j\}$ be an interior run of $w(\pi)$ (so $2\le i$ and $j+1\le k-1$). Then $w_{p+(i-1)-1}$ and $w_{p+(j+1)-1}$ are letters of $w(\sigma)$ which differ from $w_{p+i-1}$ and $w_{p+j-1}$. Therefore $\tilde{R}=\{p+i-1,\ldots,p+j-1\}$ is an interior run of $w(\sigma)$. By Theorem \ref{T:odd}, $|\tilde{R}|=|R|$ must be odd.
We conclude that if $\sigma\in\mathcal{M}_n$, then $\pi_L(\sigma)\in \mathcal{O}_k$ for every set $L$ of $k$ consecutive natural numbers.

We now divide $\{1,\ldots,n\}$ into $m$ consecutive blocks $L_1,\ldots,L_m$ of size $k$, where $m=\lfloor n/k\rfloor$, together (potentially) with a last block $T$ of size less than $k$. Sampling $\sigma$ uniformly from $\Sym_n$ is equivalent to sampling the ordered partition of values $(S_1,\ldots,S_m,S_T)$ (where $S_j=\sigma(L_j)$ for $1\le j\le m$) uniformly among all such partitions, and independently, choosing a uniform ordering of each $S_i$ and $S_T.$ Given $1\le i\le m$, determining whether $\pi_{L_i}(\sigma)\in\mathcal{O}_k$ or not depends only on the ordering chosen for $S_i$. We have
\[
\Pr[\pi_{L_i}(\sigma)\in\mathcal{O}_k \ | \ S_1,\ldots,S_m,S_T]=\frac{O(k)}{k!}=p_k.
\]
Since the events are conditionally independent, we obtain \eqref{Mnbound1}:
\[
\frac{|\mathcal{M}_n|}{n!}\le \Pr[\cap_{j=1}^m\{\pi_{L_j}(\sigma)\in \mathcal{O}_k\}]=p_k^m.
\]

To prove \eqref{Mnbound2}, we start by showing $p_{k+\ell}\le p_k p_{\ell}$ for $k,\ell\ge 1$. 
Let $\zeta\in\mathcal{O}_{k+\ell}$. Divide $\{1,\ldots,k+\ell\}$ into $L_1=\{1,\ldots,k\}$ and $L_2=\{k+1,\ldots,k+\ell\}$. The previous arguments show $\pi_{L_1}(\zeta)\in\mathcal{O}_k$ and $\pi_{L_2}(\zeta)\in\mathcal{O}_{\ell}$. Sampling $\zeta$ uniformly from $\Sym_{k+\ell}$ gives
\[
p_{k+\ell}=\Pr[\zeta \in \mathcal{O}_{k+\ell}]
\le \Pr[\pi_{L_1}(\zeta)\in\mathcal{O}_k \text{ and } \pi_{L_2}(\zeta)\in \mathcal{O}_\ell]=p_kp_{\ell}.
\]

Applying Fekete's Lemma (\cite{steele1997}) with $a_n=\log p_n$ allows us to conclude that $\lim\limits_{n\rightarrow\infty}\frac{\log(p_n)}{n}=\inf\frac{\log(p_n)}{n}$, hence $\lim\limits_{n\rightarrow\infty}p_n^{1/n}$ exists and
\begin{equation}\label{inf}
\lim\limits_{n\rightarrow\infty}p_n^{\frac{1}{n}}=\inf p_n^{\frac{1}{n}}.
\end{equation}
    By Theorem \ref{T:odd}, $\mathcal{M}_n\subseteq\mathcal{O}_n$, hence
\[
\left(\frac{|\mathcal{M}_n|}{n!}\right)^{1/n}\le p_n^{1/n}.
\]
Together with \eqref{inf}, this gives \eqref{Mnbound2}.

\end{proof}

\begin{lem}\label{L:triple} Let $S\subset\mathbb{R}$ be finite, with $|S|\ge 3$. Let $(x,y,z)$ be a uniformly random ordered triple of distinct elements of $S$. Fix $s\in\{-1,0,1\}$. Then
\[
\Pr[\text{sign}(x-2y+z)=s]\le \frac{2}{3}.
\]

\end{lem}

\begin{proof}
    Given $x,y,z\in S$, each of the six possible orderings between them is equally likely. There are two orderings where $y$ is the smallest, and in that case, $x-2y+z>0$. There are two orderings where $y$ is the largest, in which case $x-2y+z<0$. For the remaining two orderings, where $y$ is the intermediate value, the sign of $x+z-2y$ is the same. Therefore, for fixed $s\in \{-1,1\}$ at most 4 orderings realize it, and only the two orderings where $y$ is the intermediate value can give $x-2y+z=0$. Hence $\Pr[\text{sign}(x-2y+z)=s]\le 4/6=2/3$ for every $s\in\{-1,0,1\}$.
\end{proof}

We are now ready to prove the main result of this section:
\begin{thm}\label{T:estimatemax} Given $n\ge 3$,
\begin{equation}\label{Mbound3}
\frac{|\mathcal{M}_n|}{n!}\le\left(\frac{2}{3}\right)^{\lfloor n/3\rfloor-1}.
\end{equation}
In particular, 
\[
\lim\limits_{n\rightarrow\infty}\frac{|\mathcal{M}_n|}{n!}=0\]
and 
\[
\limsup_{n\rightarrow\infty}\left(\frac{|\mathcal{M}_n|}{n!}\right)^{1/n}\le (2/3)^{1/3}.
\]
    
\end{thm}

\begin{proof}
  Fix $n\ge 3$. By Theorem \ref{T:char}, $\sigma\in\mathcal{M}_n$ if and only if there exists $a\in\{-1,1\}$ so that \begin{equation}\label{secondder}
\text{sign}(\sigma''(i))=a(-1)^i \text{ for } 2\le i\le n-1.
\end{equation}
Call $k=\lfloor n/3\rfloor$. If $k=1$, then $n\in\{3,4,5\}$ and it is trivially true that $|\mathcal{M}_n|\le n!$, so we can assume $k\ge 2$.

    Given $1\le j\le k$, let $L_j=\{3j-2,3j-1,3j\}$. Notice that the blocks $L_1,\ldots, L_k$ are disjoint. Let $\sigma\in\mathcal{M}_n$ and $Q_j=\sigma''(3j-1)$. Condition \eqref{secondder} gives
    $\text{sign}(Q_j)=a(-1)^{3j-1}=a(-1)^{j-1}=(-1)^{j-1}\text{sign}(Q_1)$ for $1\le j\le k$.

   Calling $E_j=\{\sigma\in \Sym_n \ : \ \text{sign}(Q_j)=(-1)^{j-1}\text{sign}(Q_1)\}$, we have proven that
   \[
   \mathcal{M}_n\subseteq \bigcap_{j=2}^k E_j.
   \]

   Let $\sigma$ be a uniformly random permutation in $\Sym_n$. Given $1\le j\le k$, let $\mathcal{F}_j$ be the $\sigma$-algebra generated by the random vector $(\sigma(1),\ldots,\sigma(3j))$ and 
     $T_j=[n]\setminus\{\sigma(1),\ldots,\sigma(3j-3)\}$, so $|T_j|\ge n-3j+3\ge n-3k+3\ge 3$. Notice that conditioned on $\mathcal F_{j-1}$ (for $2\le j\le k$), the ordered triple $(\sigma(3j-2),\sigma(3j-1),\sigma(3j))$ is a uniformly random ordered triple of distinct elements of $T_j$.  Moreover, $\text{sign}(Q_1)$ is $\mathcal F_{1}$-measurable since  $Q_1=\sigma''(2)=\sigma(3)+\sigma(1)-2\sigma(2)$. Since $\mathcal{F}_1\subseteq\mathcal{F}_{j-1}$, $(-1)^{j-1}\text{sign}(Q_1)$ is $\mathcal{F}_{j-1}$-measurable. Hence, 
     given $2\le j\le k,$ Lemma \ref{L:triple} yields
   \begin{equation}\label{conseqlemma}
   \Pr[E_j \ | \ \mathcal{F}_{j-1}]\le \frac{2}{3}.
   \end{equation}
   Let $F_{\ell}=\cap_{j=2}^{\ell}E_{j}$. 
   Using the law of iterated expectations,
   \[
   \Pr[F_{\ell}]=\Pr[F_{\ell-1}\cap E_{\ell}]=\mathbb{E}[\mathbf 1_{F_{\ell-1}}\mathbf 1_{E_{\ell}}]=\mathbb{E}[\mathbb{E}[\mathbf 1_{F_{\ell-1}}\mathbf 1_{E_{\ell}} \ | \ \mathcal{F}_{\ell-1}]].
   \]
   Now, if $\mathcal{G}$ is a $\sigma$-algebra and $Y$ is $\mathcal G$-measurable, then $\mathbb E[YX\mid\mathcal G]=Y\,\mathbb E[X\mid\mathcal G]$.
   We apply this with $Y=\mathbf 1_{F_{\ell-1}}$ (which is $\mathcal F_{\ell-1}$-measurable since $F_{\ell-1}=\bigcap_{j=2}^{\ell-1}E_j$), and obtain

   \[
   \mathbb{E}[\mathbb{E}[\mathbf 1_{F_{\ell-1}}\mathbf 1_{E_{\ell}} \ | \ \mathcal{F}_{\ell-1}]]=  \mathbb{E}[\mathbf 1_{F_{\ell-1}}\mathbb{E}[\mathbf 1_{E_{\ell}} \ | \ \mathcal{F}_{\ell-1}]]= 
   \mathbb{E}[\mathbf 1_{F_{\ell-1}}\Pr[E_{\ell} \ | \ \mathcal{F}_{\ell-1}]].
   \]
By \eqref{conseqlemma}, $\Pr[E_\ell\mid\mathcal F_{\ell-1}]\le \frac{2}{3}$, hence by the monotonicity of expectation,
\[
\mathbb E\big[\mathbf 1_{F_{\ell-1}}\Pr[E_\ell\mid\mathcal F_{\ell-1}]\big]\le \mathbb E\big[\mathbf 1_{F_{\ell-1}}\cdot\tfrac{2}{3}\big]=\tfrac{2}{3}\,\mathbb E[\mathbf 1_{F_{\ell-1}}]=\tfrac{2}{3}\Pr[F_{\ell-1}].
\]

   Putting everything together yields 
   $\Pr[F_{\ell}]\le \frac{2}{3}\Pr[F_{\ell-1}]$. Since 
   \[
   \Pr[F_2]=\Pr[E_2]=
   \mathbb{E}[\mathbf 1_{E_2}]= \mathbb{E}\big[\,\mathbb{E}[\mathbf 1_{E_2}\mid\mathcal F_1]\,\big]= 
   \mathbb{E}\big[\Pr[E_2\mid\mathcal F_1]\big]\le\mathbb{E}\Big[\frac{2}{3}\Big]=\frac{2}{3},
   \]
   we conclude
   $\Pr[F_{\ell}]\le\left(\frac{2}{3}\right)^{\ell-1}$. Since $\mathcal{M}_n\subseteq\bigcap_{j=2}^kE_j$, we obtain \eqref{Mbound3}:
   \[
   \frac{|\mathcal{M}_n|}{n!}\le\left(\frac{2}{3}\right)^{k-1}.
   \]
   Finally, $k-1\ge n/3-2$, hence $(2/3)^{k-1}\le (9/4)(2/3)^{n/3}$, from which we conclude that 
   \[   \limsup_{n\rightarrow\infty}\left(\frac{|\mathcal{M}_n|}{n!}\right)^{1/n}\le\left(\frac{2}{3}\right)^{1/3}.
   \]
\end{proof}

\printbibliography

@article{brualdi-dahl,
  author  = {Brualdi, Richard A. and Dahl, Geir},
  title   = {Permutation matrices, their discrete derivatives and extremal properties},
  journal = {Vietnam J. Math.},
  volume  = {48},
  number  = {4},
  year    = {2020},
  pages   = {719--740},
  doi     = {10.1007/s10013-020-00392-5},
  note    = {arXiv:1908.03739}
}

@article{OmarTroyka2026,
  author  = {Omar, Mohamed and Troyka, Justin M.},
  title   = {Growth rates of permutations with given descent or peak set},
  journal = {European Journal of Combinatorics},
  volume  = {131},
  pages   = {104246},
  year    = {2026},
  doi     = {10.1016/j.ejc.2025.104246}
}

@article{billey-burdzy-sagan,
  author  = {Billey, Sara and Burdzy, Krzysztof and Sagan, Bruce E.},
  title   = {Permutations with given peak set},
  journal = {J. Integer Seq.},
  volume  = {16},
  year    = {2013},
  pages   = {Article 13.6.1}
}

@article{davis-nelson-petersen-tenner,
  author  = {Davis, Robert and Nelson, Sarah A. and Petersen, T. Kyle and Tenner, Bridget E.},
  title   = {The pinnacle set of a permutation},
  journal = {Discrete Math.},
  volume  = {341},
  number  = {11},
  year    = {2018},
  pages   = {3249--3270}
}

@article{diazlopez-pinnacle,
  author  = {Diaz-Lopez, Alexander and Harris, Pamela E. and Huang, Isabella and Insko, Erik and Nilsen, Lars},
  title   = {A formula for enumerating permutations with a fixed pinnacle set},
  journal = {Discrete Math.},
  volume  = {344},
  number  = {6},
  year    = {2021},
  pages   = {Paper No. 112375}
}

@article{elizalde-noy,
  author  = {Elizalde, Sergi and Noy, Marc},
  title   = {Consecutive patterns in permutations},
  journal = {Adv. in Appl. Math.},
  volume  = {30},
  year    = {2003},
  pages   = {110--123}
}

@article{ehrenborg-kitaev-perry,
  author  = {Ehrenborg, Richard and Kitaev, Sergey and Perry, Peter},
  title   = {A spectral approach to consecutive pattern-avoiding permutations},
  journal = {J. Comb.},
  volume  = {2},
  year    = {2011},
  pages   = {305--353}
}

@article{fewster-siemssen,
  author  = {Fewster, Christopher J. and Siemssen, Daniel},
  title   = {Enumerating permutations by their run structure},
  journal = {Electron. J. Combin.},
  volume  = {21},
  number  = {4},
  year    = {2014},
  pages   = {Paper 4.18}
}

@article {zhuang-runs,
    AUTHOR = {Zhuang, Yan},
     TITLE = {Counting permutations by runs},
   JOURNAL = {J. Combin. Theory Ser. A},
  FJOURNAL = {Journal of Combinatorial Theory. Series A},
    VOLUME = {142},
      YEAR = {2016},
     PAGES = {147--176}
}

@article{degs,
  author  = {Davis, James A. and Entringer, Roger C. and Graham, Ronald L. and Simmons, Gustavus J.},
  title   = {On permutations containing no long arithmetic progressions},
  journal = {Acta Arith.},
  volume  = {34},
  number  = {1},
  year    = {1977},
  pages   = {81--90}
}

@misc{oeis-A295370,
  author       = {{OEIS Foundation Inc.}},
  title        = {Number of permutations of $[n]$ avoiding three consecutive terms in arithmetic progression, {E}ntry {A295370} in {T}he {O}n-{L}ine {E}ncyclopedia of {I}nteger {S}equences},
  year         = {2026},
  howpublished = {\url{https://oeis.org/A295370}}
}

@book{alon-spencer,
  author    = {Alon, Noga and Spencer, Joel H.},
  title     = {The Probabilistic Method},
  edition   = {4th},
  series    = {Wiley Series in Discrete Mathematics and Optimization},
  publisher = {John Wiley \& Sons},
  address   = {Hoboken, NJ},
  year      = {2016},
  isbn      = {978-1-119-06195-3},
}

@misc{mse5132077,
  title        = {Asymptotic growth of {3AP} permutations ({A295370})},
  author       = {Garnier (https://stackexchange.com)},
  howpublished = {Mathematics Stack Exchange},
  year         = {2026},
  note         = {URL: \url{https://stackexchange.com}; accessed 23-September-2026}
}

@book{steele1997,
  author    = {Steele, J. Michael},
  title     = {Probability Theory and Combinatorial Optimization},
  series    = {CBMS-NSF Regional Conference Series in Applied Mathematics},
  volume    = {69},
  publisher = {SIAM},
  address   = {Philadelphia, PA},
  year      = {1997},
}

@book{bona2022combinatorics,
  title={Combinatorics of Permutations},
  author={B{\'o}na, Mikl{\'o}s},
  year={2022},
  publisher={CRC Press},
  edition={3rd}
}

@book{kitaev2011patterns,
  title={Patterns in Permutations and Words},
  author={Kitaev, Sergey},
  year={2011},
  publisher={Springer},
  address = {Heidelberg},
}

@incollection{stanley2010survey,
  author    = {Stanley, Richard P.},
  title     = {A Survey of Alternating Permutations},
  booktitle = {Combinatorics and Graphs},
  series    = {Contemporary Mathematics},
  volume    = {531},
  publisher = {American Mathematical Society},
  address   = {Providence, RI},
  year      = {2010},
  pages     = {165--196}
}

@article {PinnaclesTypeBD,
    AUTHOR = {Gonz\'{a}lez, Nicolle and Harris, Pamela E. and Rojas Kirby,
              Gordon and Smit Vega Garcia, Mariana and Tenner, Bridget
              Eileen},
     TITLE = {Pinnacle sets of signed permutations},
   JOURNAL = {Disc. Math.},
  FJOURNAL = {Discrete Mathematics},
    VOLUME = {346},
      YEAR = {2023},
    NUMBER = {7},
      ISSN = {0012-365X},
   MRCLASS = {Prelim},
  MRNUMBER = {4575065},
       DOI = {10.1016/j.disc.2023.113439},
       URL = {https://doi.org/10.1016/j.disc.2023.113439},
}

@article{bender-helton-richmond,
  author  = {Bender, Edward A. and Helton, William J. and Richmond, L. Bruce},
  title   = {Asymptotics of permutations with nearly periodic patterns of rises and falls},
  journal = {Electron. J. Combin.},
  volume  = {10},
  year    = {2003},
  pages   = {Research Paper 40, 27 pp.}
}

@article{kotzig-laufer,
  author  = {Kotzig, Anton and Laufer, P. J.},
  title   = {When are permutations additive?},
  journal = {Amer. Math. Monthly},
  volume  = {85},
  year    = {1978},
  pages   = {364--365}
}

@misc{oeis-A002047,
  author       = {{OEIS Foundation Inc.}},
  title        = {Number of $3\times(2n+1)$ zero-sum arrays with entries $-n,\dots,0,\dots,n$, {E}ntry {A002047} in {T}he {O}n-{L}ine {E}ncyclopedia of {I}nteger {S}equences},
  year         = {2026},
  howpublished = {\url{https://oeis.org/A002047}}
}

@article{gonzalez2023mesas,
  title={Mesas of Stirling permutations},
  author={Gonz{\'a}lez, Nicolle and Harris, Pamela E and Rojas Kirby, Gordon and Smit Vega Garcia, Mariana and Tenner, Bridget Eileen},
  journal={Journal of Combinatorics},
  volume={16},
  number={3},
  year={2025},
  eprint={2306.12158},
  archivePrefix={arXiv},
  primaryClass={math.CO},
  doi={10.4310/JOC.250830033838}
}

\end{document}